\documentclass[a4paper, 11pt, twoside, notitlepage]{amsart}

\usepackage{amsmath,amscd}
\usepackage{amssymb}
\usepackage{amsthm}
\usepackage{comment}
\usepackage{graphicx}
\usepackage{epstopdf}
\usepackage{mathrsfs}
\usepackage[ocgcolorlinks, linkcolor=blue]{hyperref}

\usepackage[ocgcolorlinks,linkcolor=blue]{hyperref}

\theoremstyle{plain}
\newtheorem{thm}{Theorem}[section]
\newtheorem{prop}{Proposition}[section]
\newtheorem{lem}[prop]{Lemma}

\newtheorem{rmk}[prop]{Remark}

\numberwithin{equation}{section}

\title[A Time-Harmonic Scattering Theory for Hyperbolic Space]{The Sommerfeld-Rellich Framework for Scattering on Hyperbolic Space: Far-Field Patterns and Inverse Problems}

\author{Lu Chen}
\address{Key Laboratory of Algebraic Lie Theory and Analysis of Ministry of Education, School of Mathematics and Statistics, Beijing Institute of Technology, Beijing 100081, P. R. China}
\email{chenlu5818804@163.com}
\author{Hongyu Liu}
\address{Department of Mathematics, City University of Hong Kong, Hong Kong SAR, China}
\email{hongyliu@cityu.edu.hk; hongyu.liuip@gmail.com}

\thanks{The first author was partly supported by the National Key Research and Development Program (No.
2022YFA1006900) and National Natural Science Foundation of China (No. 12271027), the second author was supported by NSFC/RGC Joint Research Scheme, N\_CityU101/21, ANR/RGC Joint Research Scheme, A-CityU203/19, and the Hong Kong RGC General Research Funds (projects 11311122, 12301420 and 11300821).}

\begin{document}
	
	\maketitle
	
	\begin{abstract}
		
    This paper establishes a complete time-harmonic scattering theory for hyperbolic space, formulating it within the classical Sommerfeld-Rellich paradigm centered on far-field patterns--a foundational framework that has been absent despite the well-developed spectral and time-dependent theories for this geometry. We explicitly construct the ingoing and outgoing fundamental solutions for the Helmholtz operator and perform a precise asymptotic analysis at the conformal boundary to derive a {hyperbolic Sommerfeld radiation condition}. This condition, which is a local criterion at infinity, uniquely selects physically admissible outgoing solutions. We prove a {hyperbolic Rellich theorem} guaranteeing the uniqueness of the scattered field and its far-field pattern from asymptotic data. Within this rigorous framework, we solve the direct scattering problem for compact sources, penetrable media, and impenetrable obstacles, providing explicit representations for the corresponding {far-field patterns}. As a principal application and demonstration of the framework's utility, we initiate the study of inverse scattering on hyperbolic space, formulating both the inverse obstacle and inverse medium problems where the objective is to reconstruct the scatterer from measurements of its far-field pattern. Our work lays the groundwork for a far-field-based approach to scattering and inversion in hyperbolic geometry. Due to the local feature of the {hyperbolic Sommerfeld radiation condition}, our study can be readily extended to the broader asymptotically hyperbolic manifolds.
		
		\medskip
		
		
		
	\end{abstract}

	\tableofcontents

	\section{Introduction}
	
\subsection{Background and motivation}

There is a long history of scattering theory on hyperbolic manifolds, originating from the observation that the Eisenstein series for a Fuchsian group serves as a generalized eigenfunction for the Laplacian on the associated quotient of hyperbolic space. The fundamental reference for the finite-volume case is the work of Lax and Phillips \cite{LP1}, while the infinite-volume case was initiated by Patterson \cite{Pat}. This has generated a wealth of results; we refer to \cite{Guillarmou} for a comprehensive bibliography and to \cite{Mel} for a review. Notably, scattering theory on hyperbolic space is deeply connected to the physical framework of Anti-de Sitter (AdS) theory, where hyperbolic geometry models the asymptotic structure of AdS spacetimes, encoding the interaction between bulk fields and the conformal boundary at infinity. In this setting, scattering and inverse scattering on hyperbolic space provide a natural language for describing how bulk wave propagation is reflected in boundary observables---a perspective that underlies many aspects of the AdS/CFT correspondence.
\medskip

Early foundational work by Friedlander \cite{Fre} developed a time-dependent scattering theory for wave propagation on hyperbolic space, introducing the concept of radiation fields to describe the asymptotic behavior of solutions at infinity. This framework was later placed in a broader microlocal setting by Mazzeo and Melrose \cite{Mel,MaMe}, who established a scattering theory on asymptotically hyperbolic manifolds via the \emph{resolvent} and the \emph{scattering matrix} at infinity---a formulation intrinsically linked to the time-dependent paradigm through the limiting absorption principle. Building on these developments, Joshi and S\'{a} Barreto \cite{JS} studied inverse scattering problems on hyperbolic spaces, showing how such scattering data can recover geometric information. More recently, Isozaki \cite{Iso1,Iso2} developed alternative time-dependent formulations, emphasizing connections with long-range scattering phenomena.
\medskip

This body of work has largely been developed within the {time-dependent or stationary spectral framework}, where the primary objects are wave operators, the scattering matrix, or the boundary resolvent. In contrast, the {time-harmonic framework}---centered on the explicit characterization of outgoing radiation conditions, far-field patterns, and their role as measurable data for inverse problems---has been the cornerstone of Euclidean scattering theory since the seminal work of Sommerfeld \cite{Som} and Rellich. Sommerfeld's radiation condition (1912) and Rellich's uniqueness theorem (1943) provided the rigorous foundation for the exterior Helmholtz equation. The systematic use of the far-field pattern as the central data for inverse scattering was later developed by Colton, Kress, and others \cite{CK}, creating a complete ``Sommerfeld-Rellich-Colton-Kress'' (SRCK) structure. This structure is not merely theoretical; it is the operational foundation of most numerical reconstruction algorithms in acoustics, electromagnetics, and medical imaging, where multi-static far-field data serves as the primary input.
\medskip

However, a parallel, fully-developed time-harmonic framework for hyperbolic space---complete with explicit outgoing Green's functions, a hyperbolic Sommerfeld radiation condition, a Rellich-type theorem, and a well-defined far-field pattern as the kernel of a far-field operator---remains absent. Existing theories, while profound, do not provide the concrete, frequency-domain tools necessary to pose and solve inverse problems in the manner standard in Euclidean applications.
\medskip

The main purpose of this paper is to establish this missing foundation. We develop a comprehensive time-harmonic scattering theory on hyperbolic space, explicitly constructing the SRCK structure for this geometry. Our specific contributions are:
\begin{enumerate}
    \item We derive the explicit forms of the \emph{ingoing} and \emph{outgoing} Green's functions for the Helmholtz operator on hyperbolic space and verify their fundamental properties.
    \item We perform a precise asymptotic analysis of these Green's functions at the conformal boundary (infinity) and use this to formulate {hyperbolic Sommerfeld radiation conditions} that uniquely select physically relevant outgoing solutions.
    \item We prove a {hyperbolic Rellich theorem}, ensuring the uniqueness of the scattered field and its far-field pattern from the asymptotic data.
    \item Within this rigorous framework, we formulate and solve the direct scattering problem for three canonical cases: scattering by a compact source, by a penetrable medium (inhomogeneity), and by a hard/soft obstacle. For each, we provide an accurate representation of the corresponding {far-field pattern}.
    \item As a direct application and demonstration of the framework's utility, we initiate the study of {hyperbolic inverse scattering problems}. We formulate the inverse obstacle and inverse medium scattering problems, where the goal is to recover the shape or material properties of the scatterer from measurements of its far-field pattern.
\end{enumerate}

We emphasize that the hyperbolic radiation condition and Rellich theorem are \emph{local conditions at infinity}. Consequently, our theoretical framework is not confined to pure hyperbolic space; it extends naturally to {asymptotically hyperbolic manifolds}. This work thus fills a foundational gap and paves the way for subsequent developments in both direct and inverse wave scattering on hyperbolic and asymptotically hyperbolic spaces, with potential applications ranging from geometric inverse problems to the study of wave dynamics in AdS/CFT contexts.
\medskip

The remainder of this paper is structured as follows.
\begin{itemize}
    \item Section~\ref{sec:euclidean-background} reviews the fundamental framework of time-harmonic scattering theory in the Euclidean setting, establishing the key concepts and results that motivate our hyperbolic analogue.

    \item Section~\ref{sec:hyperbolic-fourier} develops the necessary Fourier analysis on hyperbolic space and establishes several auxiliary results that are essential for the subsequent analysis.

    \item Section~\ref{sec:direct-theory} is devoted to the core development of the direct scattering theory on hyperbolic space. Here we construct the Green's functions, establish the hyperbolic Sommerfeld radiation condition, prove the Rellich-type theorem, and derive representations for the far-field patterns.

    \item Section~\ref{sec:fourier-approximation} presents several significant results on hyperbolic Fourier approximation. These constitute the analytic foundation for the inverse scattering theory developed in the following sections.

    \item Section~\ref{sec:inverse-obstacle} and Section~\ref{sec:inverse-medium} address the inverse scattering problem. Specifically, Section~\ref{sec:inverse-obstacle} treats the inverse obstacle problem, while Section~\ref{sec:inverse-medium} treats the inverse medium (or potential) scattering problem.
\end{itemize}

\section{The time-harmonic scattering theory in Euclidean space}\label{sec:euclidean-background}

In this section, we provide a brief overview of the classical time-harmonic scattering framework in Euclidean space—or more generally, on noncompact manifolds with Euclidean ends (see the remark below)—governed by the Helmholtz equation. Our aim is twofold: to establish notation and conventions, and more importantly, to highlight the essential structural features—such as radiation conditions, far-field patterns, well-posedness, and uniqueness results—that will serve as a guide for constructing an analogous theory in hyperbolic space. For comprehensive treatments of Euclidean scattering we refer to the monographs \cite{CK1, CK}; the time‑dependent theory is discussed, for example, in \cite{FP, LP}.

\subsection{Direct scattering problems}
Let $k \in \mathbb{R}_+$ denote the wavenumber and let $u^i(x)$, $x \in \mathbb{R}^n$ with $n \geq 2$, be an entire solution to the homogeneous Helmholtz equation:
\begin{equation}\label{eq:helmb1}
-\Delta u^i - k^2 u^i = 0 \quad \text{in } \mathbb{R}^n.
\end{equation}
Physically, $u^i$ represents an incident wave field. The canonical example is the time-harmonic plane wave
\begin{equation}\label{eq:pb1}
u^i(x) = e^{\mathrm{i} k x \cdot d}, \quad \mathrm{i} := \sqrt{-1}, \quad d \in \mathbb{S}^{n-1},
\end{equation}
where the unit vector $d$ indicates the direction of propagation. We note that the plane wave corresponds to a Euclidean Fourier mode $e^{\mathrm{i} x \cdot \xi}$ with wave vector $\xi \in \mathbb{R}^n$ satisfying $|\xi| = k$ and $\xi/|\xi| = d$.
\medskip

The propagation of this incident field is perturbed by the presence of a scatterer, which we model by a bounded function $q(x) = 1 + V(x)$, where $V \in L^2(\mathbb{R}^n)$ has compact support contained in a bounded domain $\Omega$ with connected complement. The real part $\Re q$ relates to the refractive index of an inhomogeneous medium, while $\Im q \geq 0$ represents damping or absorption. Although physically $\Re q > 0$, this positivity condition is not required for the mathematical analysis that follows.

The total wave field $u = u^i + u^s$, where $u^s$ denotes the scattered field, satisfies the perturbed Helmholtz equation
\begin{equation}\label{eq:helmb2}
-\Delta u - k^2(1 + V)u = 0 \quad \text{in } \mathbb{R}^n.
\end{equation}
To single out the physically relevant outgoing solution, we impose the \emph{Sommerfeld radiation condition}
\begin{equation}\label{eq:sfb1}
\lim_{r \to \infty} r^{(n-1)/2} \left( \partial_r u^s - \mathrm{i} k u^s \right) = 0, \quad r := |x|,
\end{equation}
uniformly in $\hat{x} := x/r \in \mathbb{S}^{n-1}$, where $\partial_r u^s(x) := \hat{x} \cdot \nabla_x u^s(x)$. This condition serves as a boundary condition at infinity and guarantees well-posedness of the scattering system. Indeed, for a given incident field $u^i$ and compactly supported potential $V$, there exists a unique solution $u \in H_{\text{loc}}^1(\mathbb{R}^n)$ satisfying \eqref{eq:helmb1}--\eqref{eq:sfb1}.
\medskip

From a PDE perspective, the necessity and sufficiency of the radiation condition are established through \emph{Rellich's theorem}. This theorem asserts that if $v$ satisfies $(-\Delta - k^2)v = 0$ outside a sufficiently large ball $B_R$, obeys the radiation condition \eqref{eq:sfb1}, and moreover satisfies the decay condition
\begin{equation}\label{eq:rellichb1}
\lim_{r \to \infty} \int_{\mathbb{S}^{n-1}} r^{n-1} |v(r \hat{x})|^2  ds(\hat{x}) = 0,
\end{equation}
then $v \equiv 0$ in $\mathbb{R}^n \setminus \overline{B_R}$. Thus, Rellich's theorem establishes a unique continuation property from infinity for radiating solutions of the Helmholtz equation.
\medskip

For later reference, we summarize the complete scattering system for a compactly supported potential $V \in L^2(\Omega)$. The unique scattered field $u^s \in H_{\text{loc}}^1(\mathbb{R}^n)$ satisfies:
\begin{equation}\label{eq:scatb1}
\begin{cases}
-\Delta u^i - k^2 u^i = 0 & \text{in } \mathbb{R}^n, \\[2pt]
-\Delta u - k^2(1 + V)u = 0 & \text{in } \mathbb{R}^n, \\[2pt]
u^s = u - u^i & \text{in } \mathbb{R}^n, \\[2pt]
\displaystyle\lim_{r \to \infty} r^{(n-1)/2} \left( \partial_r u^s - \mathrm{i} k u^s \right) = 0.
\end{cases}
\end{equation}
\medskip

In addition to potential scattering, one may consider scattering by impenetrable obstacles. Let $\Omega$ be a bounded domain with Lipschitz boundary, representing an obstacle that prohibits wave penetration. The physical properties of the obstacle are characterized by a boundary condition imposed on the total field $u$ along $\partial\Omega$. Common examples include:
\begin{itemize}
    \item \emph{Soft (Dirichlet) obstacle:} $u|_{\partial\Omega} = 0$,
    \item \emph{Hard (Neumann) obstacle:} $\partial_\nu u|_{\partial\Omega} = 0$,
    \item \emph{Impedance obstacle:} $(\partial_\nu u + \eta u)|_{\partial\Omega} = 0$,
\end{itemize}
where $\nu$ denotes the exterior unit normal to $\partial\Omega$ and $\eta \in L^\infty(\partial\Omega)$ with $\Im\eta \geq 0$. The forward scattering problem for an obstacle $\Omega$ is then formulated as:
\begin{equation}\label{eq:scatb2}
\begin{cases}
-\Delta u^i - k^2 u^i = 0 & \text{in } \mathbb{R}^n, \\[2pt]
-\Delta u - k^2 u = 0 & \text{in } \mathbb{R}^n \setminus \overline{\Omega}, \\[2pt]
\mathcal{B} u = 0 & \text{on } \partial\Omega, \\[2pt]
u^s = u - u^i & \text{in } \mathbb{R}^n \setminus \overline{\Omega}, \\[2pt]
\displaystyle\lim_{r \to \infty} r^{(n-1)/2} \left( \partial_r u^s - \mathrm{i} k u^s \right) = 0,
\end{cases}
\end{equation}
where $\mathcal{B}$ denotes the appropriate boundary operator. The system \eqref{eq:scatb2} is well-posed and admits a unique solution $u^s \in H_{\text{loc}}^1(\mathbb{R}^n \setminus \overline{\Omega})$.
\medskip

A fundamental property common to scattered fields from both systems \eqref{eq:scatb1} and \eqref{eq:scatb2} is their asymptotic behavior at infinity. As $r \to \infty$, $u^s$ admits the expansion
\begin{equation}\label{eq:ffb1}
u^s(x) = u^s(r\hat{x}) = \frac{e^{\mathrm{i} k r}}{r^{(n-1)/2}} \left[ u_\infty(\hat{x}) + \mathcal{O}\left(\frac{1}{r}\right) \right],
\end{equation}
uniformly in $\hat{x} \in \mathbb{S}^{n-1}$. The function $u_\infty(\hat{x})$, called the \emph{far-field pattern}, encodes complete information about the scattered field; indeed, Rellich's theorem establishes a one-to-one correspondence between $u^s$ and $u_\infty$. When emphasizing parameter dependencies, we write $u_\infty^k(\hat{x})$ for wavenumber $k$, or $u_\infty^k(\hat{x}, d)$ when including dependence on the incident direction $d$. We note that $u_\infty^k(\hat x, d)$ is real-analytic in all its arguments: $k\in\mathbb{R}_+$, $\hat x\in\mathbb{S}^{n-1}$, and $d\in\mathbb{S}^{n-1}$.

In the discussion above, we focused primarily on time-harmonic scattering from compactly supported inhomogeneities in Euclidean space. However, given that the Sommerfeld radiation condition \eqref{eq:sfb1} is formulated locally, the theory can be readily extended to noncompact manifolds with Euclidean ends (cf. \cite{BKLU}).

\subsection{Inverse scattering problems}

The direct scattering framework described above models how compactly supported perturbations—either potential inhomogeneities or impenetrable obstacles—modify wave propagation. This naturally leads to the inverse problem of recovering these perturbations from measurements of the resulting scattering patterns in the far field. Formally, the inverse scattering problem can be written as the operator equation
\begin{equation}\label{eq:ipb1}
\mathcal{F}(\mathfrak{s}) = u_\infty^k(\hat{x}, d),
\end{equation}
where $\mathfrak{s}$ denotes either the potential configuration $(\Omega, V)$ or the obstacle configuration $(\Omega, \mathcal{B})$, and $\mathcal{F}$ is the forward operator defined through the solutions of \eqref{eq:scatb1} or \eqref{eq:scatb2}. Although the forward problem is linear in the data for a fixed scatterer, the inverse problem \eqref{eq:ipb1} is inherently nonlinear and, in the sense of Hadamard, ill-posed.
\medskip

This inverse formulation serves as the mathematical foundation for numerous imaging technologies, including radar, sonar, geophysical exploration, and medical tomography. Its theoretical understanding is therefore of both mathematical and practical significance. While the literature on Euclidean inverse scattering is vast, we focus here on aspects most relevant to our subsequent development of a parallel theory on hyperbolic spaces.
\medskip

A central theoretical question for the inverse problem is that of \emph{unique identifiability}: under what conditions does the far-field data $\{u_\infty^k(\hat{x}, d)\}_{\hat{x},d \in \mathbb{S}^{n-1}}$ uniquely determine the scatterer configuration $\mathfrak{s}$? For a fixed wavenumber $k > 0$, it is known that the full far-field pattern (measured for all incident and observation directions) uniquely determines impenetrable obstacles as well as compactly supported potentials \cite{CK1,CK,DL,LP1,Isa}. The proofs typically employ either spectral arguments or the analysis of wave field singularities induced by point-source incidences. For potential scattering, a key technique establishes an equivalence between the far-field data and the Dirichlet-to-Neumann map on a sufficiently large enclosing surface, thereby reducing the inverse problem to a Calder\'on-type inverse boundary value problem \cite{SU1,U1}.
We refer to \cite{BKLU,CC,GKLU,LL} for complementary directions in inverse scattering theory, including rigorous stability analyses, numerical reconstruction algorithms, and studies of invisibility and transformation-optics based cloaking.

\section{Fourier analysis on hyperbolic space}\label{sec:hyperbolic-fourier}

\subsection{Geometric models of hyperbolic space}
We consider hyperbolic space $\mathbb{H}^n$ as a simply connected, complete, non-compact Riemannian manifold with constant sectional curvature $-1$. By the Killing-Hopf theorem, any two such spaces are isometric. In practice, three equivalent models are commonly used: the hyperboloid (or Lorentz) model, the Poincare ball model, and the half-space model. Each provides distinct computational advantages.

\subsubsection{Hyperboloid model}
Let $\mathbb{R}^{n,1}=(\mathbb{R}^{n+1},g)$ denote Minkowski space equipped with the indefinite metric (Minkowski metric)
\[
ds^2 = dx_1^2+\dots+dx_n^2 - dx_{n+1}^2 .
\]
The hyperboloid model of hyperbolic space is the hyperboloid
\[
\mathbb{H}^n = \{ x \in \mathbb{R}^{n,1} : x_1^2+\dots+x_n^2 - x_{n+1}^2 = -1,\; x_{n+1} > 0 \}
\]
equipped with the Minkowski metric.
We parametrize it by $x' = (x_1,\dots,x_n) \in \mathbb{R}^n$ via
\[
r(x') = \Bigl(x_1,\dots,x_n,\sqrt{1+|x'|^2}\Bigr), \qquad |x'|^2 = \sum_{i=1}^n x_i^2.
\]
The induced metric $g_{\mathbb{H}^n}=g_{ij}\,dx_i\otimes dx_j$ has coefficients
\[
g_{ii}=1+\frac{x_i^2}{1+|x'|^2}, \qquad
g_{ij}= \frac{x_i x_j}{1+|x'|^2} \;(i\neq j).
\]
For $x=(x',x_{n+1}),\,y=(y',y_{n+1}) \in \mathbb{H}^n$, the geodesic distance is given by
\[
\cosh\!\bigl(d_{\mathbb{H}^n}(x,y)\bigr) = x_{n+1}y_{n+1} - x'\cdot y'.
\]

The advantages of hyperboloid model lies in clear symmetric surface (parallel to Euclidean plane) and easily performing moving-plane method.

\subsubsection{Poincaré ball model}
The Poincare ball model is the unit ball $\mathbb{B}^n = \{x\in\mathbb{R}^n : |x|<1\}$ equipped with the conformal metric
\[
g = \left(\frac{2}{1-|x|^2}\right)^{\!2} g_e,
\]
where $g_e$ denotes the Euclidean metric. The hyperbolic volume element is
\[
dV_{\mathbb{B}^n} = \left(\frac{2}{1-|x|^2}\right)^{\!n} dx,
\]
and the geodesic distance from the origin to $x\in\mathbb{B}^n$ is
\[
\rho(x) = \log\frac{1+|x|}{1-|x|}.
\]
The Laplace-Beltrami operator and gradient are given by
\begin{align*}
\Delta_{\mathbb{B}^n} &= \frac{1-|x|^2}{4}\Bigl((1-|x|^2)\Delta_{\mathbb{R}^n} + 2(n-2)\sum_{i=1}^n x_i\frac{\partial}{\partial x_i}\Bigr), \\
\nabla_{\mathbb{B}^n} &= \left(\frac{1-|x|^2}{2}\right)^{\!2} \nabla_{\mathbb{R}^n}.
\end{align*}
For $u,v \in C_c^\infty(\mathbb{B}^n)$, the inner product of gradients satisfies
\[
\bigl(\nabla_{\mathbb{H}^n}u,\nabla_{\mathbb{H}^n}v\bigr)_g
= \left(\frac{1-|x|^2}{2}\right)^{\!2} (\nabla_{\mathbb{R}^n}u,\nabla_{\mathbb{R}^n}v)_{g_e}.
\]
The advantage of Poincare ball model lies in clear expression of differential operator and is used widely in the fields of analysis and PDEs. The establishment of our time-harmonic scattering theory for hyperbolic space mainly uses the Poincare ball model.
\subsubsection{Half-space model}
The half-space model is $\mathcal{R}^n_+ = \mathbb{R}^{n-1}\times(0,\infty)$ with the metric
\[
ds^2 = \frac{dx_1^2+\dots+dx_n^2}{x_n^2}.
\]
Its volume element is $dV_{\mathcal{R}^n_+}=x_n^{-n}dx$, and the hyperbolic gradient and Laplacian are
\[
\nabla_{\mathcal{R}^n_+} = x_n\nabla_{\mathbb{R}^n}, \qquad
\Delta_{\mathcal{R}^n_+} = x_n^2\Delta_{\mathbb{R}^n} - (n-2)x_n\frac{\partial}{\partial x_n}.
\]
For $x=(x',x_n),\,y=(y',y_n)\in\mathcal{R}^n_+$, the geodesic distance satisfies
\[
\cosh\!\bigl(d_{\mathcal{R}^n_+}(x,y)\bigr)
= \frac{x_n^2+y_n^2 + |x'-y'|^2}{2x_ny_n}.
\]
An explicit isometry $\mathbf{M}:\mathcal{R}^n_+\to\mathbb{B}^n$ is given by
\[
\mathbf{M}(x) = \Bigl(
\frac{2x'}{1+2x_n+|x|^2},\,
\frac{1-|x|^2}{1+2x_n+|x|^2}\Bigr), \qquad x=(x',x_n).
\]
The advantage half space model lies in clear scaling and sharp geometric inequality such as Hardy-Sobolv-Mayza inequality on the half space can be equivalently seen as Poincare-Sobolev inequality on the hyperbolic space.
\subsection{Foliations by totally geodesic hypersurfaces}
To analyze symmetric solutions on $\mathbb{H}^n$, we introduce foliations by families of totally geodesic hypersurfaces. In the hyperboloid model, fix a direction $\nu\in\mathbb{S}^{n-1}$ and define the hypersurface through the origin orthogonal to $\nu$:
\[
U_\nu = \{ x\in\mathbb{H}^n : (x_1,\dots,x_n)\cdot\nu = 0 \}.
\]
For $t\in\mathbb{R}$, let $\widetilde{A}_t$ denote the hyperbolic rotation on $\mathbb{R}^{1,1}$:
\[
\widetilde{A}_t = \begin{pmatrix}
\cosh t & \sinh t \\
\sinh t & \cosh t
\end{pmatrix},
\]
and set $A_t^\nu = \operatorname{Id}_{\nu^\perp}\oplus \widetilde{A}_t$. Then the translated hypersurfaces
\[
U_\nu^t = A_t^\nu(U_\nu), \qquad t\in\mathbb{R},
\]
are pairwise disjoint and together foliate $\mathbb{H}^n$. In coordinates, if $x=(y',x_{n+1})\in U_\nu$ with $y'\perp\nu$, then
\[
U_\nu^t = \bigl( y' + (\sinh t)\,x_{n+1}\nu,\; (\cosh t)\,x_{n+1} \bigr).
\]

Under the isometry $\phi:\mathbb{H}^n\to\mathbb{B}^n$ given by $\phi(x)=x'/(x_{n+1}+1)$, the hypersurface $U_\nu$ maps to the hyperplane in $\mathbb{B}^n$:
\[
\phi(U_\nu) = \{ x\in\mathbb{B}^n : x\cdot\nu = 0 \}.
\]
In the Poincaré ball, the full family of totally geodesic hypersurfaces is obtained by applying Mobius transformations (See 3.31). For $a\in\mathbb{B}^n$, define
\[
T_a(x) = \frac{|x-a|^2 a - (1-|a|^2)(x-a)}{1-2x\cdot a + |x|^2|a|^2},
\]
which is an isometry of $\mathbb{B}^n$. Then for any fixed $\nu$, the collection
\[
\bigl\{ T_a\bigl(\phi(U_\nu)\bigr) : a\in\mathbb{B}^n \bigr\}
\]
constitutes all totally geodesic hypersurfaces in $\mathbb{B}^n$. In particular,
\[
\bigcup_{a\in\mathbb{B}^n} T_a\bigl(\phi(U_\nu)\bigr)
= \bigcup_{a\in\mathbb{B}^n} T_a\bigl(\phi(U_{x_n})\bigr) = \mathbb{B}^n.
\]
This foliation structure will be essential for exploiting hyperbolic symmetry and will be used to solve the inverse obstacle problem on hyperbolic space using only single measurement .

\subsection{Fourier analysis on the Poincaré ball $\mathbb{B}^n$}

Under geodesic polar coordinates, the hyperbolic metric on $\mathbb{B}^n$ decomposes as
\[
g = d\rho^2 + \sinh^2\!\rho\; d\sigma,
\]
where $d\sigma$ denotes the standard metric on $\mathbb{S}^{n-1}$ and $\sinh\rho = (e^{\rho/2}-e^{-\rho/2})/2$. In these coordinates, the Laplace–Beltrami operator and gradient take the forms
\begin{align*}
\Delta_{\mathbb{B}^n} &= \frac{\partial^2}{\partial\rho^2} + (n-1)\coth\rho\,\frac{\partial}{\partial\rho} + \frac{1}{\sinh^2\!\rho}\,\Delta_{\mathbb{S}^{n-1}},\\[2mm]
\nabla_{\mathbb{B}^n} &= \Bigl(\frac{\partial}{\partial\rho},\; \frac{1}{\sinh\rho}\nabla_{\mathbb{S}^{n-1}}\Bigr).
\end{align*}
For $f\in L^1(\mathbb{B}^n)$, integration in polar coordinates becomes
\begin{equation}\label{eq:polar-integral}
\begin{split}
\int_{\mathbb{B}^n} f(x)\,dV_{\mathbb{B}^n}
&= \int_{0}^{1}\int_{\mathbb{S}^{n-1}} f(r\xi)\, r^{n-1}\Bigl(\frac{2}{1-r^2}\Bigr)^n d\sigma_\xi\,dr \\
&= \int_{0}^{\infty}\int_{\mathbb{S}^{n-1}} f\!\Bigl(\tanh\frac{\rho}{2}\,\xi\Bigr) (\sinh\rho)^{n-1}\,d\sigma_\xi\,d\rho .
\end{split}
\end{equation}

\begin{lem}[Poincaré inequality on $\mathbb{B}^n$]\label{lem:poincare-hyperbolic}
For any $u\in C_c^\infty(\mathbb{B}^n)$,
\[
\int_{\mathbb{B}^n} |\nabla_{\mathbb{B}^n}u|^2\,dV_{\mathbb{B}^n}
\ge \frac{(n-1)^2}{4}\int_{\mathbb{B}^n} |u|^2\,dV_{\mathbb{B}^n}.
\]
Equality holds only for $u\equiv 0$.
\end{lem}
\begin{proof}
Observe that
\[
\frac{d}{d\rho}(\sinh\rho)^{n-1} = (n-1)(\sinh\rho)^{n-2}\cosh\rho
\ge (n-1)(\sinh\rho)^{n-1}.
\]
Using the polar representation \eqref{eq:polar-integral} and integration by parts,
\begin{align*}
\int_{\mathbb{B}^n} |u|^2\,dV_{\mathbb{B}^n}
&\le \frac{1}{n-1}\int_0^\infty\!\int_{\mathbb{S}^{n-1}} |u|^2\,
\frac{d}{d\rho}(\sinh\rho)^{n-1}\,d\sigma_\xi\,d\rho \\
&= -\frac{2}{n-1}\int_0^\infty\!\int_{\mathbb{S}^{n-1}} u\,
\frac{\partial u}{\partial\rho}\,(\sinh\rho)^{n-1}d\sigma_\xi\,d\rho \\
&\le \frac{2}{n-1}\int_{\mathbb{B}^n} |u|\,|\nabla_{\mathbb{B}^n}u|\,dV_{\mathbb{B}^n} \\
&\le \frac{2}{n-1}
\Bigl(\int_{\mathbb{B}^n}|u|^2\Bigr)^{1/2}
\Bigl(\int_{\mathbb{B}^n}|\nabla_{\mathbb{B}^n}u|^2\Bigr)^{1/2},
\end{align*}
where we used $|\partial u/\partial\rho|\le |\nabla_{\mathbb{B}^n}u|$. Squaring gives the desired inequality. Equality forces $u\equiv0$.
\end{proof}

\begin{rmk}
Lemma~\ref{lem:poincare-hyperbolic} shows that the spectrum of $-\Delta_{\mathbb{B}^n}$ satisfies
\[
\operatorname{Spec}(-\Delta_{\mathbb{B}^n}) \subset \Bigl[\frac{(n-1)^2}{4},\infty\Bigr).
\]
In fact, it is known that $\operatorname{Spec}(-\Delta_{\mathbb{B}^n}) = \bigl[\frac{(n-1)^2}{4},\infty\bigr)$. This contrasts sharply with the Euclidean case, where $\operatorname{Spec}(-\Delta_{\mathbb{R}^n}) = [0,\infty)$ and no Poincaré inequality holds on $W^{1,2}(\mathbb{R}^n)$. This spectral gap reflects the negative curvature of hyperbolic space and underlies many analytic differences between Euclidean and hyperbolic scattering.
\end{rmk}

\subsubsection{Mobius transformations on $\mathbb{B}^n$}
For $a\in\mathbb{B}^n$, define the Möbius transformation $T_a:\mathbb{B}^n\to\mathbb{B}^n$ by
\[
T_a(x) = \frac{|x-a|^2 a - (1-|a|^2)(x-a)}{1-2x\cdot a + |x|^2|a|^2},
\]
where $x\cdot a$ denotes the Euclidean inner product. When $n=2$, writing $a=a_1+ia_2$, $x=x_1+ix_2$, one recovers the familiar form
\[
T_a(x) = \frac{a-x}{1-x\overline{a}} .
\]
These maps satisfy $T_a(0)=a$, $T_a(a)=0$, $|T_a(x)| = |T_x(a)|$, and $T_a\circ T_a = \operatorname{id}$. Crucially, they are isometries of $(\mathbb{B}^n,g)$:
\[
\frac{|dy|^2}{(1-|y|^2)^2} = \frac{|dx|^2}{(1-|x|^2)^2}, \qquad y = T_a(x).
\]
Consequently, for any $\varphi\in L^1(\mathbb{B}^n)$,
\[
\int_{\mathbb{B}^n} |\varphi\circ T_a|\,dV_{\mathbb{B}^n} = \int_{\mathbb{B}^n} |\varphi|\,dV_{\mathbb{B}^n},
\]
and the operators $-\Delta_{\mathbb{B}^n}$ and $\nabla_{\mathbb{B}^n}$ commute with composition by $T_a$ in the sense that for $\phi,\psi\in C_c^\infty(\mathbb{B}^n)$,
\begin{align*}
\int_{\mathbb{B}^n} -\Delta_{\mathbb{B}^n}(\phi\circ T_a)\,(\psi\circ T_a)\,dV_{\mathbb{B}^n}
&= \int_{\mathbb{B}^n} (-\Delta_{\mathbb{B}^n}\phi)\,\psi\,dV_{\mathbb{B}^n},\\
\int_{\mathbb{B}^n} \nabla_{\mathbb{B}^n}(\phi\circ T_a)\cdot\nabla_{\mathbb{B}^n}(\psi\circ T_a)\,dV_{\mathbb{B}^n}
&= \int_{\mathbb{B}^n} \nabla_{\mathbb{B}^n}\phi\cdot\nabla_{\mathbb{B}^n}\psi\,dV_{\mathbb{B}^n}.
\end{align*}
The geodesic distance $\rho(x,y)$ can be expressed via $T_a$ as
\[
\rho(x,y) = \rho(T_x(y),0) = \log\frac{1+|T_x(y)|}{1-|T_x(y)|}.
\]

Using Mobius transformations, we define the \emph{hyperbolic convolution} of measurable functions $f,g$ on $\mathbb{B}^n$ by
\[
(f\ast g)(x) = \int_{\mathbb{B}^n} f(y)\,g(T_x(y))\,dV_{\mathbb{B}^n}(y).
\]
Since $T_x^{-1}=T_x$, this convolution is symmetric:
\[
(f\ast g)(x) = (g\ast f)(x).
\]
If $g$ is radial, i.e., $g(x)=g(\rho(x,0))$, then
\[
(f\ast g)(x) = \int_{\mathbb{B}^n} f(y)\,g(\rho(x,y))\,dV_{\mathbb{B}^n}(y).
\]

\subsubsection{Fourier transform on $\mathbb{B}^n$}
We now recall basic facts about the Fourier transform on hyperbolic space, detailed accounts can be found in \cite{Helgason1,Helgason2,LuYangQ1,LuYangQ2,LuYangQ3}. Define the generalized eigenfunctions
\[
e_{\lambda,\xi}(x) = \Bigl(\frac{\sqrt{1-|x|^2}}{|x-\xi|}\Bigr)^{n-1+\mathrm{i}\lambda},
\qquad \lambda\in\mathbb{R},\; \xi\in\mathbb{S}^{n-1},
\]
and their spherical average
\[
\phi_\lambda(x) = |\mathbb{S}^{n-1}|^{-1}\int_{\mathbb{S}^{n-1}} e_{\lambda,\xi}(x)\,d\sigma_\xi.
\]
These satisfy
\[
-\Delta_{\mathbb{B}^n} e_{\lambda,\xi} = \frac{(n-1)^2+\lambda^2}{4}\,e_{\lambda,\xi},
\qquad
\phi_\lambda(T_x(y)) = \int_{\mathbb{S}^{n-1}} e_{\lambda,\xi}(x)\,e_{\lambda,-\xi}(y)\,d\sigma_\xi.
\]
The function $e_{-\lambda,\xi}(x)$ plays a role analogous to the Euclidean plane wave $e^{\mathrm{i}\lambda x\cdot\xi}$.
\medskip

For $f\in C_c^\infty(\mathbb{B}^n)$, the \emph{hyperbolic Fourier transform} is defined as
\[
\widehat{f}(\lambda,\xi) = \int_{\mathbb{B}^n} f(x)\,e_{-\lambda,\xi}(x)\,dV_{\mathbb{B}^n}.
\]
It satisfies the convolution identity $(\widehat{f\ast g}) = \widehat{f}\,\widehat{g}$ and admits the inversion formula
\[
f(x) = D_n\int_{-\infty}^{\infty}\int_{\mathbb{S}^{n-1}}
\widehat{f}(\lambda,\xi)\,e_{\lambda,\xi}(x)\,
|c(\lambda)|^{-2}\,d\lambda\,d\xi,
\]
where $D_n = \bigl(2^{3-n}\pi|\mathbb{S}^{n-1}|\bigr)^{-1}$ and
\[
c(\lambda) =
\frac{2^{n-1}\,\Gamma(n/2)\,\Gamma(\mathrm{i}\lambda)}
{\Gamma\bigl(\frac{n-1+\mathrm{i}\lambda}{2}\bigr)\,
\Gamma\bigl(\frac{1+\mathrm{i}\lambda}{2}\bigr)}
\]
is the Harish-Chandra function \cite{liu}. The weight $|c(\lambda)|^{-2}$ behaves like $\lambda^{n-1}$ as $\lambda\to\infty$ and vanishes as $\lambda\to0$.
For $g\in L^{2}(\mathbb{B}^n)$, $h\in L^{2}(\mathbb{B}^n)$, the Plancherel formula on the hyperbolic space
\[
\int_{\mathbb{B}^n} g(x)\,\overline{h(x)}\,dV_{\mathbb{B}^n}
= D_n\int_{-\infty}^{\infty}\int_{\mathbb{S}^{n-1}}
\widehat{g}(\lambda,\xi)\,\overline{\widehat{h}(\lambda,\xi)}\,
|c(\lambda)|^{-2}\,d\lambda\,d\xi.
\]
still holds.
Since $e_{\lambda,\xi}(x)$ ia an eigen-function of $-\Delta_{\mathbb{B}^n}$ with eigenvalue equal to $\frac{(n-1)^2+\lambda^2}{4}$, then for $f\in C^{\infty}_c(\mathbb{B}^n)$, simple calculation gives that
\begin{equation}\begin{split}
\widehat{-\Delta_{\mathbb{B}^n} f}(\lambda, \xi)&=\int_{\mathbb{B}^n}-\Delta_{\mathbb{B}^n}(f)e_{-\lambda,\xi}(x)dV_{\mathbb{B}^n}\\
&=\int_{\mathbb{B}^n} -\Delta_{\mathbb{B}^n}(e_{-\lambda,\xi})fdV_{\mathbb{B}^n}\\
&=\frac{(n-1)^2+\lambda}{4}\widehat{f}(\lambda, \xi).
\end{split}\end{equation}

\subsection{Green's formulas on bounded hyperbolic domains}

\begin{lem}[Green's identities on a hyperbolic ball]\label{lem:green-ball}
Let $B_{\mathbb{B}^n}(0,R)$ be the hyperbolic ball of radius $R$ centered at $0$. For $u,v\in C^2(\overline{B_{\mathbb{B}^n}(0,R)})$, there holds
\begin{align}
&\int_{B_{\mathbb{B}^n}(0,R)} -\Delta_{\mathbb{B}^n}u\;v\,dV_{\mathbb{B}^n}
= \int_{B_{\mathbb{B}^n}(0,R)} u\;(-\Delta_{\mathbb{B}^n}v)\,dV_{\mathbb{B}^n}
+ \int_{\partial B_{\mathbb{B}^n}(0,R)}
\Bigl(\frac{\partial v}{\partial\rho}u - \frac{\partial u}{\partial\rho}v\Bigr)\,d\sigma_{\mathbb{B}^n},
\label{eq:green1}\\
&\int_{B_{\mathbb{B}^n}(0,R)} -\Delta_{\mathbb{B}^n}u\;v\,dV_{\mathbb{B}^n}
= \int_{B_{\mathbb{B}^n}(0,R)} (\nabla_{\mathbb{B}^n}u,\nabla_{\mathbb{B}^n}v)_g\,dV_{\mathbb{B}^n}
- \int_{\partial B_{\mathbb{B}^n}(0,R)} \frac{\partial u}{\partial\rho}\,v\,d\sigma_{\mathbb{B}^n},
\label{eq:green2}
\end{align}
where $d\sigma$ denotes the surface measure on the boundary of Euclidean ball $B^{n}(0, \tanh(\frac{R}{2}))$ and $d\sigma_{\mathbb{B}^n}$ denotes the surface measure on the boundary of hyperbolic ball $B_{\mathbb{B}^n}(0,R)$.
\end{lem}
\begin{proof}
Recall the conformal laws between the hyperbolic Laplacian and the Euclidean Laplacian:
\[
\Bigl(-\Delta_{\mathbb{B}^n}-\frac{n(n-2)}{4}\Bigr)u
= -\Bigl(\frac{2}{1-|x|^2}\Bigr)^{-\frac{n}{2}-1}
\Delta_{\mathbb{R}^n}\!\Bigl[\Bigl(\frac{2}{1-|x|^2}\Bigr)^{\frac{n}{2}-1}u\Bigr].
\]
Set $\widetilde{u} = \bigl(\frac{2}{1-|x|^2}\bigr)^{\frac{n}{2}-1}u$ and similarly for $\widetilde{v}$. In Euclidean polar coordinates $|x| = r = \tanh\frac{\rho}{2}$, we compute
\begin{align*}
\int_{B_{\mathbb{B}^n}(0,R)} -\Delta_{\mathbb{B}^n}u\;v\,dV_{\mathbb{B}^n}
&= \int_{B(0,\tanh\frac{R}{2})} -\Delta_{\mathbb{R}^n}\widetilde{u}\;\widetilde{v}\,dx
   + \frac{n(n-2)}{4}\int_{B_{\mathbb{B}^n}(0,R)} uv\,dV_{\mathbb{B}^n} \\
&= \int_{B(0,\tanh\frac{R}{2})} \widetilde{u}\;(-\Delta_{\mathbb{R}^n}\widetilde{v})\,dx
   + \int_{\partial B(0,\tanh\frac{R}{2})}
      \Bigl(\frac{\partial\widetilde{v}}{\partial r}\widetilde{u}
      - \frac{\partial\widetilde{u}}{\partial r}\widetilde{v}\Bigr)d\sigma \\
&\quad+ \frac{n(n-2)}{4}\int_{B_{\mathbb{B}^n}(0,R)} uv\,dV_{\mathbb{B}^n}\\
&=\int_{B_{\mathbb{B}^n}(0,R)}-\Delta_{\mathbb{B}^n}(v)udV_{\mathbb{B}^n}+\int_{\partial B^{n}(0, \tanh(\frac{R}{2}))}\left(\frac{\partial \widetilde{v}}{\partial r}\widetilde{u}-\frac{\partial \widetilde{u}}{\partial r}\tilde{v}\right)d\sigma.
\end{align*}
The boundary term transforms as follows. Noting that $\frac{\partial\rho}{\partial r}=2\cosh^2\frac{\rho}{2}$,
\[
\frac{\partial\widetilde{u}}{\partial r}
= 2^{\frac{n}{2}-1}(n-2)\cosh^{n-1}\!\frac{\rho}{2}\,\sinh\frac{\rho}{2}\,u
+ 2^{\frac{n}{2}}\cosh^n\!\frac{\rho}{2}\,\frac{\partial u}{\partial\rho},
\]
and similarly for $\widetilde{v}$. Hence the boundary integral $$\int_{\partial B^{n}(0, \tanh(\frac{R}{2}))}\left(\frac{\partial \widetilde{v}}{\partial r}\widetilde{u}-\frac{\partial \widetilde{u}}{\partial r}\tilde{v}\right)d\sigma$$
can be written as
\begin{align*}
&\int_{\partial B(0,\tanh\frac{R}{2})}
\Bigl(\frac{\partial\widetilde{v}}{\partial r}\widetilde{u}
- \frac{\partial\widetilde{u}}{\partial r}\widetilde{v}\Bigr)d\sigma \\
&= 2^{n-1}\Bigl(\cosh\frac{R}{2}\Bigr)^{2n-2}
   \int_{\partial B(0,\tanh\frac{R}{2})}
   \Bigl(\frac{\partial v}{\partial\rho}u - \frac{\partial u}{\partial\rho}v\Bigr)d\sigma \\
&= \int_{\partial B_{\mathbb{B}^n}(0,R)}
   \Bigl(\frac{\partial v}{\partial\rho}u - \frac{\partial u}{\partial\rho}v\Bigr)d\sigma_{\mathbb{B}^n},
\end{align*}
 where $d\sigma_{\mathbb{B}^n}|_{\partial B_{\mathbb{B}^n}(0,R)}=\left(\frac{2}{1-|x|^2}\right)^{n-1}d\sigma|_{\partial B^{n}(0, \tanh(\frac{R}{2}))}$ denotes the hyperbolic surface measure in $\partial B_{\mathbb{B}^n}(0,R)$. This proves \eqref{eq:green1}. Identity \eqref{eq:green2} follows similarly by applying the Euclidean Green's formula to $\widetilde{u},\widetilde{v}$ and transforming the boundary term.
\end{proof}

\subsection{Outward normal derivatives on bounded domains}

Let $\Omega\subset\mathbb{B}^n$ be a smooth bounded domain, we try to give the accurate calculation formula for the hyperbolic out-normal derivatives $\frac{\partial}{\partial \nu_{\mathbb{B}^n}}$ along the boundary $\partial \Omega$ through integration by parts. In fact, for any $u\in W^{1,2}(\mathbb{B}^n)$ and $v\in W^{1,2}(\mathbb{B}^n)$, careful computations gives that
\begin{equation}\begin{split}
&\int_{\Omega}-\Delta_{\mathbb{B}^n}(u)vdV_{\mathbb{B}^n}\\
&\ \ =\int_{\Omega}\left(-\Delta_{\mathbb{B}^n} u-\frac{n(n-2)}{4}u\right)vdV_{\mathbb{B}^n}+\int_{\Omega}\frac{n(n-2)}{4}uvdV_{\mathbb{B}^n}\\
&\ \ =\int_{\Omega}\left(\frac{2}{1-|x|^2}\right)^{-\frac{n}{2}-1}\left(-\Delta_{\mathbb{R}^n}\left(\left(\frac{2}{1-|x|^2}\right)^{\frac{n}{2}-1}u\right)v\left(\frac{2}{1-|x|^2}\right)^{n}\right)dx\\
&\ \ \ \ +\int_{\Omega}\frac{n(n-2)}{4}uvdV_{\mathbb{B}^n}\\
&\ \ =\int_{\Omega}-\Delta_{\mathbb{R}^n}(\widetilde{u})\tilde{v}dx+\int_{\Omega}\frac{n(n-2)}{4}uvdV_{\mathbb{B}^n}\\
&\ \ =\int_{\Omega}-\Delta_{\mathbb{B}^n}(v)udV_{\mathbb{B}^n}+\int_{\partial \Omega}\left(\frac{\partial \widetilde{v}}{\partial \nu}\widetilde{u}-\frac{\partial \widetilde{u}}{\partial \nu}\widetilde{v}\right)d\sigma
\end{split}\end{equation}
where $\widetilde{u}=\left(\frac{2}{1-|x|^2}\right)^{\frac{n}{2}-1}u$, $\widetilde{v}=\left(\frac{2}{1-|x|^2}\right)^{\frac{n}{2}-1}v$.
\vskip 0.1cm

We define the \emph{hyperbolic outward normal derivative} $\partial/\partial\nu_{\mathbb{B}^n}$ on $\partial\Omega$ via the conformal factor. For $u\in W^{1,2}(\Omega)$, set $\widetilde{u} = \bigl(\frac{2}{1-|x|^2}\bigr)^{\frac{n}{2}-1}u$ and let $\nu$ denote the Euclidean outward normal on $\partial\Omega$. Define
\[
\frac{\partial u}{\partial\nu_{\mathbb{B}^n}}
:= \Bigl(\frac{2}{1-|x|^2}\Bigr)^{-\frac{n}{2}}
\frac{\partial\widetilde{u}}{\partial\nu}.
\]
Then for $u,v\in C^2(\overline{\Omega})$, there holds
\begin{align}
\int_{\Omega} -\Delta_{\mathbb{B}^n}u\;v\,dV_{\mathbb{B}^n}
&= \int_{\Omega} u\;(-\Delta_{\mathbb{B}^n}v)\,dV_{\mathbb{B}^n}
+ \int_{\partial\Omega}
   \Bigl(\frac{\partial v}{\partial\nu_{\mathbb{B}^n}}u
   - \frac{\partial u}{\partial\nu_{\mathbb{B}^n}}v\Bigr)d\sigma_{\mathbb{B}^n},
\label{eq:green-gen1}\\
\int_{\Omega} -\Delta_{\mathbb{B}^n}u\;v\,dV_{\mathbb{B}^n}
&= \int_{\Omega} (\nabla_{\mathbb{B}^n}u,\nabla_{\mathbb{B}^n}v)_g\,dV_{\mathbb{B}^n}
- \int_{\partial\Omega} \frac{\partial u}{\partial\nu_{\mathbb{B}^n}}\,v\,d\sigma_{\mathbb{B}^n},
\label{eq:green-gen2}
\end{align}
where $d\sigma_{\mathbb{B}^n} = \bigl(\frac{2}{1-|x|^2}\bigr)^{n-1}d\sigma$. These formulas extend the classical Green identities to the hyperbolic setting and will be essential for deriving integral representations of scattered fields.

\section{Time-harmonic scattering theory on hyperbolic space}\label{sec:direct-theory}

\subsection{Heat kernel and Green's function of the shifted Laplacian}

We begin with the heat equation on the Poincaré ball:
\[
\begin{cases}
\partial_t u(x,t) - \Delta_{\mathbb{B}^n} u(x,t) = 0, & (x,t)\in \mathbb{B}^n\times \mathbb{R}^+,\\[2mm]
u(x,0) = f(x)\in C_c^\infty(\mathbb{B}^n).
\end{cases}
\]
Applying the hyperbolic Fourier transform yields
\[
\partial_t \widehat{u}(\lambda,\xi,t) + \frac{(n-1)^2+\lambda^2}{4}\,\widehat{u}(\lambda,\xi,t) = 0,
\qquad \widehat{u}(\lambda,\xi,0)=\widehat{f}(\lambda,\xi),
\]
which is solved by
\[
\widehat{u}(\lambda,\xi,t) = e^{-\frac{(n-1)^2+\lambda^2}{4}t}\,\widehat{f}(\lambda,\xi).
\]
Inverting the transform gives the solution as a convolution with the heat kernel:
\[
u(x,t) = (P_t * f)(x) := \int_{\mathbb{B}^n} P_t(T_x(y))\,f(y)\,dV_{\mathbb{B}^n}(y) = e^{t\Delta_{\mathbb{B}^n}}f,
\]
where the heat kernel $P_t$ admits the spectral representation
\[
P_t(x) = D_n\int_{-\infty}^{\infty}\int_{\mathbb{S}^{n-1}}
e^{-\frac{(n-1)^2+\lambda^2}{4}t}
e_{\lambda,\xi}(x)\,|c(\lambda)|^{-2}\,d\lambda\,d\xi.
\]
Since $P_t$ depends only on the hyperbolic distance $\rho(x)$, we write $P_t(\rho)$. Explicit formulas are known (see \cite{Davis,Gri}): for $n=2m$ even,
\[
P_t(\rho) = (2\pi)^{-\frac{n+1}{2}} t^{-\frac12} e^{-\frac{(n-1)^2}{4}t}
\int_\rho^\infty \frac{\sinh r}{\sqrt{\cosh r-\cosh\rho}}
\Bigl(-\frac{1}{\sinh r}\frac{\partial}{\partial r}\Bigr)^m e^{-\frac{r^2}{4t}}\,dr,
\]
and for $n=2m+1$ odd,
\[
P_t(\rho) = 2^{-m-1}\pi^{-m-\frac12} t^{-\frac12} e^{-\frac{(n-1)^2}{4}t}
\Bigl(-\frac{1}{\sinh\rho}\frac{\partial}{\partial\rho}\Bigr)^m e^{-\frac{\rho^2}{4t}}.
\]

Because $\operatorname{Spec}(-\Delta_{\mathbb{B}^n})\subset[\frac{(n-1)^2}{4},\infty)$, the operator $-\Delta_{\mathbb{B}^n}-\frac{(n-1)^2}{4}+k^2$ is invertible for $k>0$. By the Mellin (or Laplace) transform,
\begin{equation}\label{eq:inverse-shifted}
\Bigl(-\Delta_{\mathbb{B}^n}-\frac{(n-1)^2}{4}+k^2\Bigr)^{-1}
= \int_0^\infty e^{(\frac{(n-1)^2}{4}-k^2)t}\, e^{t\Delta_{\mathbb{B}^n}}\,dt.
\end{equation}
An explicit formula for the corresponding Green's function is available (see \cite{li, Mat, LuYangQ2, LuYangQ3}): for $n\ge 3$,
\begin{align}
\Bigl(-\Delta_{\mathbb{B}^n}-\frac{(n-1)^2}{4}+k^2\Bigr)^{-1}
&= (2\pi)^{-\frac{n}{2}} (\sinh\rho)^{-\frac{n-2}{2}}
e^{-\frac{(n-2)\pi}{2}\mathrm{i}}
Q_{k-\frac12}^{\frac{n-2}{2}}(\cosh\rho) \label{eq:legendre-form} \\
&= \frac{A_{n,k}}{(\sinh\rho)^{n-2}}
\int_0^\pi (\cosh\rho+\cos t)^{\frac{n-3}{2}-k} (\sin t)^{2k}\,dt, \label{eq:integral-form}
\end{align}
where $Q_{\nu}^{\mu}$ denotes the associated Legendre function of the second kind \cite{Erd}, and
\[
A_{n,k} = (2\pi)^{-\frac{n}{2}}
\frac{\Gamma\bigl(\frac{n-1}{2}+k\bigr)}
{2^{k+\frac12}\,\Gamma\bigl(k+\frac12\bigr)}.
\]
The integral representation \eqref{eq:integral-form} follows from the identity
\[
e^{-\frac{(n-2)\pi}{2}\mathrm{i}}
Q_{k-\frac12}^{\frac{n-2}{2}}(\cosh\rho)
= \frac{\Gamma\bigl(\frac{n-1}{2}+k\bigr)}
{2^{k+\frac12}\,\Gamma\bigl(k+\frac12\bigr)\,\sinh^{\frac{n-2}{2}}\rho}
\int_0^\pi (\cosh\rho+\cos t)^{\frac{n-3}{2}-k} (\sin t)^{2k}\,dt.
\]

\begin{rmk}
The formulas \eqref{eq:legendre-form} and \eqref{eq:integral-form} are valid for $k>0$. They will serve as the starting point for constructing the Green's functions of the Helmholtz operator in the next subsection via analytic continuation in the spectral parameter.
\end{rmk}

\subsection{Green's function for the Helmholtz operator}

We now seek fundamental solutions for the time-harmonic Helmholtz equation on hyperbolic space. Consider the time-dependent wave equation
\begin{equation}\label{eq:wave-eq}
\partial_{tt} w(x,t) - \Bigl(\Delta_{\mathbb{B}^n} w(x,t) + \frac{(n-1)^2}{4} w(x,t)\Bigr) = 0,
\qquad (x,t)\in \mathbb{B}^n\times\mathbb{R}^+.
\end{equation}
Looking for solutions of the form $w(x,t)=e^{\mathrm{i}\mu t}u(x)$ leads to the stationary Helmholtz equation
\begin{equation}\label{eq:helmholtz}
\Bigl(-\Delta_{\mathbb{B}^n} - \frac{(n-1)^2}{4} - \mu^2\Bigr)u = 0,
\qquad x\in\mathbb{B}^n,\; \mu>0.
\end{equation}
Our goal is to construct the Green's function (fundamental solution) for the operator
\[
\mathcal{L}_\mu := -\Delta_{\mathbb{B}^n} - \frac{(n-1)^2}{4} - \mu^2.
\]

\subsubsection{Analytic continuation from the shifted Laplacian\\}

Through Helgason-Fourier transform on hyperbolic space, we see that
\begin{equation}\begin{split}
\left(-\Delta_{\mathbb{B}^n}-z\right)^{-1}&=\int_{0}^{+\infty}e^{tz}e^{t\Delta_{\mathbb{B}^n}}dt\\
\end{split}\end{equation}
holds if $\Re(z)\leq \frac{(n-1)^2}{4}$ and fails if $\Re(z)>\frac{(n-1)^2}{4}$. Hence one can not apply the technique of Fourier transform to derive the formula of $(-\Delta_{\mathbb{B}^n}-\frac{(n-1)2}{4}-\mu^2)^{-1}$. It is well known that the spectrum of the $-\Delta_{\mathbb{B}^n}$ is $(\frac{(n-1)^2}{4}, +\infty)$.  From \cite{Guillarmou}, we know that $(-\Delta_{\mathbb{B}^n}-z)^{-1}$ is holomorphic with respect to $z$ when $z$ is far from the line $(\frac{(n-1)^2}{4}, +\infty)$.
\medskip

Recall from the previous subsection that for $k>0$, we have the explicit Green's function
\[
\Bigl(-\Delta_{\mathbb{B}^n}-\frac{(n-1)^2}{4}+k^2\Bigr)^{-1}
= \frac{A_{n,k}}{(\sinh\rho)^{n-2}}
\int_0^\pi (\cosh\rho+\cos t)^{\frac{n-3}{2}-k} (\sin t)^{2k}\,dt,
\]
with $A_{n,k} = (2\pi)^{-n/2} \Gamma(\frac{n-1}{2}+k) / \bigl(2^{k+1/2}\Gamma(k+\frac12)\bigr)$.
The right-hand side depends on $k$ only through the parameters $A_{n,k}$ and the exponents in the integrand. Both
\[
A_{n,z} = (2\pi)^{-n/2}
\frac{\Gamma\bigl(\frac{n-1}{2}+z\bigr)}
{2^{z+1/2}\,\Gamma\bigl(z+\frac12\bigr)},
\]
and the integral
\[
I_z(\rho) := \int_0^\pi (\cosh\rho+\cos t)^{\frac{n-3}{2}-z} (\sin t)^{2z}\,dt
\]
are well-defined and holomorphic for $\Re z > -\frac12$. By the uniqueness of holomorphic extension, we deduce that $$(-\Delta_{\mathbb{B}^n}-\frac{(n-1)^2}{4}+z^2)^{-1}=\frac{A_{n,z}}{(\sinh \rho)^{n-2}}\int_{0}^{\pi}(\cosh \rho+\cos t)^{\frac{n-3}{2}-z}(\sin t)^{2z}dt$$
for any complex $z$ satisfying $z^2\in \mathbb{C}\setminus (-\infty, 0)$ and $Re(z)>-\frac{1}{2}$.
\vskip0.3cm

Observe that the Helmholtz operator $\mathcal{L}_\mu$ corresponds formally to taking $z = \pm\mathrm{i}\mu$, which lies outside the original domain of definition $k>0$. However, the analytically continued family $$F_z(\rho):=\frac{A_{n,z}}{(\sinh\rho)^{n-2}}I_{z}(\rho)$$ allows us to define the desired Green's functions via boundary values.

\subsubsection{Outgoing and ingoing Green's functions}

For $\mu>0$ and $\epsilon>0$ small, consider the two families
\[
\Bigl(-\Delta_{\mathbb{B}^n}-\frac{(n-1)^2}{4}+(\epsilon\pm\mathrm{i}\mu)^2\Bigr)^{-1}
= F_{\epsilon\pm\mathrm{i}\mu}(\rho),
\]
Letting $\epsilon\to 0^+$ yields two distinct limits:

\begin{align*}
G_{-\mu\mathrm{i}}(\rho) &:=
\lim_{\epsilon\to0^+} F_{\epsilon-\mathrm{i}\mu}(\rho) \\[1mm]
&= A_{n,-\mu\mathrm{i}}\,
\frac{(\cosh\rho)^{\frac{n-3}{2}+\mu\mathrm{i}}}{(\sinh\rho)^{n-2}}
\int_0^\pi \Bigl(1+\frac{\cos t}{\cosh\rho}\Bigr)^{\frac{n-3}{2}+\mu\mathrm{i}}
(\sin t)^{-2\mu\mathrm{i}}\,dt, \\[3mm]
G_{\mu\mathrm{i}}(\rho) &:=
\lim_{\epsilon\to0^+} F_{\epsilon+\mathrm{i}\mu}(\rho) \\[1mm]
&= A_{n,\mu\mathrm{i}}\,
\frac{(\cosh\rho)^{\frac{n-3}{2}-\mu\mathrm{i}}}{(\sinh\rho)^{n-2}}
\int_0^\pi \Bigl(1+\frac{\cos t}{\cosh\rho}\Bigr)^{\frac{n-3}{2}-\mu\mathrm{i}}
(\sin t)^{2\mu\mathrm{i}}\,dt,
\end{align*}
where $A_{n,\pm\mu\mathrm{i}}$ are obtained by replacing $z$ with $\pm\mu\mathrm{i}$ in the formula for $A_{n,z}$.
\medskip

These two functions are candidates for the fundamental solutions of $\mathcal{L}_\mu$. Their asymptotic behaviour at infinity will distinguish them as \emph{outgoing} ($G_{-\mu\mathrm{i}}$) and \emph{ingoing} ($G_{\mu\mathrm{i}}$) waves, in analogy with the Euclidean case $e^{\pm\mathrm{i}k|x|}/|x|^{(n-1)/2}$.

\begin{thm}[Fundamental solutions of the Helmholtz operator]\label{thm:helmholtz-green}
Define
\begin{align*}
G_{-\mu\mathrm{i}}(\rho(x)) &:=
A_{n,-\mu\mathrm{i}}\,
\frac{(\cosh\rho(x))^{\frac{n-3}{2}+\mu\mathrm{i}}}{(\sinh\rho(x))^{n-2}}
\int_0^\pi \Bigl(1+\frac{\cos t}{\cosh\rho(x)}\Bigr)^{\frac{n-3}{2}+\mu\mathrm{i}}
(\sin t)^{-2\mu\mathrm{i}}\,dt, \\[2mm]
G_{\mu\mathrm{i}}(\rho(x)) &:=
A_{n,\mu\mathrm{i}}\,
\frac{(\cosh\rho(x))^{\frac{n-3}{2}-\mu\mathrm{i}}}{(\sinh\rho(x))^{n-2}}
\int_0^\pi \Bigl(1+\frac{\cos t}{\cosh\rho(x)}\Bigr)^{\frac{n-3}{2}-\mu\mathrm{i}}
(\sin t)^{2\mu\mathrm{i}}\,dt.
\end{align*}
Then $G_{-\mu\mathrm{i}}$ satisfies
\[
\Bigl(-\Delta_{\mathbb{B}^n} - \frac{(n-1)^2}{4} - \mu^2\Bigr) G_{-\mu\mathrm{i}}(\rho(x)) = \delta_0(x),
\qquad x\in\mathbb{B}^n,
\]
with the decay estimate $G_{-\mu\mathrm{i}}(\rho(x)) = O\bigl(\sinh^{-\frac{n-1}{2}}\rho(x)\bigr)$ as $\rho(x)\to\infty$.

Moreover, for any $x,y\in\mathbb{B}^n$,
\begin{align*}
\Bigl(-\Delta_{\mathbb{B}^n}\Bigr)_y G_{-\mu\mathrm{i}}(\rho(x,y))
- \frac{(n-1)^2}{4} G_{-\mu\mathrm{i}}(\rho(x,y))
- \mu^2 G_{-\mu\mathrm{i}}(\rho(x,y)) &= \delta_x(y), \\[1mm]
\Bigl(-\Delta_{\mathbb{B}^n}\Bigr)_y G_{\mu\mathrm{i}}(\rho(x,y))
- \frac{(n-1)^2}{4} G_{\mu\mathrm{i}}(\rho(x,y))
- \mu^2 G_{\mu\mathrm{i}}(\rho(x,y)) &= \delta_x(y),
\end{align*}
where $\rho(x,y)$ denotes the hyperbolic distance. Both $G_{\pm\mu\mathrm{i}}(\rho(x,y))$ satisfy
\[
G_{\pm\mu\mathrm{i}}(\rho(x,y)) = O\bigl(\sinh^{-\frac{n-1}{2}}\rho(y)\bigr)
\quad\text{as }\rho(y)\to\infty.
\]
\end{thm}

\begin{proof}
For $\epsilon>0$, let $G_{\epsilon-\mu\mathrm{i}}$ denote the Green's function of the operator
$-\Delta_{\mathbb{B}^n}-\frac{(n-1)^2}{4}+(\epsilon-\mu\mathrm{i})^2$. For any $\phi\in C_c^\infty(\mathbb{B}^n)$,
\[
\int_{\mathbb{B}^n} G_{\epsilon-\mu\mathrm{i}}(\rho(x))
\Bigl(-\Delta_{\mathbb{B}^n}\phi-\frac{(n-1)^2}{4}\phi+(\epsilon-\mu\mathrm{i})^2\phi\Bigr)dV_{\mathbb{B}^n}
= \phi(0).
\]
Since $G_{\epsilon-\mu\mathrm{i}}(\rho(x)) \lesssim (\cosh\rho)^{\frac{n-3}{2}}/(\sinh\rho)^{n-2}$ and
$\sinh\rho \sim |x|$ as $\rho\to0$, we may apply the dominated convergence theorem as $\epsilon\to0^+$ to obtain
\[
\int_{\mathbb{B}^n} G_{-\mu\mathrm{i}}(\rho(x))
\Bigl(-\Delta_{\mathbb{B}^n}\phi-\frac{(n-1)^2}{4}\phi-\mu^2\phi\Bigr)dV_{\mathbb{B}^n}
= \phi(0),
\]
which establishes that $G_{-\mu\mathrm{i}}$ is the Green's function with singularity at the origin.

The translation property follows from the Mobius invariance of $\Delta_{\mathbb{B}^n}$. For any $\phi\in C_c^\infty(\mathbb{B}^n)$ and fixed $x$, let $\tau_x$ denote the Mobius transformation that maps $0$ to $x$. Then
\begin{align*}
&\int_{\mathbb{B}^n} G_{-\mu\mathrm{i}}(\rho(x,y))
\Bigl(-\Delta_{\mathbb{B}^n}\phi(y)-\frac{(n-1)^2}{4}\phi(y)-\mu^2\phi(y)\Bigr)dV_{\mathbb{B}^n}(y)\\
&= \int_{\mathbb{B}^n} G_{-\mu\mathrm{i}}(\rho(y))
\Bigl(-\Delta_{\mathbb{B}^n}(\phi\circ\tau_x)(y)-\frac{(n-1)^2}{4}(\phi\circ\tau_x)(y)-\mu^2(\phi\circ\tau_x)(y)\Bigr)dV_{\mathbb{B}^n}(y)\\
&= (\phi\circ\tau_x)(0) = \phi(x),
\end{align*}
proving that $G_{-\mu\mathrm{i}}(\rho(x,y))$ satisfies $\mathcal{L}_\mu G_{-\mu\mathrm{i}}(\rho(x,y)) = \delta_x(y)$. The argument for $G_{\mu\mathrm{i}}$ is identical.
\end{proof}

\begin{rmk}
In view of the time-dependent equation \eqref{eq:wave-eq}, solutions of the form $e^{\mathrm{i}\mu t}G_{-\mu\mathrm{i}}$ represent outgoing spherical waves, while $e^{\mathrm{i}\mu t}G_{\mu\mathrm{i}}$ represent ingoing waves. For scattering problems we will always select the outgoing Green's function $G_{-\mu\mathrm{i}}$, which corresponds to the physically relevant condition that energy radiates toward infinity.
\end{rmk}

\subsection{Accurate asymptotic behaviors of Green's function $G_{-\mu \mathrm{i}}(\rho(x,y))$}
\begin{thm}\label{adlem1}
For fixed $y\in \mathbb{B}^n$, there holds
\begin{equation*}\begin{split}
G_{-\mu \mathrm{i}}(\rho(x,y))&=\frac{A_{n,-\mu \mathrm{i}}}{2^{\frac{n-1}{2}}}\frac{\left(\cosh(\frac{\rho(x)}{2})\right)^{2\mu \mathrm{i}}}{\cosh^{n-1}(\frac{\rho(x)}{2})}\frac{\left(\cosh(\frac{\rho(y)}{2})\right)^{2\mu \mathrm{i}}}{\cosh^{n-1}(\frac{\rho(y)}{2})}\\
&\ \ \times   \left(1-2\hat{x}\cdot y+|y|^2\right)^{\mu \mathrm{i}-\frac{n-1}{2}}\int_{0}^{\pi}(\sin t)^{-2\mu \mathrm{i}}dt+ O\left(\frac{1}{\sinh^{\frac{n+1}{2}}(\rho(x))}\right)
\end{split}\end{equation*}
as $\rho(x)\rightarrow +\infty$.
\end{thm}
\begin{proof}
Careful calculation gives that
\begin{equation}\begin{split}
G_{-\mu \mathrm{i}}(\rho(x,y))&=A_{n,-\mu \mathrm{i}} \frac{\left(2\sinh^2(\frac{\rho(x,y)}{2})+1\right)^{\frac{n-3}{2}+\mu \mathrm{i}}}{\left(2\sinh(\frac{\rho(x,y)}{2})\cosh(\frac{\rho(x,y)}{2})\right)^{n-2}}\int_{0}^{\pi}\left(1+\frac{\cos t}{\cosh \rho} \right)^{\frac{n-3}{2}+\mu \mathrm{i}}(\sin t)^{-2\mu \mathrm{i}}dt\\
&=\frac{A_{n,-\mu \mathrm{i}}}{2^{\frac{n-1}{2}}}\frac{\left(\sinh(\frac{\rho(x,y)}{2})\right)^{2\mu \mathrm{i}}}{\sinh(\frac{\rho(x,y)}{2})\left(\cosh(\frac{\rho(x,y)}{2})\right)^{n-2}}\int_{0}^{\pi}\left(1+\frac{\cos t}{\cosh \rho} \right)^{\frac{n-3}{2}+\mu \mathrm{i}}(\sin t)^{-2\mu \mathrm{i}}dt\\
&\ \ \ +O\left(\frac{1}{\sinh^{n-2}(\rho(x))}\right)\\
&=\frac{A_{n,-\mu \mathrm{i}}}{2^{\frac{n-1}{2}}}\frac{\left(\sinh(\frac{\rho(x,y)}{2})\right)^{2\mu \mathrm{i}}}{\sinh(\frac{\rho(x,y)}{2})\left(\cosh(\frac{\rho(x,y)}{2})\right)^{n-2}}\int_{0}^{\pi}(\sin t)^{-2\mu \mathrm{i}}dt\\
&\ \ \ +O\left(\frac{1}{\sinh^{\frac{n-1}{2}+1}(\rho(x))}\right).
\end{split}\end{equation}
From \cite{Ahlfors,Hua}, we know that
$$\sinh(\frac{\rho(x,y)}{2})=\frac{|x-y|}{\sqrt{(1-|x|^2)(1-|y|^2)}}=|x-y|\cosh\left(\frac{\rho(x)}{2}\right)\cosh\left(\frac{\rho(y)}{2}\right),$$
and $$\cosh\left(\frac{\rho(x,y)}{2}\right)=\frac{\sqrt{1-2x\cdot y+|x|^2|y|^2}}{\sqrt{(1-|x|^2)(1-|y|^2)}}=\sqrt{1-2x\cdot y+|x|^2|y|^2}\cosh\left(\frac{\rho(x)}{2}\right)\cosh\left(\frac{\rho(y)}{2}\right).$$
Hence we can furthermore write
\begin{equation}\label{expansion}\begin{split}
G_{-\mu \mathrm{i}}(\rho(x,y))&=\frac{A_{n,-\mu \mathrm{i}}}{2^{\frac{n-1}{2}}}\frac{\left(\cosh(\frac{\rho(x)}{2})\right)^{2\mu i}}{\cosh^{n-1}(\frac{\rho(x)}{2})}\frac{\left(\cosh(\frac{\rho(y)}{2})\right)^{2\mu \mathrm{i}}}{\cosh^{n-1}(\frac{\rho(y)}{2})}\\
&\ \ \times   \frac{|x-y|^{2\mu \mathrm{i}}}{|x-y|(\sqrt{1-2x\cdot y+|x|^2|y|^2})^{n-2}} \int_{0}^{\pi}(\sin t)^{-2\mu \mathrm{i}}dt+ O\left(\frac{1}{\sinh^{\frac{n+1}{2}}(\rho(x))}\right)\\
&=\frac{A_{n,-\mu \mathrm{i}}}{2^{\frac{n-1}{2}}}\frac{\left(\cosh(\frac{\rho(x)}{2})\right)^{2\mu \mathrm{i}}}{\cosh^{n-1}(\frac{\rho(x)}{2})}\frac{\left(\cosh(\frac{\rho(y)}{2})\right)^{2\mu\mathrm{i}}}{\cosh^{n-1}(\frac{\rho(y)}{2})}\\
&\ \ \times   \left(1-2\hat{x}\cdot y+|y|^2\right)^{\mu\mathrm{i}-\frac{n-1}{2}}\int_{0}^{\pi}(\sin t)^{-2\mu \mathrm{i}}dt+ O\left(\frac{1}{\sinh^{\frac{n+1}{2}}(\rho(x))}\right)\\
\end{split}\end{equation}
as $\rho(x)\rightarrow +\infty$ for any fixed $y\in \mathbb{B}^n$, where $\hat{x}=\frac{x}{|x|}$. Then the proof of Theorem \ref{adlem1} is accomplished.
\end{proof}

\subsection{Radiation condition and far-field pattern}
For fixed $y\in \mathbb{B}^n$, through the polar coordinates, we can also see the $G(\rho(x,y))$ as the function of the variable $\rho$ and $\hat{x}$. According to the asymptotic formula \eqref{expansion}, we can write
\begin{equation}\begin{split}
G_{-\mu \mathrm{i}}(\rho(x,y))=G(\rho,\hat{x}, y)&=\frac{A_{n,-\mu i}}{2^{\frac{n-1}{2}}}\frac{\left(\cosh(\frac{\rho}{2})\right)^{2\mu \mathrm{i}}}{\cosh^{n-1}(\frac{\rho}{2})}\frac{\left(\cosh(\frac{\rho(y)}{2})\right)^{2\mu \mathrm{i}}}{\cosh^{n-1}(\frac{\rho(y)}{2})}\\
&\ \ \times   \left(1-2\hat{x}\cdot y+|y|^2\right)^{\mu\mathrm{i}-\frac{n-1}{2}}\int_{0}^{\pi}(\sin t)^{-2\mu \mathrm{i}}dt+ O\left(\frac{1}{\sinh^{\frac{n+1}{2}}(\rho)}\right)
\end{split}\end{equation}
when $\rho\rightarrow +\infty$. Direct computation gives that
\begin{equation*}\begin{split}
\frac{\partial G}{\partial \rho}&=\frac{1}{2}\frac{A_{n,-\mu \mathrm{i}}}{2^{\frac{n-1}{2}}}\left(2\mu \mathrm{i}-(n-1)\right)\left(\cosh(\frac{\rho}{2})\right)^{2\mu \mathrm{i}-n}\sinh(\frac{\rho}{2})\frac{\left(\cosh(\frac{\rho(y)}{2})\right)^{2\mu \mathrm{i}}}{\cosh^{n-1}(\frac{\rho(y)}{2})}\\
&\ \ \times   \left(1-2\hat{x}\cdot y+|y|^2\right)^{\mu i-\frac{n-1}{2}}\int_{0}^{\pi}(\sin t)^{-2\mu \mathrm{i}}dt+ O\left(\frac{1}{\sinh^{\frac{n+1}{2}}(\rho)}\right)\\
&=(\mu \mathrm{i}-\frac{n-1}{2})\tanh(\frac{\rho}{2})G+O\left(\frac{1}{\sinh^{\frac{n+1}{2}}(\rho)}\right)
\end{split}\end{equation*}
as $\rho\rightarrow +\infty$. Hence the Green function $G_{-\mu \mathrm{i}}(\rho,\hat{x}, y)$ satisfies the condition
\begin{equation}\label{radiation}
\frac{\partial G}{\partial \rho}-(\mu \mathrm{i}-\frac{n-1}{2})\tanh(\frac{\rho}{2})G=O\left(\frac{1}{\sinh^{\frac{n+1}{2}}(\rho)}\right)
\end{equation}
as $\rho\rightarrow +\infty$. We called the condition \eqref{radiation} \textbf{hyperbolic radiation condition} and prove that this radiation condition will guarantee the existence and uniqueness for solutions of the wave equation involving the source or potential later.
\begin{rmk}
If we consider the ingoing wave, the Green function $G_{\mu \mathrm{i}}$ satisfies inner radiation condition
$$\frac{\partial G}{\partial \rho}+(\mu \mathrm{i}+\frac{n-1}{2})\tanh(\frac{\rho}{2})G=O\left(\frac{1}{\sinh^{\frac{n+1}{2}}(\rho)}\right).$$
One can similar prove the existence and uniqueness for the wave equation involving the source or potential under the inner radiation condition.
\end{rmk}
\vskip0.1cm
\subsection{Scattering theory for Helmholtz equation with a source}
Consider the Helmholtz equation
\begin{equation}\label{wave-f}
\left(-\Delta_{\mathbb{B}^n}-\frac{(n-1)^2}{4}-\mu^2\right)u=f(x),\ \  x\in \mathbb{B}^n.
\end{equation} with the source $f\in C^{\infty}_{c}(\mathbb{B}^n)$.

Obviously, the above equation has a solution  $$u(x)=\int_{\mathbb{B}^n}G_{-\mu \mathrm{i}}(\rho(x,y))f(y)dV_{\mathbb{B}^n}(y)$$
since $$(-\Delta_{\mathbb{B}^n})_yG_{-\mu i}(\rho(x,y))-\frac{(n-1)^2}{4}G_{-\mu \mathrm{i}}(\rho(x,y))-\mu^2G_{-\mu i}(\rho(x,y))=\delta_{x}(y).$$
According to the accurate asymptotic behavior of $G_{-\mu \mathrm{i}}(\rho(x,y))$ when $\rho(x)\rightarrow +\infty$ from \eqref{expansion}, we derive that
\begin{equation}\begin{split}
u(x)&=\frac{\left(\cosh(\frac{\rho(x)}{2})\right)^{2\mu \mathrm{i}}}{\cosh^{n-1}(\frac{\rho(x)}{2})}\int_{0}^{\pi}(\sin t)^{-2\mu \mathrm{i}}dt\\
&\ \ \times \int_{\mathbb{B}^n}\frac{\left(\cosh(\frac{\rho(y)}{2})\right)^{2\mu \mathrm{i}}}{\cosh^{n-1}(\frac{\rho(y)}{2})}\left(1-2\hat{x}\cdot y+|y|^2\right)^{\mu \mathrm{i}-\frac{n-1}{2}}f(y)dV_{\mathbb{B}^n}(y)+O\left(\frac{1}{\sinh^{\frac{n+1}{2}}(\rho(x))}\right).
\end{split}\end{equation}
The \textbf{far-field pattern} of $u$ is defined through
\begin{equation}\begin{split}
u^{\mu}_{\infty}(\hat{x})&=\frac{A_{n,-\mu \mathrm{i}}}{2^{\frac{n-1}{2}}}\int_{0}^{\pi}(\sin t)^{-2\mu \mathrm{i}}dt\int_{\mathbb{B}^n}\frac{\left(\cosh(\frac{\rho(y)}{2})\right)^{2\mu \mathrm{i}}}{\cosh^{n-1}(\frac{\rho(y)}{2})}\left(1-2\hat{x}\cdot y+|y|^2\right)^{\mu i-\frac{n-1}{2}}f(y)dV_{\mathbb{B}^n}(y).\\
\end{split}\end{equation}
Careful computation gives
\begin{equation*}\begin{split}
u^{\mu}_{\infty}(\hat{x})&=\frac{A_{n,-\mu \mathrm{i}}}{2^{\frac{n-1}{2}}}\int_{0}^{\pi}(\sin t)^{-2\mu \mathrm{i}}dt\int_{\mathbb{B}^n}e_{-2\mu, \hat{x}}(y)f(y)dV_{\mathbb{B}^n}(y)\\
&=\frac{A_{n,-\mu \mathrm{i}}}{2^{\frac{n-1}{2}}}\int_{0}^{\pi}(\sin t)^{-2\mu \mathrm{i}}dt\widehat{f}(2\mu,\hat{x}),
\end{split}\end{equation*}
which gives $$\widehat{f}(\hat{x},\mu)=\left(\frac{A_{n,-\mu \mathrm{i}}}{2^{\frac{n-1}{2}}}\right)^{-1}\left(\int_{0}^{\pi}(\sin t)^{-2\mu \mathrm{i}}dt\right)^{-1}u_{\infty}^{\mu/2}(\hat{x}).$$
Hence, through inverse Fourier transform, we derive that
$$f(y)=D_{n}\left(\frac{A_{n,-\mu \mathrm{i}}}{2^{\frac{n-1}{2}}}\right)^{-1}\left(\int_{0}^{\pi}(\sin t)^{-2\mu \mathrm{i}}dt\right)^{-1}\int_{-\infty}^{\infty}\int_{\mathbb{S}^{n-1}}u_{\infty}^{\mu/2}(\hat{x})e_{\mu,\hat{x}}(y)|c(\mu)|^{-2}d\mu d\sigma.$$
That is to say that if we know the information for far-field pattern $u_{\infty}(\hat{x}, \frac{\mu}{2})$ for all observable direction $\hat{x}$ and frequency $\mu$,
then the source $f(x)$ can be uniquely determined by inverse Fourier transform in hyperbolic space. This is one of the rare cases that unknown source can be expressed explicitly.

\vskip0.1cm

Next we will prove that if the solution of Helmholtz equation \eqref{wave-f} satisfies the hyperbolic radiation condition, then the solution is unique and hence $$u(x)=\int_{\mathbb{B}^n}G_{-\mu \mathrm{i}}(\rho(x,y))f(y)dV_{\mathbb{B}^n}(y).$$

\begin{thm}\label{thm2}
Assume that $u\in C^{2}(\mathbb{B}^n)$ satisfies the Helmholtz equation with the source
\begin{equation}
\left(-\Delta_{\mathbb{B}^n}-\frac{(n-1)^2}{4}-\mu^2\right)u=f(x),\ \  x\in \mathbb{B}^n
\end{equation}
and
the hyperbolic radiation condition $$\frac{\partial u}{\partial \rho}-\left(\mu i-\frac{n-1}{2}\right)\tanh(\frac{\rho}{2})u=o\left(\frac{1}{\sinh^{\frac{n-1}{2}}(\rho)}\right),$$
where $f\in C^{\infty}_c(\mathbb{B}^n)$.
Then the solution $u$ is unique and can be written as
\begin{equation*}\begin{split}
&u(x)=\int_{\mathbb{B}^n}G_{-\mu \mathrm{i}}(\rho(x,y))f(y)dV_{\mathbb{B}^n}(y)\\
&=\frac{\left(\cosh(\frac{\rho(x)}{2})\right)^{2\mu \mathrm{i}}}{\cosh^{n-1}(\frac{\rho(x)}{2})}u_{\infty}(\hat{x}, \mu)+O\left(\frac{1}{\sinh^{\frac{n+1}{2}}(\rho(x))}\right)£¬
\end{split}\end{equation*}
where
$$u_{\infty}(\hat{x}, \mu)=\frac{A_{n,-\mu \mathrm{i}}}{2^{\frac{n-1}{2}}}\int_{0}^{\pi}(\sin t)^{-2\mu \mathrm{i}}dt\int_{\mathbb{B}^n}e_{-2\mu, \hat{x}}(y)f(y)dV_{\mathbb{B}^n}(y).$$
Especially, the solution of Helmholtz equation $$\left(-\Delta_{\mathbb{B}^n}-\frac{(n-1)^2}{4}-\mu^2\right)u=0,\ \  x\in \mathbb{B}^n$$ satisfying the hyperbolic radiation condition must be equal to zero in $\mathbb{B}^n$.

\end{thm}

\begin{proof}
We first claim that the radiation condition implies that
$$\lim\limits_{R\rightarrow +\infty}\sinh^{n-1}(R)\int_{\partial B^n(0,\tanh(\frac{R}{2}))}|u|^2d\sigma=O(1),$$
which is equivalent to $$\lim\limits_{R\rightarrow +\infty}\int_{\partial B_\mathbb{B}^n(0,R)}|u|^2d\sigma_{\mathbb{B}^n}=O(1).$$
From the hyperbolic radiation condition, we deduce that
\begin{equation}\label{adintegral}\begin{split}
&\lim\limits_{R\rightarrow +\infty}\int_{\partial B_{\mathbb{B}^n}(0,R)}\bigg|\frac{\partial u}{\partial \rho}-\left(\mu i-\frac{n-1}{2}\right)\tanh(\frac{\rho}{2})u\bigg|^2d\sigma_{\mathbb{B}^n}\\
&\ \  =\lim\limits_{R\rightarrow +\infty}\int_{\partial B^n(0,\tanh(\frac{R}{2}))}\left|\frac{\partial u}{\partial \rho}-\left(\mu i-\frac{n-1}{2}\right)\tanh(\frac{\rho}{2})u\right|^2\left(\frac{2}{1-|x|^2}\right)^{n-1}d\sigma\\
&\ \ =\lim\limits_{R\rightarrow +\infty}\left(2\cosh(\frac{R}{2})\right)^{n-1}\int_{\partial B^n(0,\tanh(\frac{R}{2}))}o\left(\frac{1}{\sinh^{n-1}(\rho)}\right)d\sigma=0.
\end{split}\end{equation}
On the other hand, careful calculation gives that
\begin{equation}\label{expansion1}\begin{split}
&\int_{\partial B_{\mathbb{B}^n}(0,R)}\left|\frac{\partial u}{\partial \rho}-\left(\mu i-\frac{n-1}{2}\right)\tanh(\frac{\rho}{2})u\right|^2d\sigma_{\mathbb{B}^n}\\
&\ \ =\int_{\partial B_{\mathbb{B}^n}(0,R)}\left|\frac{\partial u}{\partial\rho}\right|^2d\sigma_{\mathbb{B}^n}+\int_{\partial B_{\mathbb{B}^n}(0,R)}\left(\frac{(n-1)^2}{4}+\mu^2\right)\tanh^2(\frac{\rho}{2})u^2d\sigma_{\mathbb{B}^n}\\
&\ \ \ \ +2\Im\left(\int_{\partial B_{\mathbb{B}^n}(0,R)}\tanh(\frac{\rho}{2})\frac{\partial \overline{u}}{\partial\rho}ud\sigma_{\mathbb{B}^n}\right)+(n-1)\Re\left(\int_{\partial B_{\mathbb{B}^n}(0,R)}\tanh(\frac{\rho}{2})\frac{\partial \overline{u}}{\partial\rho}ud\sigma_{\mathbb{B}^n}\right).
\end{split}\end{equation}
Since $$\int_{\partial B_{\mathbb{B}^n}(0,R)}\left(\left|\frac{\partial u}{\partial\rho}\right|^2+\frac{(n-1)^2}{4}\tanh^2(\frac{\rho}{2})u^2\right)d\sigma_{\mathbb{B}^n}\geq (n-1)\Re\left(\int_{\partial B_{\mathbb{B}^n}(0,R)}\tanh(\frac{\rho}{2})\frac{\partial \overline{u}}{\partial\rho}ud\sigma_{\mathbb{B}^n}\right),$$
then we deduce from \eqref{adintegral} and \eqref{expansion1} that
\begin{equation}\label{eq1}\begin{split}
o(1)\geq \int_{\partial B_{\mathbb{B}^n}(0,R)}\mu^2\tanh^2(\frac{\rho}{2})u^2d\sigma_{\mathbb{B}^n}+2\Im\left(\int_{\partial B_{\mathbb{B}^n}(0,R)}\tanh(\frac{\rho}{2})\frac{\partial \overline{u}}{\partial\rho}ud\sigma_{\mathbb{B}^n}\right).
\end{split}\end{equation} when $R\rightarrow +\infty$.
Next, we start to calculate the $2\Im\left(\int_{\partial B_{\mathbb{B}^n}(0,R)}\tanh(\frac{\rho}{2})\frac{\partial \overline{u}}{\partial\rho}ud\sigma_{\mathbb{B}^n}\right)$.
Since $$\int_{\partial B_{\mathbb{B}^n}(0,R)}\tanh(\frac{\rho}{2})\frac{\partial \overline{u}}{\partial\rho}ud\sigma_{\mathbb{B}^n}=2^{n-1}\cosh^{2n-2}(\frac{R}{2})\tanh(\frac{R}{2})\int_{\partial B^{n}(0, \tanh(\frac{R}{2}))}\frac{\partial \overline{u}}{\partial\rho}ud\sigma$$ and $u$ satisfies equation $$(-\Delta_{\mathbb{B}^n}-\frac{(n-1)^2}{4}-\mu^2)u=f(x),\ \  x\in \mathbb{B}^n.$$
 Applying the Green formula \eqref{lem:green-ball} in $B_{\mathbb{B}^n}(0,R)\setminus  B_{\mathbb{B}^n}(0,R_0)$, we derive that
\begin{equation}\label{eq2}\begin{split}
&\int_{\partial B_{\mathbb{B}^n}(0,R)}\tanh(\frac{\rho}{2})\frac{\partial \overline{u}}{\partial\rho}ud\sigma_{\mathbb{B}^n}\\
&\ \ =\tanh(\frac{R}{2})\int_{\partial B_{\mathbb{B}^n}(0,R_0)}\frac{\partial \overline{u}}{\partial\rho}ud\sigma_{\mathbb{B}^n}+\tanh(\frac{R}{2})\int_{B_{\mathbb{B}^n}(0,R)\setminus  B_{\mathbb{B}^n}(0,R_0)}\left(\Delta_{\mathbb{B}^n}(\overline{u})u+(\nabla_{\mathbb{B}^n}\overline{u},\nabla_{\mathbb{B}^n}u)_g\right)dV_{\mathbb{B}^n}\\
&\ \ =\tanh(\frac{R}{2})\int_{\partial B_{\mathbb{B}^n}(0,R_0)}\frac{\partial \overline{u}}{\partial\rho}ud\sigma_{\mathbb{B}^n}+\tanh(\frac{R}{2})\int_{B_{\mathbb{B}^n}(0,R)\setminus  B_{\mathbb{B}^n}(0,R_0)}\left(|\nabla_{\mathbb{B}^n}u|^2-\left(\frac{(n-1)^2}{4}+\mu^2\right)|u|^2\right)dV_{\mathbb{B}^n}\\
&\ \ \ \ -\tanh(\frac{R}{2})\int_{supp(f)}\overline{f}udV_{\mathbb{B}^n}.
\end{split}\end{equation}
Hence $$\Im\left(\int_{\partial B_{\mathbb{B}^n}(0,R)}\tanh(\frac{\rho}{2})\frac{\partial \overline{u}}{\partial\rho}ud\sigma_{\mathbb{B}^n}\right)=\Im\left(\tanh(\frac{R}{2})\int_{\partial B_{\mathbb{B}^n}(0,R_0)}\frac{\partial \overline{u}}{\partial\rho}ud\sigma_{\mathbb{B}^n}\right)-\Im\left(\tanh(\frac{R}{2})\int_{supp(f)}\overline{f}udV_{\mathbb{B}^n}\right).$$
This together with \eqref{eq1} and \eqref{eq2} yields that
\begin{equation}\label{eq3}\begin{split}
o(1)\geq \int_{\partial B_{\mathbb{B}^n}(0,R)}\mu^2\tanh^2(\frac{\rho}{2})u^2d\sigma_{\mathbb{B}^n}+\Im\left(\tanh(\frac{R}{2})\int_{\partial B_{\mathbb{B}^n}(0,R_0)}\frac{\partial \overline{u}}{\partial\rho}ud\sigma_{\mathbb{B}^n}\right)-\Im\left(\tanh(\frac{R}{2})\int_{supp(f)}\overline{f}udV_{\mathbb{B}^n}\right).
\end{split}\end{equation}
Let $R\rightarrow +\infty$, we conclude that
$$2^{n-1}\cosh^{2n-2}(\frac{R}{2})\int_{\partial B^{n}(0, \tanh(\frac{R}{2}))}|u|^2d\sigma=\int_{\partial B_\mathbb{B}^n(0,R)}|u|^2d\sigma_{g}=O(1).$$
Now, we are in position to prove that $$u(x)=\int_{\mathbb{B}^n}G_{-\mu i}(\rho(x,y))f(y)dV_{\mathbb{B}^n}(y).$$
We only need to prove that $$u(0)=\int_{\mathbb{B}^n}G_{-\mu i}(\rho(y))f(y)dV_{\mathbb{B}^n}(y).$$
In fact, define $v(y)=(u\circ\tau_{x})(y)$, through the commutative property of the operator $-\Delta_{\mathbb{B}^n}$ and $\tau_{x}$, then $v(y)$ satisfies equation
\begin{equation}
(-\Delta_{\mathbb{B}^n}-\frac{(n-1)^2}{4}-\mu^2)v=f(\tau_x(y)),\ \ y\in \mathbb{B}^n.
\end{equation}
Hence \begin{equation}\begin{split}
u(x)=v(0)&=\int_{\mathbb{B}^n}G_{-\mu \mathrm{i}}(\rho(y))f(\tau_x(y))dV_{\mathbb{B}^n}(y)\\
&=\int_{\mathbb{B}^n}G_{-\mu \mathrm{i}}(\rho(x,y))f(y)dV_{\mathbb{B}^n}(y).
\end{split}\end{equation}
Now, we start to prove that $$u(0)=\int_{\mathbb{B}^n}G_{-\mu \mathrm{i}}(\rho(y))f(y)dV_{\mathbb{B}^n}(y).$$

According to Lemma \ref{lem:green-ball}, we can write
\small{\begin{equation}\label{eq4}\begin{split}
&\int_{B_{\mathbb{B}^n}(0,R)}(-\Delta_{\mathbb{B}^n}-\frac{(n-1)^2}{4}-\mu^2)\left(u\right)\cdot G_{-\mu \mathrm{i}}(\rho(y))dV_{\mathbb{B}^n}(y)\\
&\ \ =\int_{B_{\mathbb{B}^n}(0,R)}(-\Delta_{\mathbb{B}^n}-\frac{(n-1)^2}{4}-\mu^2)\left(G_{-\mu \mathrm{i}}(\rho)\right)\cdot u(y)dV_{\mathbb{B}^n}(y)+\int_{B_{\mathbb{B}^n}(0,R)}\left(\frac{\partial u}{\partial \rho}G_{-\mu\mathrm{i}}-\frac{\partial G_{-\mu\mathrm{i}}}{\partial \rho}u\right)d\sigma_{\mathbb{B}^n}\\
&=\int_{B_{\mathbb{B}^n}(0,R)}(-\Delta_{\mathbb{B}^n}-\frac{(n-1)^2}{4}-\mu^2)(G_{-\mu\mathrm{i}}(\rho))u(y)dV_{\mathbb{B}^n}(y)+\int_{\partial B_{\mathbb{B}^n}(0,R)}\left(\frac{\partial u}{\partial \rho}-\left(\mu\mathrm{i}-\frac{n-1}{2}\right)\tanh(\frac{\rho}{2})u\right)G_{-\mu \mathrm{i}}d\sigma_{\mathbb{B}^n}\\
&\ \ -\int_{\partial B_{\mathbb{B}^n}(0,R)}\left(\frac{\partial G_{-\mu \mathrm{i}}}{\partial \rho}-\left(\mu i-\frac{n-1}{2}\right)\tanh(\frac{\rho}{2})G_{-\mu \mathrm{i}}\right)ud\sigma_{\mathbb{B}^n}\\
&=u(0)+\int_{\partial B_{\mathbb{B}^n}(0,R)}\left(\frac{\partial u}{\partial \rho}-\left(\mu i-\frac{n-1}{2}\right)\tanh(\frac{\rho}{2})u\right)G_{-\mu\mathrm{i}}d\sigma_{\mathbb{B}^n}\\
&\ \ -\int_{\partial B_{\mathbb{B}^n}(0,R)}\left(\frac{\partial G_{-\mu \mathrm{i}}}{\partial \rho}-\left(\mu i-\frac{n-1}{2}\right)\tanh(\frac{\rho}{2})G_{-\mu \mathrm{i}}\right)ud\sigma_{\mathbb{B}^n}
\end{split}\end{equation}}
Since $u$ and $G_{-\mu\mathrm{i}}$ satisfies the radiation condition
$$\frac{\partial u}{\partial \rho}-\left(\mu i-\frac{n-1}{2}\right)\tanh(\frac{\rho}{2})u=o\left(\frac{1}{\sinh^{\frac{n-1}{2}}(\rho)}\right),$$  $$\frac{\partial G_{-\mu \mathrm{i}}}{\partial \rho}-\left(\mu \mathrm{i}-\frac{n-1}{2}\right)\tanh(\frac{\rho}{2})G_{-\mu i}=o\left(\frac{1}{\sinh^{\frac{n-1}{2}}(\rho)}\right),$$
and $\int_{\partial B_{\mathbb{B}^n}(0,R)}|u|^2d\sigma_{\mathbb{B}^n}$ is bounded,
hence \small{\begin{equation}\begin{split}
&\lim\limits_{R\rightarrow +\infty}\int_{\partial B_{\mathbb{B}^n}(0,R)}\left(\left(\frac{\partial u}{\partial \rho}-(\mu \mathrm{i}-\frac{n-1}{2})\tanh(\frac{\rho}{2})u\right)G_{-\mu \mathrm{i}}-\left(\frac{\partial G_{-\mu \mathrm{i}}}{\partial \rho}-(\mu \mathrm{i}-\frac{n-1}{2})\tanh(\frac{\rho}{2})G_{-\mu \mathrm{i}}\right)u\right)d\sigma_{\mathbb{B}^n}\\
&\ \ \leq \lim\limits_{R\rightarrow +\infty}\left(\int_{\partial B_{\mathbb{B}^n}(0,R)}\left|\left(\frac{\partial u}{\partial \rho}-(\mu \mathrm{i}-\frac{n-1}{2})\tanh(\frac{\rho}{2})\right)u\right|^2d\sigma_{\mathbb{B}^n}\right)^{\frac{1}{2}}\left(\int_{\partial B_{\mathbb{B}^n}(0,R)}|G_{-\mu \mathrm{i}}|^2d\sigma_{\mathbb{B}^n}\right)^{\frac{1}{2}}\\
&\ \ \ \ +\lim\limits_{R\rightarrow +\infty}\left(\int_{\partial B_{\mathbb{B}^n}(0,R)}\left|\left(\frac{\partial G_{-\mu \mathrm{i}}}{\partial \rho}-(\mu \mathrm{i}-\frac{n-1}{2})\tanh(\frac{\rho}{2})\right)G_{-\mu \mathrm{i}}\right|^2d\sigma_{\mathbb{B}^n}\right)^{\frac{1}{2}}\left(\int_{\partial B_{\mathbb{B}^n}(0,R)}|u|^2d\sigma_{\mathbb{B}^n}\right)^{\frac{1}{2}}\\
&\ \ =0.
\end{split}\end{equation}}
This together with \eqref{eq4} implies that
$$u(0)=\lim\limits_{R\rightarrow +\infty}\int_{B_{\mathbb{B}^n}(0,R)}(-\Delta_{\mathbb{B}^n}-\frac{(n-1)^2}{4}-\mu^2)(u)G(\rho)dV_{\mathbb{B}^n}(y)=\int_{\mathbb{B}^n}G_{-\mu i}(\rho(y))f(y)dV_{\mathbb{B}^n}(y).$$
\end{proof}
\subsection{Rellich Theorem in hyperbolic space}
\begin{thm}\label{thm3}(Rellich Theorem in hyperbolic space)
Let $\Omega$ be a bounded domain in hyperbolic space $\mathbb{B}^n$. Assume that $u\in C^{2}(\mathbb{B}^n\setminus \Omega)$ satisfies the Helmholtz equation in hyperbolic space
\begin{equation}
-\Delta_{\mathbb{B}^n}u-\frac{(n-1)^2}{4}u-\mu^2 u=0
\end{equation}
in $\mathbb{B}^n\setminus \Omega$. If $\lim\limits_{R\rightarrow +\infty}\int_{\partial B_{\mathbb{B}^n}(0,R)}|u|^2d\sigma_{\mathbb{B}^n}=0$, then $u$ is equal to zero in $\mathbb{B}^n\setminus \Omega$.
\end{thm}

\begin{proof}
The proof is based on spherical harmonic expansion. We first recall some basic facts about spherical harmonics and a detailed discussion on spherical harmonics can be found in, e.g., \cite{M} and \cite{SteinWeiss}. In brief, spherical harmonics are restrictions to the
 unit sphere $\mathbb{S}^n$ of polynomial $Y(x)$, where $Y(x)$ satisfies $\Delta_{\mathbb{R}^n}Y=0$ and $\Delta_{\mathbb{R}^n}$ is the Laplacian operator in $\mathbb{R}^{n+1}$. The space of all spherical harmonics of degree $l$ on $\mathbb{S}^n$, denoted by $\mathcal{H}_{l}^n$, has an orthonormal basis
$$\{Y_{l,m}:m=1,\cdots,N(n,l),\ l=1, 2, \cdot\cdot\cdot \},$$
where
$$N(n,0)=1~~ \text{and}~~ N(n,l)=\frac{(2l+n-1)(\Gamma(l+n-1))}{\Gamma(l+1)\Gamma(n)}.$$
Furthermore, $Y_{l,m}$ is in fact the spectrum of $\Delta_{\mathbb{S}^n}$ and satisfies equation $-\Delta_{\mathbb{S}^n}Y_{m,l}(x')=l(n+l-1)$.
For $u\in{L^2(\mathbb{S}^n)}$, it can be expanded in terms of spherical harmonics
$$u=\sum_{l=0}^{+\infty}\sum_{m=1}^{N(n,l)}u_{m,l}Y_{m,l}, ~~\text{where}~~u_{m,l}=\int_{\mathbb{S}^n}uY_{m,l}d\omega,$$
and $d\omega$ denotes the surface measure in $\mathbb{S}^n$.
By Parseval's identity, we have
$$\|u\|^2_{L^2(\mathbb{S}^n)}=\int_{\mathbb{S}^n}|u|^2d\omega=\sum_{l=0}^{+\infty}\sum_{m=1}^{N(n,l)}|u_{m,l}|^2.$$
\vskip0.1cm

Now, we are in position to give the proof of hyperbolic Rellich theorem. Let $B_{\mathbb{B}^n}(0,R_0)$ denote the hyperbolic ball centered at origin with the radius equal to $R_0$. Since $\Omega$ is a bounded domain, one choose sufficiently large $R_0$ such that $\Omega\subseteq B_{\mathbb{B}^n}(0,R_0)$. Then for any $x\in \mathbb{B}^n\setminus B_{\mathbb{B}^n}(0,R_0)$, $u(\tanh \frac{\rho}{2}x')$ is a continuous function defined on the sphere $S^{n-1}$, where $|x|=\tanh(\frac{\rho}{2})$ and $x'\in \mathbb{S}^{n-1}$. According to the spherical harmonic expansion formula, we can write
$$u(\tanh \frac{\rho}{2}x')=\sum_{l=1}^{+\infty}\sum_{m=1}^{N(l,n)}a_{m,l}(\rho)Y_{m,l}(x'),$$
where $Y_{m,l}(x')$ is the normalized orthogonal basis of $L^{2}(\mathbb{S}^{n-1})$ and $a_{m,l}(\rho)=\int_{\mathbb{S}^{n-1}}u(\tanh \frac{\rho}{2}x')Y_{m,l}(x')d\omega$. Since the Laplace operator $-\Delta_{\mathbb{B}^n}$ and gradient operator $\nabla_{\mathbb{B}^n}$ on the hyperbolic space can be written as $$\Delta_{\mathbb{B}^n}=\frac{\partial^2}{\partial \rho^2}+(n-1)\coth\rho\frac{\partial}{\partial \rho}+\frac{1}{\sinh^2 \rho}\Delta_{\mathbb{S}^{n-1}},\ \ \nabla_{\mathbb{B}^n}=(\frac{\partial}{\partial \rho},\  \frac{1}{\sinh \rho}\nabla_{\mathbb{S}^{n-1}}),$$
then it follows that for $x\in \mathbb{B}^n\setminus B_{\mathbb{B}^n}(0,R_0)$,
\begin{equation}\begin{split}
\Delta_{\mathbb{B}^n}(u)&=\sum_{l=1}^{+\infty}\sum_{m=1}^{N(l,n)}\Delta_{\mathbb{B}^n}(a_{m,l}(\rho)Y_{m,l}(x'))\\
&=\sum_{l=1}^{+\infty}\sum_{m=1}^{N(l,n)}\left((\Delta_{\mathbb{B}^n}a_{m,l}(\rho))Y_{m,l}(x')+a_{m,l}(\Delta_{\mathbb{B}^n}Y_{m,l}(x'))-\nabla_{\mathbb{B}^n}(a_{m,l}(\rho))\nabla_{\mathbb{B}^n}(Y_{m,l}(x'))\right)\\
&=\sum_{l=1}^{+\infty}\sum_{m=1}^{N(l,n)}\left((\frac{\partial^2}{\partial \rho^2}+(n-1)\coth \rho \frac{\partial}{\partial \rho})a_{m,l}(\rho)\times Y_{m,l}(x')-\frac{a_{m,l}(\rho)}{\sinh^2 \rho}l(l+n-2)Y_{m,l}(x')\right)\\
&=\sum_{l=1}^{+\infty}\sum_{m=1}^{N(l,n)}\left(\frac{\partial^2}{\partial \rho^2}+(n-1)\coth \rho \frac{\partial}{\partial \rho}-l(l+n-2)\frac{1}{\sinh^2 \rho}\right)a_{m,l}(\rho).
\end{split}\end{equation}
Combining this and $u$ satisfying equation $-\Delta_{\mathbb{B}^n}u-\frac{(n-1)^2}{4}u-\mu^2 u=0$ in $\mathbb{B}^n\setminus B_{\mathbb{B}^n}(0,R_0)$, we derive that $a_{m,l}(\rho)$ satisfies the following ordinary equation
$$\left(\frac{\partial^2}{\partial \rho^2}+(n-1)\coth \rho \frac{\partial}{\partial \rho}+\frac{(n-1)^2}{4}+\mu^2-l(l+n-2)\frac{1}{\sinh^2 \rho}\right)a_{m,l}(\rho)=0.$$
When $\rho$ is sufficiently large, this equation can be compared with the following equation
$$\left(\frac{\partial^2}{\partial \rho^2}+(n-1)\coth \rho \frac{\partial}{\partial \rho}+\frac{(n-1)^2}{4}+\mu^2\right)a_{m,l}(\rho)=0,$$
which is a standard Jacobi equation. The general Jacobi equation (see \cite{liu}) is the following second-order ordinary equation
\begin{equation}\label{jocabi}
\left(\frac{\partial^2}{\partial \rho^2}+(2\alpha+1)\coth \rho \frac{\partial}{\partial \rho}+(2\beta+1)\tanh\rho\frac{\partial}{\partial \rho}+\lambda^2+(\alpha+\beta+1)^2\right)f(\rho)=0.
\end{equation}
When $\lambda\neq -i, -2i,\cdot\cdot\cdot$, this Jacobi equation has two different fundamental solutions $\Psi^{\lambda}_{\alpha,\beta}(\rho)$ and $\Psi^{-\lambda}_{\alpha,\beta}(\rho),$
where $$\Psi^{\lambda}_{\alpha,\beta}(\rho)=(2\sinh\rho)^{i\lambda-\alpha-\beta-1}F(\frac{\alpha+\beta+1-i\lambda}{2}, \frac{-\alpha+\beta+1-i\lambda}{2}, 1-i\lambda, -\sinh^{-2}\rho)$$
and $F$ denotes the hypergeometric function. The hypergeometric function $F(a;b;c;z)$
is given by
$$F(a;b;c;z)=\sum_{k=0}^{\infty}\frac{(a)_k(b)_k}{(c)_k}\frac{z^k}{k!}$$
with $a$, $b\in \mathbb{C}$, $c\neq 0, -1, 2,\cdot\cdot\cdot,$ where
$$(a)_0=1,\ \ (a)_k=a(a+1)\cdot\cdot\cdot(a+k-1)\ {\rm for}\ k\geq 1.$$
This series gives an analytic function for $|z|<1$ and it can be continued to the $\mathbb{C}\setminus [1,+\infty)$. Furthermore $F(a,b,c,0)=1$. We refer to Chapter II of \cite{Erd} for more details about hypergeometric function. Coming back to the equation which $a_{m,l}(\rho)$ satisfies
$$\left(\frac{\partial^2}{\partial \rho^2}+(n-1)\coth \rho \frac{\partial}{\partial \rho}+\frac{(n-1)^2}{4}+\mu^2\right)a_{m,l}(\rho)=0,$$
this equation has two different solution $\Psi^{\mu}_{\frac{n-2}{2},-\frac{1}{2}}(\rho)$ and $\Psi^{-\mu}_{\frac{n-2}{2},-\frac{1}{2}}(\rho)$. Then according to second-order linear ODE theory,  we deduce that there exist $C_1$ and $C_2$ such that $$a_{m,l}(\rho)=C_1\Psi^{\mu}_{\frac{n-2}{2},-\frac{1}{2}}(\rho)+C_2\Psi^{-\mu}_{\frac{n-2}{2},-\frac{1}{2}}(\rho).$$
If $C_1^2+C_2^2\neq 0$, in consideration of $$\Psi^{\mu}_{\frac{n-2}{2},-\frac{1}{2}}(\rho)=(2\sinh \rho)^{-\frac{n-1}{2}+i\mu}F(\frac{n-1-2i\mu}{4},\frac{-(n-3)-2i\mu}{4}, 1-i\mu,-\sinh^{-2}\rho)$$
and $F(\frac{n-1-2i\mu}{4},\frac{-(n-3)-2i\mu}{4}, 1-i\mu,0)=1$ (see Chapter II of \cite{Erd}), we deduce that
\begin{equation}\label{asymptotic}
a_{l,m}(\rho)=O\left(\sinh^{-\frac{n-1}{2}}\rho\right).
\end{equation}
On the other hand, direct calculation gives that
\begin{equation}\begin{split}
\int_{\partial B_{\mathbb{B}^n}(0,R)}|u|^2d\sigma_{\mathbb{B}_n}&=\int_{\partial B(0,\tanh(\frac{R}{2}))}|u|^2\left(\frac{2}{1-|x|^2}\right)^{n-1}d\sigma\\
&=\int_{\partial B(0,\tanh(\frac{R}{2}))}|u|^2\left(2\cosh^2 \frac{\rho}{2}\right)^{n-1}d\sigma\\\
&=\sum_{l=1}^{+\infty}\sum_{m=1}^{N(l,n)}a_{m,l}^2(R)\left(2\cosh^2 \frac{R}{2}\right)^{n-1}\tanh^{n-1}(\frac{R}{2})\int_{\mathbb{S}^{n-1}}|Y_{m,l}(x')|^2d\omega\\
&=\sum_{l=1}^{+\infty}\sum_{m=1}^{N(l,n)}a_{m,l}^2(R)\left(2\cosh^2 \frac{R}{2}\right)^{n-1}\tanh^{n-1}(\frac{R}{2}),\\
\end{split}\end{equation}
where $d\sigma_{\mathbb{B}_n}$ denotes the hyperbolic surface measure in $\partial B_{\mathbb{B}^n}(0,R)$ and $d\sigma$ denotes the
surface measure in $\partial B(0,\tanh(\frac{R}{2}))$. This together with $\lim\limits_{R\rightarrow +\infty}\int_{\partial B_{\mathbb{B}^n}(0,R)}|u|^2d\sigma_{g}=0$ yields that
$$\lim\limits_{R\rightarrow +\infty}\sum_{l=1}^{+\infty}\sum_{m=1}^{N(l,n)}a_{m,l}^2(R)\left(2\cosh^2 \frac{R}{2}\right)^{n-1}\tanh^{n-1}(\frac{R}{2})=0.$$
This implies that $a_{m,l}(R)=o(\sinh^{-\frac{n-1}{2}} R)$ when $R\rightarrow +\infty$. Then we arrive at a contradiction with
\eqref{asymptotic}. Hence, we conclude that $u$ is equal to zero in $\mathbb{B}^n\setminus B_{\mathbb{B}^n}(0,R_0)$. Finally, we show that
$u$ must be equal to zero in $\mathbb{B}^n\setminus \Omega$ as a consequence of strong unique continuation property for elliptic PDEs of $\mathbb{R}^n$. In fact, choose $R_1>R_0$, Helmholtz equation
\begin{equation*}
-\Delta_{\mathbb{B}^n}u-\frac{(n-1)^2}{4}u-\mu^2 u=0,\ \ x\in B_{\mathbb{B}^n}(0,R_1)\setminus \Omega
\end{equation*}
can be written as the following elliptic equation in $B^n(0, \tanh(R_1/2))\setminus \Omega$
\begin{equation*}
\left(\frac{2}{1-|x|^2}\right)^{-\frac{n}{2}-1} \Delta_{\mathbb{R}^n}\left(\left(\frac{2}{1-|x|^2}\right)^{\frac{n}{2}-1}u\right)+\left(\frac{1}{4}+\mu^2\right)u=0,\ \ x\in B^n(0, \tanh(R_1/2)\setminus \Omega.
\end{equation*}
Let $\tilde{u}(x)=\left(\frac{2}{1-|x|^2}\right)^{\frac{n}{2}-1}u$, the above equation can be simplified as
\begin{equation}\label{elliptic}
\Delta_{\mathbb{R}^n}\tilde{u}+V(x)\tilde{u}=0,\ \ x\in B^n(0, \tanh(\frac{R_1}{2}))\setminus \Omega,
\end{equation}
where $V(x)=\left(\frac{1}{4}+\mu^2\right)\left(\frac{2}{1-|x|^2}\right)^{2}$ is bounded in $L^{\infty}(B^n(0, \tanh(\frac{R_1}{2}))\setminus \Omega)$. Recall we have proved that $u$ is equal to zero in $\mathbb{B}^n\setminus B_{\mathbb{B}^n}(0,R_0)$, hence $$\tilde{u}=0,\ \ x\in B^n(0, \tanh(R_1/2))\setminus B^n(0, \tanh(R_0/2)).$$  Applying the strong unique continuation property into equation \eqref{elliptic}, we derive that $\widetilde{u}$ must be equal to zero in $B^n(0, \tanh(\frac{R_1}{2}))\setminus \Omega$, that is $u=0$ in $B_{\mathbb{B}^n}(0,R_1)\setminus \Omega$. In summary, we prove that $u$ must be equal to zero in $\mathbb{B}^n\setminus \Omega$.
Then the proof of Theorem \ref{thm3} is accomplished.

\end{proof}

\subsection{Scattering theory for Helmholtz equation involving potential }
In this subsection, we mainly describe the scattering theory for the Helmholtz equation from the potential on hyperbolic space $\mathbb{B}^n$. In Euclidean space, an incident wave $u^{i}$ is defined through $u^i=e^{i\mu x\cdot \xi}$. Obviously $u^{i}$ is a periodic function and oscillate rapidly when the frequency $\mu$ is bigger. On the other hand, $e^{i\mu x\cdot \xi}$ is bounded and plays an important role in the fields of harmonic analysis. Indeed, $e^{i\mu x\cdot \xi}$ is spectrum of Laplacian operator in $\mathbb{R}^n$. Any function $f(x)\in C^{\infty}_c(\mathbb{R}^n)$ can be written as
$$f(x)=\frac{1}{2\pi}\int_{\mathbb{R}^n}\hat{f}(\xi)e^{ix\cdot \xi}d\xi,$$
where $$\hat{f}(\xi)=\frac{1}{2\pi}\int_{\mathbb{R}^n}f(x)e^{ix\cdot \xi}dx.$$  If we consider the scattering theory in hyperbolic space $\mathbb{B}^n$, we first need to define an suitable incident wave in hyperbolic space. Recall the Fourier analysis in hyperbolic space.
The Fourier transform on hyperbolic space of a function $f\in C^{\infty}_{c}(\mathbb{B}^n)$ is defined as
$$\widehat{f}(\mu,\xi)=\int_{\mathbb{B}^n}f(x)e_{-\mu, \xi}(x)dV_{\mathbb{B}^n},$$
where $$e_{-\mu,\xi}(x)=\left(\frac{\sqrt{1-|x|^2}}{|x-\xi|}\right)^{n-1-\mathrm{i}\mu},\ x\in \mathbb{B}^n,\ \mu\in \mathbb{R},\ \  \xi \in \mathbb{S}^{n-1}.$$
Furthermore, we have the inversion formula:
$$f(x)=D_{n}\int_{-\infty}^{\infty}\int_{\mathbb{S}^{n-1}}\hat{f}(\mu,\xi)e_{\lambda,\xi}(x)|c(\mu)|^{-2}d\mu d\xi,$$
where $D_n=\frac{1}{2^{3-n}\pi |\mathbb{S}^{n-1}|}$ and $c(\mu)$ is the Harish-Chandra function. Hence it is reasonable to define \textbf{an incident wave in hyperbolic space $\mathbb{B}^n$ through $u^{i}=e_{-2\mu,\xi}(x)$}. Recall that $e_{-2\mu, \xi}(x)$ is the spectrum of the operator $-\Delta_{\mathbb{B}^n}-\frac{(n-1)^2}{4}$ in hyperbolic space $\mathbb{B}^n$, that is $e_{-2\mu, \xi}(x)$ satisfying equation $$-\Delta_{\mathbb{B}^n}(e_{-2\mu,\xi})-\frac{(n-1)^2}{4}e_{-2\mu,\xi}=\mu^2 e_{-2\mu,\xi}.$$  Obviously,
for any $\hat{x}\neq \xi$, $$\lim\limits_{\rho(x)\rightarrow +\infty}e_{-2\mu,\xi}(x)=\lim\limits_{\rho(x)\rightarrow +\infty}\left(\frac{\sqrt{1-|x|^2}}{|x-\xi|}\right)^{n-1-2\mathrm{i}\mu}=0$$
and for $\hat{x}=\xi$, $$\lim\limits_{\rho(x)\rightarrow +\infty}|e_{-2\mu,\xi}|=\lim\limits_{\rho(x)\rightarrow +\infty}|\left(\frac{\sqrt{1-|x|^2}}{|x-\xi|}\right)^{n-1-2\mathrm{i}\mu}|=+\infty,$$
which is obviously different from an incident wave $e^{-i\mu x\cdot \xi}$ in Euclidean space since absolute value of $e^{-i\mu x\cdot \xi}$ is still equal to one.
\medskip

 On the other hand, careful computation gives that
\begin{equation}\begin{split}
\frac{\partial e_{-2\mu,\xi}(x)}{\partial \rho}&=\left(\mu \mathrm{i}-\frac{n-1}{2}\right)\left(\cosh(\frac{\rho}{2})\right)^{2\mu\mathrm{i}-(n-1)}\tanh(\frac{\rho}{2})|\tanh(\frac{\rho}{2})\hat{x}-\xi|^{2\mu\mathrm{i}-(n-1)}\\
&\ \ +\left(\cosh(\frac{\rho}{2})\right)^{2\mu\mathrm{i}-(n+1)}\left(\mu\mathrm{i}-\frac{n-1}{2}\right) |\tanh(\frac{\rho}{2})\hat{x}-\xi|^{2\mu\mathrm{i}-(n+1)}\left(\tanh(\frac{\rho}{2})\hat{x}-\xi\right)\cdot\hat{x}.
\end{split}\end{equation}
Then it follows that
\begin{equation}\begin{split}
&\frac{\partial  e_{-2\mu,\xi}}{\partial \rho}-\left(\mu\mathrm{i}-\frac{n-1}{2}\right)\tanh(\frac{\rho}{2})  e_{-2\mu,\xi}\\
&\ \ =\left(\cosh(\frac{\rho}{2})\right)^{2\mu\mathrm{i}-(n+1)}\left(\mu\mathrm{i}-\frac{n-1}{2}\right) |\tanh(\frac{\rho}{2})\hat{x}-\xi|^{2\mu\mathrm{i}-(n+1)}\left(\tanh(\frac{\rho}{2})\hat{x}-\xi\right)\cdot\hat{x}.
\end{split}\end{equation}
This gives
\begin{equation}\label{p-wave}
\frac{\partial  e_{-2\mu,\xi}}{\partial \rho}-\left(\mu\mathrm{i}-\frac{n-1}{2}\right)\tanh(\frac{\rho}{2})  e_{-2\mu,\xi}=O(\sinh^{-\frac{n+1}{2}} \rho)
\end{equation}
for any $\hat{x}\neq \xi$. It deduces that $e_{-2\mu,\xi}$ satisfies hyperbolic radiation condition for any direction $\hat{x}\neq \xi$. This is the second main difference between incident wave $e_{-2\mu,\xi}$ of hyperbolic space and the counterpart $e^{-i\mu x\cdot \xi}$ of Euclidean space since $e^{-i\mu x\cdot \xi}$ does not satisfies radiation condition of Euclidean space for any direction $\hat{x}$.
\medskip

Furthermore, in Euclidean space, $e^{-i\mu x\cdot \xi}$ and $e^{i\mu x\cdot \xi}$ is essential same. However, in the setting of hyperbolic space, $e_{-2\mu,\xi}$ and $e_{2\mu,\xi}$ has the essential difference. Indeed, we can check that $e_{-2\mu,\xi}$ satisfies the hyperbolic radiation condition
$$\frac{\partial  e_{-2\mu,\xi}}{\partial \rho}-\left(\mu\mathrm{i}-\frac{n-1}{2}\right)\tanh(\frac{\rho}{2})  e_{-2\mu,\xi}=O(\sinh^{-\frac{n+1}{2}} \rho)$$
for any $\hat{x}\neq \xi$. Careful computations gives that $e_{2\mu,\xi}$ satisfies the condition
$$\frac{\partial  e_{2\mu,\xi}}{\partial \rho}-\left(-\mu\mathrm{i}-\frac{n-1}{2}\right)\tanh(\frac{\rho}{2})  e_{2\mu,\xi}=O(\sinh^{-\frac{n+1}{2}} \rho)$$
for any direction $\hat{x}\neq \xi$. This is the third main difference between the incident wave in hyperbolic space and Euclidean space.
\medskip

Now, we start to consider the scattering theory from potential. A $L^{\infty}$ potential $q(x)$ satisfying $V(x):=q(x)-1$ is equal to zero in the supplement of $B_{\mathbb{B}^n}(0,R_0)$ will create the new wave
$u$ solving the equation
$$\left(-\Delta_{\mathbb{B}^n}-\frac{(n-1)^2}{4}-\mu^2q\right)u=0,\ \ x\in \mathbb{B}^n.$$
The wave $u$ and $u^{i}$ is linked, the essential connection being that the scattered wave $u^{s}=u-u^{i}$ satisfies the hyperbolic radiation condition
$$\frac{\partial u^s}{\partial \rho}-\left(\mu i-\frac{n-1}{2}\right)\tanh(\frac{\rho}{2})u^s=o\left(\frac{1}{\sinh^{\frac{n-1}{2}}(\rho)}\right)$$
as $\rho\rightarrow +\infty$. Since the incident wave $u^{i}$ satisfies equation
$$\left(-\Delta_{\mathbb{B}^n}-\frac{(n-1)^2}{4}-\mu^2\right)u^i=0.$$
Hence it is not difficult to check that the scattered wave $u^{s}$ satisfies equation
\begin{equation}\label{scattering}
\left(-\Delta_{\mathbb{B}^n}-\frac{(n-1)^2}{4}-\mu^2\right)u^{s}=\mu^2V(x) u^{s}+\mu^2V(x) u^{i},\ \ x\in \mathbb{B}^n.
\end{equation}
We will prove that the total wave $u$ is existing and unique.
\begin{thm}\label{thm4}
Assume that $u$ satisfies equation
\begin{equation}\label{scattering}
\left(-\Delta_{\mathbb{B}^n}-\frac{(n-1)^2}{4}-\mu^2q\right)u=0,\ \ x\in \mathbb{B}^n,
\end{equation}
and $u-u^i$ satisfies the hyperbolic radiation condition
$$\frac{\partial (u-u^i)}{\partial \rho}-(\mu \mathrm{i}-\frac{n-1}{2})\tanh(\frac{\rho}{2})(u-u^i)=o\left(\frac{1}{\sinh^{\frac{n-1}{2}}(\rho)}\right).$$
If imaginary part $\mathbb{\Im}(q)$of the potential $q$ is non-negative, then there exits a unique solution $u\in H^1_{loc}(\mathbb{B}^n)$.
\end{thm}

\begin{proof}
Since $u(x)=e_{-2\mu, \xi}(x)+u^s(x)$, we only need to prove that $u^s$ is existing and unique. Obviously $u^s$ satisfies
equation
\begin{equation}\begin{cases}\label{scattering}
\left(-\Delta_{\mathbb{B}^n}-\frac{(n-1)^2}{4}-\mu^2\right)u^s=f(x),\ \ x\in \mathbb{B}^n,\\
\frac{\partial u^s}{\partial \rho}-(\mu \mathrm{i}-\frac{n-1}{2})\tanh(\frac{\rho}{2})u^s=o\left(\frac{1}{\sinh^{\frac{n-1}{2}}(\rho)}\right)\\
f(x)=\mu^2 V(x)u^{s}+\mu^2 V(x) u^{i},\ \ x\in \mathbb{B}^n.
\end{cases}\end{equation}

Since $u^s$ satisfies the hyperbolic radiation condition, applying Theorem \ref{thm2}, the existence and uniqueness of equation \eqref{scattering} is equivalent to the existence and uniqueness of the following integral equation
\begin{equation}\label{int1}
u^{s}(x)=\int_{\mathbb{B}^n}G_{-\mu \mathrm{i}}(\rho(x,y))f(y)dV_{\mathbb{B}^n}(y),\ \ x\in \mathbb{B}^n.
\end{equation}
Define the integral operator $K$ by $$K(u^{s})=\int_{\mathbb{B}^n}G_{-\mu \mathrm{i}}(\rho(x,y))\mu^2V(y) u^{s}(y)dV_{\mathbb{B}^n}(y).$$
The integral equation \eqref{int1} can be rewritten as
$$(I-K)u^{s}(x)=\int_{\mathbb{B}^n}G_{-\mu \mathrm{i}}(\rho(x,y))V(y)\mu^2 u^{i}(y)dV_{\mathbb{B}^n}(y),\ \ x\in \mathbb{B}^n.$$
Now we claim $\int_{\mathbb{B}^n}G_{-\mu \mathrm{i}}(\rho(x,y))V(y)\mu^2 u^{i}(y)dV_{\mathbb{B}^n}(y)$ is bounded in $L^2(B^n)$.
Recall the Hardy-Littlewood-Sobolev inequality in $\mathbb{R}^n$ which states that
$$\int_{\mathbb{R}^n}\int_{\mathbb{R}^n}\frac{h(x)g(y)}{|x-y|^{\lambda}}dxdy\lesssim \left(\int_{\mathbb{R}^n}|h|^pdx\right)^{\frac{1}{p}}\left(\int_{\mathbb{R}^n}|g|^qdx\right)^{\frac{1}{q}}$$
for any $0<\lambda<n$, $1<p,q<+\infty$ with $\frac{1}{p}+\frac{1}{q}+\frac{\lambda}{n}=2$.
Since
\begin{equation}\begin{split}
G_{-\mu i}(\rho(x,y))&=A_{n,-\mu \mathrm{i}} \frac{(\cosh \rho)^{\frac{n-3}{2}+\mu \mathrm{i}}}{(\sinh \rho)^{n-2}}\int_{0}^{\pi}\left(1+\frac{\cos t}{\cosh \rho}\right )^{\frac{n-3}{2}+\mu \mathrm{i}}(\sin t)^{-2\mu \mathrm{i}}dt\\
&\lesssim \left(\frac{1}{\sinh \frac{\rho(x,y)}{2}}\right)^{n-2},
\end{split}\end{equation}
then it follows that
\begin{equation}\begin{split}
&\int_{\mathbb{B}^n}G_{-\mu \mathrm{i}}(\rho(x,y))V(y)\mu^2 u^{i}(y)dV_{\mathbb{B}^n}(y)\\
&\ \ =\int_{\mathbb{B}^n}\left(|x-y|\cosh(\frac{\rho(x)}{2})\cosh(\frac{\rho(y)}{2})\right)^{-(n-2)}V(y)\mu^2 u^{i}(y)dV_{\mathbb{B}^n}(y)\\
&\ \ \lesssim \int_{B^n}V(y)\frac{\mu^2 u^{i}(y)}{|x-y|^{n-2}}dy
\end{split}\end{equation}
where we use $\sinh(\frac{\rho(x,y)}{2})=|x-y|\cosh(\frac{\rho(x)}{2})\cosh(\frac{\rho(y)}{2})$ and $\mathcal{V}(x)=1$ on $\mathbb{B}^n\setminus B_{\mathbb{B}^n}(0,R_0)$.
Applying the Hardy-Littlewood-Sobolev inequality, one can get
\begin{equation}\begin{split}
&\int_{B^n}|\int_{\mathbb{B}^n}G_{-\mu i}(\rho(x,y))\left(\mathcal{V}(y)-1\right)\mu^2 u^{i}(y)dV_{\mathbb{B}^n}(y)|^2dx\\
&\ \ \lesssim \int_{B^n}|\int_{B^n}V(y)\frac{\mu^2 u^{i}(y)}{|x-y|^{n-2}}dy|^2dx\\
&\ \ \lesssim \left(\int_{B^n}|u^i|^{\frac{2n}{n+4}}dx\right)^{\frac{n+4}{n}}
\end{split}\end{equation}
if $n\geq 5$ and
\begin{equation}\begin{split}
&\int_{B^n}|\int_{\mathbb{B}^n}G_{-\mu \mathrm{i}}(\rho(x,y))\left(\mathcal{V}(y)-1\right)\mu^2 u^{i}(y)dV_{\mathbb{B}^n}(y)|^2dx\\
&\ \ \lesssim \left(\int_{B^n}|\int_{\mathbb{B}^n}G_{-\mu \mathrm{i}}(\rho(x,y))\left(\mathcal{V}(y)-1\right)\mu^2 u^{i}(y)dV_{\mathbb{B}^n}(y)|^4dx\right)^{\frac{1}{2}}\\
&\ \ \lesssim \left(\int_{\Omega}|u^i|^{\frac{4n}{n+8}}dx\right)^{\frac{n+8}{2n}}
\end{split}\end{equation}
if $n=3, 4$, where $B^n$ denotes the unit ball in $\mathbb{R}^n$.
Now, we prove that the integral equation
\begin{equation}\label{integral}\begin{split}
(I-K)u^{s}(x)=\int_{\mathbb{B}^n}G_{-\mu \mathrm{i}}(\rho(x,y))V(y)\mu^2 u^{i}(y)dV_{\mathbb{B}^n}(y)\in L^{2}(B^n),\ \ x\in \mathbb{B}^n
\end{split}\end{equation} has a unique solution $u^{s}\in L^2(B^n)$. We first claim that $K$ is a compact operator from $L^2(B^n)$ to $L^2(B^n)$.
Recalling the definition of the operator $K$, since $$G_{-\mu \mathrm{i}}(\rho(x,y))\simeq |x-y|^{-(n-2)}\  {\rm when}\  \rho(x,y)\rightarrow 0$$ and $V(x)$ has the compact support, hence one can obtain for any $h\in L^2(B^n)$, there holds
$$\int_{B^n}\left(|\nabla K(h)|^2+|K(h)|^2\right)dx\lesssim \int_{B^n}|h|^2dx,$$
which together with the Sobolev imbedding in $W^{1,2}(B^n)$ implies that
$K$ is a compact operator from $L^2(B^n)$ to $L^2(B^n)$. Hence according to Fredholm theorem for compact operator (see \cite{Evans}), we know that the integral equation \eqref{integral} has a unique solution $u^s\in L^2(B^n)$ if and only if the linear equation $$(I-K)u^{s}(x)=0,\ \ \ \ x\in \mathbb{B}^n$$ has just the zero solution. This is equivalent to prove that the differential equation
$$\left(-\Delta-\frac{(n-1)^2}{4}-\mu^2q(x)\right)u^s=0,\ \ x\in \mathbb{B}^n$$ has only zero solution through Green representative formula in Theorem \ref{thm2}.
Noticing $q(x)=1$ in $\mathbb{B}^n\setminus B_{\mathbb{B}^n}(0,R_0)$, if we can prove that $u^s$ satisfies the condition of Rellich theorem in hyperbolic space (Theorem \ref{thm3}), we conclude that $u^s(x)$ is equal to zero in $\mathbb{B}^n\setminus B_{\mathbb{B}^n}(0,R_0)$. This together with the strong unique continuation theorem yields that $u^s$ must be equal to zero on the whole hyperbolic space $\mathbb{B}^n$. We now turn to proving $$\lim\limits_{R\rightarrow +\infty}\int_{\partial B_{\mathbb{B}^n}(0,R)}|u^s|^2d\sigma_{g}=0.$$
Recalling the proof of Theorem \ref{thm2}, we have obtained that
\begin{equation}\label{eq6}\begin{split}
o(1)\geq \int_{\partial B_{\mathbb{B}^n}(0,R)}\mu^2\tanh^2(\frac{\rho}{2})|u^s|^2d\sigma_{\mathbb{B}^n}+2\Im\left(\tanh(\frac{R}{2})\int_{\partial B_{\mathbb{B}^n}(0,R)}\frac{\partial \overline{u^s}}{\partial\rho}u^sd\sigma_{\mathbb{B}^n}\right)
\end{split}\end{equation}
when $R\rightarrow +\infty$.
Applying the Green formula in Lemma \ref{lem:green-ball}, one has
\begin{equation}\begin{split}
\int_{\partial B_{\mathbb{B}^n}(0,R)}\frac{\partial \overline{u^s}}{\partial\rho}u^sd\sigma_{\mathbb{B}^n}&=\int_{B_{\mathbb{B}^n}(0,R)\setminus B_{\mathbb{B}^n}(0,R_0)}\Delta_{\mathbb{B}^n}(\overline{u^s})u^sdV_{\mathbb{B}^n}+\int_{B_{\mathbb{B}^n}(0,R)\setminus B_{\mathbb{B}^n}(0,R_0)}|\nabla_{\mathbb{B}^n}u^s|^2dV_{\mathbb{B}^n}\\
&\ \ +\int_{\partial B_{\mathbb{B}^n}(0,R_0)}\frac{\partial \overline{u^s}}{\partial\rho}u^sd\sigma_{\mathbb{B}^n}\\
&=-\int_{B_{\mathbb{B}^n}(0,R)\setminus B_{\mathbb{B}^n}(0,R_0)}\left(\frac{(n-1)^2}{4}+\mu^2\right)|u^s|^2dV_{\mathbb{B}^n}+\int_{B_{\mathbb{B}^n}(0,R)\setminus B_{\mathbb{B}^n}(0,R_0)}|\nabla_{\mathbb{B}^n}u^s|^2dV_{\mathbb{B}^n}\\
&\ \ +\int_{\partial B_{\mathbb{B}^n}(0,R_0)}\frac{\partial \overline{u^s}}{\partial\rho}u^sd\sigma_{\mathbb{B}^n}
\end{split}\end{equation}
and
\begin{equation}\begin{split}
\int_{\partial B_{\mathbb{B}^n}(0,R_0)}\frac{\partial \overline{u^s}}{\partial\rho}u^sd\sigma_{\mathbb{B}^n}&=\int_{ B_{\mathbb{B}^n}(0,R_0)}\Delta_{\mathbb{B}^n}(\overline{u^s})u^sdV_{\mathbb{B}^n}+\int_{ B_{\mathbb{B}^n}(0,R_0)}|\nabla_{\mathbb{B}^n}u^s|^2dV_{\mathbb{B}^n}\\
&=-\int_{B_{\mathbb{B}^n}(0,R_0)}\left(\frac{(n-1)^2}{4}+\mu^2q(x)\right)|u^s|^2dV_{\mathbb{B}^n}+\int_{B_{\mathbb{B}^n}(0,R_0)}|\nabla_{\mathbb{B}^n}u^s|^2dV_{\mathbb{B}^n}.
\end{split}\end{equation}
This implies that
$$\Im\left(\tanh(\frac{R}{2})\int_{\partial B_{\mathbb{B}^n}(0,R)}\frac{\partial \overline{u^s}}{\partial\rho}u^sd\sigma_{\mathbb{B}^n}\right)\geq 0$$ if $\Im\left(q\right)\geq 0$. Combining this and \eqref{eq6}, we conclude that $$\lim\limits_{R\rightarrow +\infty}\int_{\partial B_{\mathbb{B}^n}(0,R)}|u^s|^2d\sigma_{\mathbb{B}^n}=0.$$ Then the proof of Theorem \ref{thm4} is accomplished.
\end{proof}
\begin{rmk}\label{acouse}
We have proved that $u$ exists and is unique. Furthermore, the total wave $u$ can be written as far-field pattern form.
In fact, for any $x\in \mathbb{B}^n$, applying Green's representation formula and the accurate expansion of $G_{-\mu i}$, we have
\begin{equation}\begin{split}
u(x)&=u^{i}(x)+\int_{\mathbb{B}^n}G_{-\mu \mathrm{i}}(\rho(x,y))f(y)dV_{\mathbb{B}^n}(y)\\
&=u^{i}(x)+\frac{\left(\cosh(\frac{\rho(x)}{2})\right)^{2\mu \mathrm{i}}}{\cosh^{n-1}(\frac{\rho(x)}{2})}\int_{0}^{\pi}(\sin t)^{-2\mu \mathrm{i}}dt\\
&\ \ \times \frac{A_{n,-\mu \mathrm{i}}}{2^{\frac{n-1}{2}}}\int_{\mathbb{B}^n}\frac{\left(\cosh(\frac{\rho(y)}{2})\right)^{2\mu \mathrm{i}}}{\cosh^{n-1}(\frac{\rho(y)}{2})}\left(1-2\hat{x}\cdot y+|y|^2\right)^{\mu \mathrm{i}-\frac{n-1}{2}}f(y)dV_{\mathbb{B}^n}(y)+O\left(\frac{1}{\sinh^{\frac{n+1}{2}}(\rho(x))}\right)\\
&=u^{i}(x)+\left(\frac{\left(\cosh(\frac{\rho(x)}{2})\right)^{2\mu \mathrm{i}}}{\cosh^{n-1}(\frac{\rho(x)}{2})} \right) u^{\mu}_{\infty}(\hat{x}, \xi)+O\left(\frac{1}{\sinh^{\frac{n+1}{2}}(\rho(x))}\right),
\end{split}\end{equation}
where $f(x)=\mu^2 V(x)u(x)$ has the compact support in $B_{\mathbb{B}^n}(0,R_0)$ and
$$u^{\mu}_{\infty}(\hat{x}, \xi)=\frac{A_{n,-\mu \mathrm{i}}}{2^{\frac{n-1}{2}}}\int_{0}^{\pi}(\sin t)^{-2\mu \mathrm{i}}dt\int_{\mathbb{B}^n}\frac{\left(\cosh(\frac{\rho(y)}{2})\right)^{2\mu \mathrm{i}}}{\cosh^{n-1}(\frac{\rho(y)}{2})}\left(1-2\hat{x}\cdot y+|y|^2\right)^{\mu \mathrm{i}-\frac{n-1}{2}}f(y)dV_{\mathbb{B}^n}(y).$$
\end{rmk}

\subsection{Scattering theory for obstacle problems in hyperbolic space}

Consider scattering by impenetrable obstacles in hyperbolic space. Let $\Omega\in \mathbb{B}^n$ be a bounded domain, representing an obstacle that prohibits wave penetration. The physical properties of the obstacle are characterized by a boundary condition imposed on the total field $u$ along $\partial\Omega$. Common examples include:
\begin{itemize}
    \item \emph{Soft (Dirichlet) obstacle:} $u|_{\partial\Omega} = 0$,
    \item \emph{Hard (Neumann) obstacle:} $\frac{\partial u}{\partial\nu_{\mathbb{B}^n}}|_{\partial\Omega} = 0$,
    \item \emph{Impedance obstacle:} $(\frac{\partial u}{\partial\nu_{\mathbb{B}^n}} + \eta u)|_{\partial\Omega} = 0$,
\end{itemize}
where $\eta \in L^\infty(\partial\Omega)$ with $\Im\eta \geq 0$. The forward scattering problem for an obstacle $\Omega$ is then formulated as:
\begin{equation}\label{obs}\begin{cases}
&-\Delta_{\mathbb{B}^n}u-\frac{(n-1)^2}{4}u=\mu_0^2u,\ \ x\in \mathbb{B}^n\setminus \Omega\\
&\mathcal{B}(u)=0,\ \ x\in \partial \Omega,
\end{cases}\end{equation}
where $\mathcal{B}$ denotes the boundary operator, $u(x)=u^{i}+u^s(x)$, $u^i(x)=e_{-2\mu_0, \xi}(x)$ is an incident wave satisfying
$$-\Delta_{\mathbb{B}^n}u^{i}-\frac{(n-1)^2}{4}u^i=\mu_0^2u^i$$ and $u^s$ is the scattered wave satisfying the hyperbolic radiation condition
$$\frac{\partial u^s}{\partial \rho}-(\mu_0 \mathrm{i}-\frac{n-1}{2})\tanh(\frac{\rho}{2})u^s=o\left(\frac{1}{\sinh^{\frac{n-1}{2}}(\rho)}\right),$$

It is easy to check that $u^{s}$ satisfies equation
\begin{equation}\label{scas}\begin{cases}
&-\Delta_{\mathbb{B}^n}u^s-\frac{(n-1)^2}{4}u^s=\mu_0^2u^s,\ \ x\in \mathbb{B}^n\setminus \Omega.\\
&\mathcal{B}(u^s+u^{i})=0,\ \ x\in \partial \Omega\\
&\frac{\partial u_s}{\partial \rho}-(\mu_0 \mathrm{i}-\frac{n-1}{2})\tanh(\frac{\rho}{2})u_s=o\left(\frac{1}{\sinh^{\frac{n-1}{2}}(\rho)}\right).
\end{cases}\end{equation}
We claim that
\begin{thm}\label{obstacle1}
There exists a unique solution $u\in H_{loc}^1(\mathbb{B}^n)$ to the scattering system \eqref{scas}.
\end{thm}
\begin{proof}
Without loss of generality, we only consider the soft boundary case, other case is similar. We first show the uniqueness. Assume that $u_1$ and $u_2$ satisfies equation
\begin{equation}\begin{cases}
&-\Delta_{\mathbb{B}^n}u-\frac{(n-1)^2}{4}u=\mu_0^2u,\ \ x\in \mathbb{B}^n\setminus \Omega\\
&u=0,\ \ x\in \partial \Omega,\\
&\frac{\partial (u-u^i)}{\partial \rho}-(\mu_0 \mathrm{i}-\frac{n-1}{2})\tanh(\frac{\rho}{2})(u-u^i)=o\left(\frac{1}{\sinh^{\frac{n-1}{2}}(\rho)}\right).
\end{cases}\end{equation}
Then $u_1-u_2$ satisfies equation
\begin{equation}\label{adeq1}\begin{cases}
&-\Delta_{\mathbb{B}^n}(u_1-u_2)-\frac{(n-1)^2}{4}(u_1-u_2)=\mu_0^2(u_1-u_2),\ \ x\in \mathbb{B}^n\setminus \Omega.\\
&u_1-u_2=0,\ \ x\in \partial \Omega,\\
&\frac{\partial (u_1-u_2)}{\partial \rho}-(\mu_0 \mathrm{i}-\frac{n-1}{2})\tanh(\frac{\rho}{2})(u_1-u_2)=o\left(\frac{1}{\sinh^{\frac{n-1}{2}}(\rho)}\right).
\end{cases}\end{equation}

Careful calculation gives that
\begin{equation}\label{addexpansion1}\begin{split}
&\int_{\partial B_{\mathbb{B}^n}(0,R)}\left|\frac{\partial (u_1-u_2)}{\partial \rho}-\left(\mu i-\frac{n-1}{2}\right)\tanh(\frac{\rho}{2})(u_1-u_2)\right|^2d\sigma_{\mathbb{B}^n}\\
&\ \ =\int_{\partial B_{\mathbb{B}^n}(0,R)}\left|\frac{\partial (u_1-u_2)}{\partial\rho}\right|^2d\sigma_{\mathbb{B}^n}+\int_{\partial B_{\mathbb{B}^n}(0,R)}\left(\frac{(n-1)^2}{4}+\mu^2\right)\tanh^2(\frac{\rho}{2})|u_1-u_2|^2d\sigma_{\mathbb{B}^n}\\
&\ \ \ \ +2\Im\left(\int_{\partial B_{\mathbb{B}^n}(0,R)}\tanh(\frac{\rho}{2})\frac{\partial \left(\overline{u_1-u_2}\right)}{\partial\rho}(u_1-u_2)d\sigma_{\mathbb{B}^n}\right)\\
&\ \ \ \ +(n-1)\Re\left(\int_{\partial B_{\mathbb{B}^n}(0,R)}\tanh(\frac{\rho}{2})\frac{\partial \left(\overline{u_1-u_2}\right)}{\partial\rho}(u_1-u_2)d\sigma_{\mathbb{B}^n}\right).
\end{split}\end{equation}
Since
\begin{equation}\begin{split}
&\int_{\partial B_{\mathbb{B}^n}(0,R)}\left(\left|\frac{\partial \left(u_1-u_2\right)}{\partial\rho}\right|^2+\frac{(n-1)^2}{4}\tanh^2(\frac{\rho}{2})|u_1-u_2|^2\right)d\sigma_{\mathbb{B}^n}\\
&\ \ \geq (n-1)\Re\left(\int_{\partial B_{\mathbb{B}^n}(0,R)}\tanh(\frac{\rho}{2})\frac{\partial \left(\overline{u_1-u_2}\right)}{\partial\rho}\left(u_1-u_2\right)d\sigma_{\mathbb{B}^n}\right),
\end{split}\end{equation}
combining this and \eqref{addexpansion1} and hyperbolic radiation condition, we derive that
\begin{equation}\label{addeq1}\begin{split}
o(1)\geq \int_{\partial B_{\mathbb{B}^n}(0,R)}\mu^2\tanh^2(\frac{\rho}{2})|u_1-u_2|^2d\sigma_{\mathbb{B}^n}+2\Im\left(\int_{\partial B_{\mathbb{B}^n}(0,R)}\tanh(\frac{\rho}{2})\frac{\partial \left(\overline{u_1-u_2}\right)}{\partial\rho}\left(u_1-u_2\right)d\sigma_{\mathbb{B}^n}\right).
\end{split}\end{equation} when $R\rightarrow +\infty$.

Applying the Green formula \eqref{eq:green-gen1} and \eqref{eq:green-gen2}, using the boundary condition $u_1-u_2|_{\partial \Omega}=0$, one has
\begin{equation}\begin{split}
&\int_{\partial B_{\mathbb{B}^n}(0,R)}\frac{\partial \overline{u}}{\partial\nu_{\mathbb{B}^n}}\left(u_1-u_2\right)d\sigma_{\mathbb{B}^n}\\
&\ \ =\int_{B_{\mathbb{B}^n}(0,R)\setminus\Omega}\Delta_{\mathbb{B}^n}(\left(\overline{u_1-u_2}\right))\left(u_1-u_2\right)dV_{\mathbb{B}^n}\\
&\ \ \ \ +\int_{B_{\mathbb{B}^n}(0,R)\setminus\Omega}|\nabla_{\mathbb{B}^n}\left(u_1-u_2\right)|^2dV_{\mathbb{B}^n}
+\int_{\partial \Omega}\frac{\partial \left(\overline{u_1-u_2}\right)}{\partial\nu_{\mathbb{B}^n}}\left(u_1-u_2\right)d\sigma\\
&\ \ =-\int_{B_{\mathbb{B}^n}(0,R)\setminus\Omega}\left(\frac{(n-1)^2}{4}+\mu_0^2\right)|u_1-u_2|^2dV_{\mathbb{B}^n}+\int_{B_{\mathbb{B}^n}(0,R)\setminus\Omega}|\nabla_{\mathbb{B}^n}\left(u_1-u_2\right)|^2dV_{\mathbb{B}^n},\\
\end{split}\end{equation}
which implies
$$\Im\left(\tanh(\frac{R}{2})\int_{\partial B_{\mathbb{B}^n}(0,R_0)}\frac{\partial \left(\overline{u_1-u_2}\right)}{\partial\rho}\left(u_1-u_2\right)d\sigma_{\mathbb{B}^n}\right)=0.$$  Combining this and \eqref{addeq1}, we conclude that $$\lim\limits_{R\rightarrow +\infty}\int_{\partial B_{\mathbb{B}^n}(0,R)}|u_1-u_2|^2d\sigma_{\mathbb{B}^n}=0.$$ Applying the hyperbolic Rellich theorem, we derive that $u_1-u_2$ is equal to zero in $\mathbb{B}^n\setminus \Omega$. This accomplishes
the proof of uniqueness. Now, we start to show the existence. Since $u(x)=u^s(x)+e_{-2\mu, \xi}(x)$, we only need to show that the existence of the scattered wave $u^s$. We need the following lemma.
\begin{lem}
The following scattering problem
\begin{equation}\label{sca1}\begin{cases}
&-\Delta_{\mathbb{B}^n}u^s-\frac{(n-1)^2}{4}u^s=\mu_0^2u^s,\ \ x\in \mathbb{B}^n\setminus \Omega.\\
&u^s=-e_{2\mu, \xi}(x),\ \ x\in \partial \Omega\\
& \frac{\partial u_s}{\partial \rho}-(\mu_0 \mathrm{i}-\frac{n-1}{2})\tanh(\frac{\rho}{2})u_s=o\left(\frac{1}{\sinh^{\frac{n-1}{2}}(\rho)}\right)
\end{cases}\end{equation}
is solvable and it is equivalent to the following truncated system: find $u_1$ such that
\begin{equation}\label{trun}\begin{cases}
&-\Delta_{\mathbb{B}^n}u_1-\frac{(n-1)^2}{4}u_1=\mu_0^2u_1,\ \ x\in B_{\mathbb{B}^n}(0,R_0)\setminus \Omega.\\
&u_1=-e_{2\mu, \xi}(x),\ \ x\in \partial \Omega\\
& \frac{\partial u_1}{\partial \nu_{\mathbb{B}^n}}=\Lambda u_1,\ \ x\in \partial B_{\mathbb{B}^n}(0,R_0)
\end{cases}\end{equation}
where $\Lambda: H^{\frac{1}{2}}(\partial B_{\mathbb{B}^n}(0,R_0))\rightarrow H^{-\frac{1}{2}}(\partial B_{\mathbb{B}^n}(0,R_0))$ is the Dirichlet-to-Neuman map defined by $\Lambda \psi=\frac{\partial W}{\partial \nu_{\mathbb{B}^n} }|_{\partial B_{\mathbb{B}^n}(0,R_0)}$ with $W\in H^{1}_{loc}(\mathbb{B}^n\setminus \overline{B_{\mathbb{B}^n}(0,R_0)})$ being the unique solution to the equation
\begin{equation*}\label{DNM}\begin{cases}
&-\Delta_{\mathbb{B}^n}W-\frac{(n-1)^2}{4}W=\lambda_0^2W,\ \ x\in \mathbb{B}^n \setminus B_{\mathbb{B}^n}(0,R_0).\\
&W=\psi\in H^{\frac{1}{2}}(\partial B_{\mathbb{B}^n}(0,R_0))\\
& \frac{\partial W}{\partial \rho}-(\mu_0 \mathrm{i}-\frac{n-1}{2})\tanh(\frac{\rho}{2})W=o\left(\frac{1}{\sinh^{\frac{n-1}{2}}(\rho)}\right)
\end{cases}\end{equation*}
\end{lem}
\begin{proof}
By the definition of $\Lambda$, we see that if $u^s$ solves the equation \eqref{sca1}, then $u_1=u^s|_{B_{\mathbb{B}^n}(0,R_0)\setminus\Omega}$ is the solution of equation \eqref{trun}. On the other hand, by applying the Green's representation formula to the solution $u_1$ of equation \eqref{trun}, we obtain that
\begin{equation*}\begin{split}
u_1(x)&=-\int_{\partial\Omega}\left(\frac{\partial u_1(y)}{\partial \nu_{\mathbb{B}^n}(y)}G_{-\mu_0 \mathrm{i}}(\rho(x,y))-u_1(y)\frac{\partial G_{-\mu_0 \mathrm{i}}(\rho(x,y))}{\partial \nu_{\mathbb{B}^n}(y)}\right)d\sigma_{\mathbb{B}^n}(y)\\
&\ \ +\int_{\partial B_{\mathbb{B}^n}(0,R)}\left(\Lambda u_1(y)G_{-\mu_0 \mathrm{i}}(\rho(x,y))-u_1(y)\frac{\partial G_{-\mu_0 \mathrm{i}}(\rho(x,y))}{\partial \nu_{\mathbb{B}^n}(y)}\right)d\sigma_{\mathbb{B}^n}(y)
\end{split}\end{equation*}
for $x\in B_{\mathbb{B}^n}(0,R_0)\setminus\overline{\Omega}$. Now, we claim that $$\int_{\partial B_{\mathbb{B}^n}(0,R_0)}\left(\Lambda u_1(y)G_{-\mu_0 \mathrm{i}}(\rho(x,y))-u_1(y)\frac{\partial G_{-\mu_0 \mathrm{i}}(\rho(x,y))}{\partial \nu_{\mathbb{B}^n}(y)}\right)d\sigma_{\mathbb{B}^n}(y)=0.$$

For $x\in B_{\mathbb{B}^n}(0,R_0)$, since both $G_{-\mu_0 \mathrm{i}}(\rho(x,.))$ and $W(.)$ satisfy equation
$$-\Delta_{\mathbb{B}^n}u-\frac{(n-1)^2}{4}u=\mu_0^{2}u,\ \ y\in \mathbb{B}^n\setminus B_{\mathbb{B}^n}(0,R_0).$$ Green representation and the definition of $\Lambda$ gives that
\begin{equation*}\begin{split}
&\int_{\partial B_{\mathbb{B}^n}(0,R_0)}\left(\Lambda u_1(y)G_{-\mu_0 \mathrm{i}}(\rho(x,y))-u_1(y)\frac{\partial G_{-\mu \mathrm{i}}(\rho(x,y))}{\partial \nu_{\mathbb{B}^n}(y)}\right)d\sigma_{\mathbb{B}^n}(y)\\
&\ \ =\int_{\partial B_{\mathbb{B}^n}(0,R_0)}\left( \frac{\partial W(y)}{\partial \nu_{\mathbb{B}^n}(y)}G_{-\mu_0 \mathrm{i}}(\rho(x,y))-W(y)\frac{\partial G_{-\mu_0 \mathrm{i}}(\rho(x,y))}{\partial \nu_{\mathbb{B}^n}(y)}\right)d\sigma_{\mathbb{B}^n}(y)\\
&\ \ =\lim\limits_{R\rightarrow +\infty}\int_{\partial B_{\mathbb{B}^n}(0,R)}\left( \frac{\partial W(y)}{\partial \nu_{\mathbb{B}^n}(y)}G_{-\mu \mathrm{i}}(\rho(x,y))-W(y)\frac{\partial G_{-\mu_0 \mathrm{i}}(\rho(x,y))}{\partial \nu_{\mathbb{B}^n}(y)}\right)d\sigma_{\mathbb{B}^n}(y)\\
&\ \ =\lim\limits_{R\rightarrow +\infty}\int_{\partial B_{\mathbb{B}^n}(0,R)}\left( \frac{\partial W(y)}{\partial \rho}G_{-\mu_0 \mathrm{i}}(\rho(x,y))-W(y)\frac{\partial G_{-\mu_0 \mathrm{i}}(\rho(x,y))}{\partial \rho}\right)d\sigma_{\mathbb{B}^n}(y)\\
&\ \ =\lim\limits_{R\rightarrow +\infty} \int_{\partial B_{\mathbb{B}^n}(0,R)}\left(\frac{\partial W}{\partial \rho}-(\mu \mathrm{i}-\frac{n-1}{2})\tanh(\frac{\rho}{2})W\right)G_{-\mu_0 \mathrm{i}}d\sigma_{\mathbb{B}^n}\\
&\ \ \ \ \ -\lim\limits_{R\rightarrow +\infty}\left(\frac{\partial G_{-\mu_0 \mathrm{i}}}{\partial \rho}-(\mu_0 \mathrm{i}-\frac{n-1}{2})\tanh(\frac{\rho}{2})G_{-\mu_0 \mathrm{i}}\right)Wd\sigma_{\mathbb{B}^n}\\
&\ \ \leq \lim\limits_{R\rightarrow +\infty}\left(\int_{\partial B_{\mathbb{B}^n}(0,R)}\left|\left(\frac{\partial W}{\partial \rho}-(\mu_0 \mathrm{i}-\frac{n-1}{2})\tanh(\frac{\rho}{2})\right)W\right|^2d\sigma_{\mathbb{B}^n}\right)^{\frac{1}{2}}\left(\int_{\partial B_{\mathbb{B}^n}(0,R)}|G_{-\mu_0 \mathrm{i}}|^2d\sigma_{\mathbb{B}^n}\right)^{\frac{1}{2}}\\
&\ \ \ \ +\lim\limits_{R\rightarrow +\infty}\left(\int_{\partial B_{\mathbb{B}^n}(0,R)}\left|\left(\frac{\partial G_{-\mu_0 \mathrm{i}}}{\partial \rho}-(\mu_0 \mathrm{i}-\frac{n-1}{2})\tanh(\frac{\rho}{2})\right)G_{-\mu_0 \mathrm{i}}\right|^2d\sigma_{\mathbb{B}^n}\right)^{\frac{1}{2}}\left(\int_{\partial B_{\mathbb{B}^n}(0,R)}|W|^2d\sigma_{\mathbb{B}^n}\right)^{\frac{1}{2}}\\
\end{split}\end{equation*}
Since $W$ and $G_{-\mu_0 \mathrm{i}}$ satisfy hyperbolic radiation condition, then we deduce that
$$\lim\limits_{R\rightarrow +\infty}\int_{\partial B_{\mathbb{B}^n}(0,R)}\left|\left(\frac{\partial W}{\partial \rho}-(\mu_0 \mathrm{i}-\frac{n-1}{2})\tanh(\frac{\rho}{2})\right)W\right|^2d\sigma_{\mathbb{B}^n}=0$$
and $$\lim\limits_{R\rightarrow +\infty}\int_{\partial B_{\mathbb{B}^n}(0,R)}\left|\left(\frac{\partial G_{-\mu_0 \mathrm{i}}}{\partial \rho}-(\mu_0 \mathrm{i}-\frac{n-1}{2})\tanh(\frac{\rho}{2})\right)G_{-\mu_0 \mathrm{i}}\right|^2d\sigma_{\mathbb{B}^n}=0.$$
Combining the above estimate, we conclude that $$\int_{\partial B_{\mathbb{B}^n}(0,R)}\left(\Lambda u_1(y)G_{-\mu_0 \mathrm{i}}(\rho(x,y))-u_1(y)\frac{\partial G_{-\mu_0 \mathrm{i}}(\rho(x,y))}{\partial \nu_{\mathbb{B}^n}(y)}\right)d\sigma_{\mathbb{B}^n}(y)=0.$$
Hence, $$u(x)=-\int_{\partial\Omega}\left(\frac{\partial u_1(y)}{\partial \nu_{\mathbb{B}^n}(y)}G_{-\mu_0 \mathrm{i}}(\rho(x,y))-u_1(y)\frac{\partial G_{-\mu_0 \mathrm{i}}(\rho(x,y))}{\partial \nu_{\mathbb{B}^n}(y)}\right)d\sigma_{\mathbb{B}^n}(y).$$
It is easy to check that $u_1$ can be readily extended to an $H^{1}_{loc}(\mathbb{B}^n\setminus B_{\mathbb{B}^n}(0,R))$ function, which together with the uniqueness of solution of equation \eqref{sca1} implies that $u_1=u_s$. This proves the equivalence of equation \eqref{sca1} and equation \eqref{trun}. The existence of equation \eqref{sca1} follows from standard Fredholm alternative theorem if we choose $B_{R_0}$  such that $\mu_0^2$ is not the eigenvalue of the operator $-\Delta_{\mathbb{B}^n}-\frac{(n-1)^2}{4}$ in $B_{\mathbb{B}^n}(0,R_0)\setminus \Omega$ with the Dirichlet boundary on $\partial \Omega$ and Neuman boundary condition on $\partial B_{\mathbb{B}^n}(0,R)$. 
\end{proof}

Through the above lemma, we prove that there exists a unique total field wave $u\in H^1_{loc}(\mathbb{B}^n)$ satisfies equation \eqref{obs}. Then we accomplish the proof of Theorem \ref{obstacle1}.
\end{proof}

\begin{rmk}\label{rem4.5}
Recall that there exists a unique solution $u_s\in H_{loc}^1(\mathbb{B}^n)$ to the scattering system \eqref{sca1}, we claim that $u$ can be written as far-field pattern form. Since $u^s$ satisfies hyperbolic radiation condition, applying Green representative formula, for any $x\in \mathbb{B}^n\setminus B_{\mathbb{B}^n}(0,R)$, we can write
\begin{equation*}\begin{split}
& u^s=\int_{\partial B_{\mathbb{B}^n}(0,R)}\left(\frac{\partial u^s}{\partial \nu_{\mathbb{B}^n}(y)}G_{-\mu \mathrm{i}}(\rho(x,y))-\frac{\partial G_{-\mu \mathrm{i}}(\rho(x,y))}{\partial \nu_{\mathbb{B}^n}(y)}u^s\right)d\sigma_{\mathbb{B}^n}(y)\\
&\ \ + \lim\limits_{t\rightarrow +\infty}\int_{\partial B_{\mathbb{B}^n}(0,t)}\left(\frac{\partial u^s}{\partial \nu_{\mathbb{B}^n}(y)}G_{-\mu \mathrm{i}}(\rho(x,y))-\frac{\partial G_{-\mu \mathrm{i}}(\rho(x,y))}{\partial \nu_{\mathbb{B}^n}(y)}u^s\right)d\sigma_{\mathbb{B}^n}(y)\\
&=\int_{\partial B_{\mathbb{B}^n}(0,R)}\left(\frac{\partial u^s}{\partial \rho(y)}G_{-\mu \mathrm{i}}(\rho(x,y))-\frac{\partial G_{-\mu \mathrm{i}}(\rho(x,y))}{\partial \rho(y)}u^s\right)d\sigma_{\mathbb{B}^n}(y)\\
&\ \ + \lim\limits_{t\rightarrow +\infty}\int_{\partial B_{\mathbb{B}^n}(0,t)}\left(\frac{\partial u^s}{\partial \rho(y)}G_{-\mu \mathrm{i}}(\rho(x,y))-\frac{\partial G_{-\mu \mathrm{i}}(\rho(x,y))}{\partial \rho(y)}u^s\right)d\sigma_{\mathbb{B}^n}(y).
\end{split}\end{equation*}
Recall that the hyperbolic radiation condition implies that
$$\lim\limits_{t\rightarrow +\infty}\int_{\partial B_{\mathbb{B}^n}(0,t)}|u^s|^2d\sigma_{\mathbb{B}^n}(y)=O(1),$$
hence \small{\begin{equation*}\begin{split}
&\lim\limits_{t\rightarrow +\infty}\int_{\partial B_{\mathbb{B}^n}(0,t)}\left(\frac{\partial u^s}{\partial \rho(y)}G_{-\mu \mathrm{i}}(\rho(x,y))-\frac{\partial G_{-\mu \mathrm{i}}(\rho(x,y))}{\partial \rho(y)}u^s\right)d\sigma_{\mathbb{B}^n}(y)\\
&\ \ =\lim\limits_{t\rightarrow +\infty}\int_{\partial B_{\mathbb{B}^n}(0,t)}\left(\left(\frac{\partial u^s}{\partial \rho(y)}-(\mu \mathrm{i}-\frac{n-1}{2})\tanh(\frac{\rho}{2})u^s\right)G_{-\mu \mathrm{i}}-\left(\frac{\partial G_{-\mu \mathrm{i}}}{\partial \rho(y)}-(\mu \mathrm{i}-\frac{n-1}{2})\tanh(\frac{\rho}{2})G_{-\mu \mathrm{i}}\right)u^s\right)d\sigma_{\mathbb{B}^n}(y)\\
&\ \ \leq \lim\limits_{t\rightarrow +\infty}\left(\int_{\partial B_{\mathbb{B}^n}(0,t)}\left|\left(\frac{\partial u^s}{\partial \rho(y)}-(\mu \mathrm{i}-\frac{n-1}{2})\tanh(\frac{\rho}{2})\right)u^s\right|^2d\sigma_{\mathbb{B}^n}(y)\right)^{\frac{1}{2}}\left(\int_{\partial B_{\mathbb{B}^n}(0,t)}|G_{-\mu \mathrm{i}}|^2d\sigma_{\mathbb{B}^n}(y)\right)^{\frac{1}{2}}\\
&\ \ \ \ +\lim\limits_{t\rightarrow +\infty}\left(\int_{\partial B_{\mathbb{B}^n}(0,t)}\left|\left(\frac{\partial G_{-\mu \mathrm{i}}}{\partial \rho(y)}-(\mu \mathrm{i}-\frac{n-1}{2})\tanh(\frac{\rho}{2})\right)G_{-\mu \mathrm{i}}\right|^2d\sigma_{\mathbb{B}^n}(y)\right)^{\frac{1}{2}}\left(\int_{\partial B_{\mathbb{B}^n}(0,t)}|u^s|^2d\sigma_{\mathbb{B}^n}(y)\right)^{\frac{1}{2}}\\
&\ \ =0.
\end{split}\end{equation*}}
Combining the above estimate, we derive that
$$ u^s=\int_{\partial B_{\mathbb{B}^n}(0,R)}\left(\frac{\partial u^s}{\partial \rho(y)}G_{-\mu \mathrm{i}}(\rho(x,y))-\frac{\partial G_{-\mu \mathrm{i}}(\rho(x,y))}{\partial \rho(y)}u^s\right)d\sigma_{\mathbb{B}^n}(y).$$
For $y\in \partial B_{\mathbb{B}^n}(0,R)$, when $\rho(x)\rightarrow +\infty$ (i.e. $\rho(x,y)\rightarrow +\infty$), we have
\begin{equation*}\begin{split}
G_{-\mu \mathrm{i}}(\rho(x,y))&=\frac{A_{n,-\mu \mathrm{i}}}{2^{\frac{n-1}{2}}}\frac{\left(\cosh(\frac{\rho(x)}{2})\right)^{2\mu \mathrm{i}}}{\cosh^{n-1}(\frac{\rho(x)}{2})}\frac{\left(\cosh(\frac{\rho(y)}{2})\right)^{2\mu \mathrm{i}}}{\cosh^{n-1}(\frac{\rho(y)}{2})}\\
&\ \ \times   \left(1-2\hat{x}\cdot y+|y|^2\right)^{\mu \mathrm{i}-\frac{n-1}{2}}\int_{0}^{\pi}(\sin t)^{-2\mu \mathrm{i}}dt+ O\left(\frac{1}{\sinh^{\frac{n+1}{2}}(\rho(x))}\right)\\
&=\frac{A_{n,-\mu \mathrm{i}}}{2^{\frac{n-1}{2}}}\frac{\left(\cosh(\frac{\rho(x)}{2})\right)^{2\mu \mathrm{i}}}{\cosh^{n-1}(\frac{\rho(x)}{2})}\frac{\left(\cosh(\frac{\rho(y)}{2})\right)^{2\mu \mathrm{i}}}{\cosh^{n-1}(\frac{\rho(y)}{2})}\\
&\ \ \times   \left(1-2\hat{x}\cdot \tanh(\frac{\rho(y)}{2})\hat{y}+\tanh^2(\frac{\rho(y)}{2})\right)^{\mu \mathrm{i}-\frac{n-1}{2}}\int_{0}^{\pi}(\sin t)^{-2\mu \mathrm{i}}dt+ O\left(\frac{1}{\sinh^{\frac{n+1}{2}}(\rho(x))}\right)\\
\end{split}\end{equation*}
and
\begin{equation*}\begin{split}
&\frac{\partial G_{-\mu \mathrm{i}}(\rho(x,y))}{\partial \rho(y)}\\
&\ \ =\frac{A_{n,-\mu \mathrm{i}}}{2^{\frac{n-1}{2}}}\frac{\left(\cosh(\frac{\rho(x)}{2})\right)^{2\mu \mathrm{i}}}{\cosh^{n-1}(\frac{\rho(x)}{2})}\int_{0}^{\pi}(\sin t)^{-2\mu \mathrm{i}}dt \\
&\ \  \ \ \times \frac{\partial}{\partial\rho(y)}
\left(\frac{\left(\cosh(\frac{\rho(y)}{2})\right)^{2\mu \mathrm{i}}}{\cosh^{n-1}(\frac{\rho(y)}{2})}
\left(1-2\hat{x}\cdot \tanh(\frac{\rho(y)}{2})+|\tanh(\frac{\rho(y)}{2})|^2\right)^{\mu \mathrm{i}-\frac{n-1}{2}}\right)+ O\left(\frac{1}{\sinh^{\frac{n+1}{2}}(\rho(x))}\right)
\end{split}\end{equation*}

Then
\begin{equation*}\begin{split}
u^s&=\frac{\left(\cosh(\frac{\rho(x)}{2})\right)^{2\mu \mathrm{i}}}{\cosh^{n-1}(\frac{\rho(x)}{2})}\int_{0}^{\pi}(\sin t)^{-2\mu \mathrm{i}}dt\times \frac{A_{n,-\mu \mathrm{i}}}{2^{\frac{n-1}{2}}}\int_{\partial B_{\mathbb{B}^n}(0,R)}\frac{\left(\cosh(\frac{\rho(y)}{2})\right)^{2\mu \mathrm{i}}}{\cosh^{n-1}(\frac{\rho(y)}{2})}\left(1-2\hat{x}\cdot y+|y|^2\right)^{\mu \mathrm{i}-\frac{n-1}{2}}\frac{\partial u^s}{\partial \rho(y)}d\sigma_{\mathbb{B}^n}(y)\\
&\ \ +\frac{A_{n,-\mu \mathrm{i}}}{2^{\frac{n-1}{2}}}\frac{\left(\cosh(\frac{\rho(x)}{2})\right)^{2\mu \mathrm{i}}}{\cosh^{n-1}(\frac{\rho(x)}{2})}\int_{0}^{\pi}(\sin t)^{-2\mu \mathrm{i}}dt \\
&\ \  \ \ \times \int_{\partial B_{\mathbb{B}^n}(0,R)}\frac{\partial}{\partial\rho(y)}
\left(\frac{\left(\cosh(\frac{\rho(y)}{2})\right)^{2\mu \mathrm{i}}}{\cosh^{n-1}(\frac{\rho(y)}{2})}
\left(1-2\hat{x}\cdot \tanh(\frac{\rho(y)}{2})+|\tanh(\frac{\rho(y)}{2})|^2\right)^{\mu \mathrm{i}-\frac{n-1}{2}}\right)d\sigma_{\mathbb{B}^n}(y)\\
&\ \ \ \ +O\left(\frac{1}{\sinh^{\frac{n+1}{2}}(\rho(x))}\right)\\
&=\frac{\left(\cosh(\frac{\rho(x)}{2})\right)^{2\mu \mathrm{i}}}{\cosh^{n-1}(\frac{\rho(x)}{2})} u_{\infty}^{\mu}(\hat{x}, \xi)+O\left(\frac{1}{\sinh^{\frac{n+1}{2}}(\rho(x))}\right),
\end{split}\end{equation*}
where
\begin{equation*}\begin{split}
u_{\infty}^{\mu}(\hat{x}, \xi)&=\left(\int_{0}^{\pi}(\sin t)^{-2\mu \mathrm{i}}dt\times \frac{A_{n,-\mu \mathrm{i}}}{2^{\frac{n-1}{2}}}\right)\left(\int_{\partial B_{\mathbb{B}^n}(0,R)}\frac{\left(\cosh(\frac{\rho(y)}{2})\right)^{2\mu \mathrm{i}}}{\cosh^{n-1}(\frac{\rho(y)}{2})}\left(1-2\hat{x}\cdot y+|y|^2\right)^{\mu \mathrm{i}-\frac{n-1}{2}}\frac{\partial u^s}{\partial \rho(y)}d\sigma_{\mathbb{B}^n}(y)\right.\\
&\ \ + \left.\int_{\partial B_{\mathbb{B}^n}(0,R)}\frac{\partial}{\partial\rho(y)}
\left(\frac{\left(\cosh(\frac{\rho(y)}{2})\right)^{2\mu \mathrm{i}}}{\cosh^{n-1}(\frac{\rho(y)}{2})}
\left(1-2\hat{x}\cdot \tanh(\frac{\rho(y)}{2})+|\tanh(\frac{\rho(y)}{2})|^2\right)^{\mu \mathrm{i}-\frac{n-1}{2}}\right)d\sigma_{\mathbb{B}^n}(y)\right).
\end{split}\end{equation*}
In summary, the total wave $u$ can be written as
$$u(x)=e_{-2\mu, \xi}(x)+\frac{\left(\cosh(\frac{\rho(x)}{2})\right)^{2\mu \mathrm{i}}}{\cosh^{n-1}(\frac{\rho(x)}{2})}u_{\infty}^{\mu}(\hat{x}, \xi)+O\left(\frac{1}{\sinh^{\frac{n+1}{2}}(\rho(x))}\right).$$
\end{rmk}

Finally, we note that in Theorems~\ref{thm2}, \ref{thm4}, and~\ref{obstacle1}, we considered scattering from compactly supported sources, potentials, and obstacles in the hyperbolic space $\mathbb{B}^n$, and established the well-posedness of our scattering theory. It is worth emphasizing that, owing to the local nature of the hyperbolic Sommerfeld radiation condition, the theory can be readily extended to more general asymptotically hyperbolic manifolds (cf. \cite{JS, Mel}).

\section{Hyperbolic Fourier approximation}\label{sec:fourier-approximation}

\subsection{Density theorem for Helmholtz equation in hyperbolic space}
The eigenvalue function $e_{-2\mu_0, \xi}(x)$ is a smooth solution of the Helmholtz equation in hyperbolic space $$\big(-\Delta_{\mathbb{B}^n}-\frac{(n-1)^2}{4}\big)v=\mu^2 v,\ \ x\in \mathbb{B}^n,$$
where $\mu\geq 0$.
For any $g(\xi)\in L^{2}(\mathbb{S}^{n-1})$, it is also easy to check that the function $$u(x)\triangleq\int_{S^{n-1}}e_{-2\mu, \xi}(x)g(\xi)d\xi$$ belongs to $C^{\infty}(\mathbb{B}^n)$ and also satisfies the Helmholtz equation in hyperbolic space. We define
$$W_{\mathbb{B}^n}^{\mu}(\Omega)=\{u|_{\Omega}\ :\ u(x)=\int_{\mathbb{S}^{n-1}}e_{-2\mu, \xi}(x)g(\xi)d\xi,\ g(\xi)\in L^{2}(\mathbb{S}^{n-1})\},$$
where $\Omega$ is a some domain in $\mathbb{B}^n$.
We also introduce the space $\widetilde{W_{\mathbb{B}^n}^{\mu}}(\Omega)$ by
$$\widetilde{W_{\mathbb{B}^n}^{\mu}}(\Omega)=\{u\in C^{\infty}(\Omega)\ |\ \big(-\Delta_{\mathbb{B}^n}-\frac{(n-1)^2}{4}\big)u=\mu^2 u\}.$$
Obviously, $\widetilde{W_{\mathbb{B}^n}}^{\mu}(\Omega)$ is a closed subspace of $C^{\infty}(\Omega)$ with respect to the topology of $C^{\infty}(\Omega)$.
We prove the Herglotz approximation theorem in hyperbolic space.
\begin{thm}\label{density1}
$W_{\mathbb{B}^n}^{\mu}(\Omega)$ is dense in $\widetilde{W_{\mathbb{B}^n}^{\mu}}(\Omega)$ with respect to the topology of $C^{\infty}(\Omega)$.
\end{thm}
\begin{proof}
For any $u\in \widetilde{W_{\mathbb{B}^n}^{\mu}}(\Omega)$, we need to prove that there exists $u_n\in W_{\mathbb{B}^n}^{\mu}(\Omega)$ such that
$$\lim_{n\rightarrow +\infty}\|u_n-u\|_{C^k(D_j)}=0,$$ for any fixed $k$,
where $D_j$ is increasing and satisfies $\bigcup_{j=1}^{\infty}D_j=\Omega$.
According to the elliptic regularity, it is equivalent to prove that
$$\lim_{n\rightarrow +\infty}\int_{D_j}|u_n-u|^2dV_{\mathbb{B}^n}=0,$$ for any fixed $k$.
Let $V_{w}^{\mu}(D_j)$ denote the completion of $W_{\mathbb{B}^n}^{\mu}(D_j)$ under the $L^2(D_j)$ norm. We only need to prove that $u\in V_{w}(D_j)$.
According to the orthogonal decomposition theorem, it is equivalent to prove that if
$$\int_{D_j}uvdV_{\mathbb{B}^n}=0,\ \ \forall\ v\in W_{\mathbb{B}^n}^{\mu}(D_j),$$ then $u$ must be equal to zero in $D_j$. Define $u'=u\chi_{D_j}$ and
careful computation gives that
\begin{equation}\begin{split}
\int_{\mathbb{B}^n}u'vdV_g&=\int_{\mathbb{B}^n}u'(x)e_{-2\mu, \xi}(x)g(\xi)d\sigma_{\xi} dV_{\mathbb{B}^n}(x)\\
&=\int_{\mathbb{S}^{n-1}}\int_{\mathbb{B}^n}u'(x)e_{-2\mu, \xi}(x)dV_{\mathbb{B}^n}(x)d\sigma_{\xi}\\
&=\int_{\mathbb{S}^{n-1}} \widehat{u}(2\mu, \xi)g(\xi)d\sigma_{\xi}.
\end{split}\end{equation}
This together with $\int_{D_j}uvdV_g=0$ yields that $\widehat{u'}(2\mu, \xi)=0$.
Since $$\widehat{u'}(2\lambda, \xi)=\int_{\mathbb{B}^n}u'(x)\big(\frac{\sqrt{1-|x|^2}}{|x-\xi|}\big)^{n-1-2i\lambda}dV=$$ and $u'$ has the compact support, it is not difficult to check that $\widehat{u'}(2\lambda, \xi)\in C^{\infty}(\mathbb{R}\times \mathbb{S}^{n-1})$ and
$|\lambda|^m \widehat{u'}(2\lambda,\xi)\in L^{\infty}(\mathbb{R}\times \mathbb{S}^{n-1})$ for any $m>0$. We introduce the new function $v$ by
$$v(x)=D_{n}\int_{-\infty}^{\infty}\int_{\mathbb{S}^{n-1}}(\lambda^2-\mu^2)^{-1}\widehat{u'}(\lambda,\xi)\big(\frac{\sqrt{1-|x|^2}}{|x-\xi|}\big)^{n-1-2i\lambda}|c(\lambda)|^{-2}d\lambda d\sigma_{\xi}.$$ and the $v$ is well defined since $(\lambda^2-\mu^2)^{-1}\widehat{u'}(2\lambda,\xi)\in C^{\infty}(\mathbb{R}\times \mathbb{S}^{n-1})$. By the property of Fourier transform in hyperbolic space, we see that $v$ satisfies the Helmholtz equation in hyperbolic space
\begin{equation}\label{adeq1}
\big(-\Delta_{\mathbb{B}^n}-\frac{(n-1)^2}{4}-\mu^2\big)v=u'=0,\ \ x\in \mathbb{B}^n\setminus \overline{D_j}.
\end{equation}
Now, we will show that $v(x)$ satisfies the condition of hyperbolic Rellich's type theorem. Denote by $h(\lambda, \xi)=(\lambda^2-\mu^2)^{-1}\widehat{u'}(2\lambda,\xi)$ and Careful computations give that
\begin{equation}\begin{split}
&\frac{v(x)}{(1-|x|^2)^{\frac{n-1}{2}}}|\ln(1-|x|^2)|\\
&\ \ \lesssim D_{n}\int_{-\infty}^{\infty}\int_{\mathbb{S}^{n-1}}h(\lambda, \xi)(\Delta_{\mathbb{R}^n})_{\xi}(|x-\xi|^{-{n-2}})|c(\lambda)|^{-2}\frac{d}{d\lambda}(\sqrt{1-|x|^2})^{-2i\lambda}d\sigma_{\xi} d\lambda\\
&\ \ \lesssim\int_{-\infty}^{\infty}\int_{\mathbb{S}^{n-1}}(\Delta_{\mathbb{R}^n})_{\xi}\frac{d}{d\lambda}\big(h(\lambda, \xi)\big)|x-\xi|^{-{n-2}}|c(\lambda)|^{-2}(\sqrt{1-|x|^2})^{-2i\lambda}d\sigma_{\xi} d\lambda\lesssim 1.
\end{split}\end{equation}
This proves that
$$v(x)=o\big(\sinh^{-\frac{n-1}{2}} \rho(x))\big), \  {\rm as}\ \rho(x)\rightarrow +\infty.$$ According to Rellich theorem of hyperbolic space, we can immediately derive that $v$ must be equal to zero in $\mathbb{B}^n\setminus \overline{D_j}$. Now, we are in position to prove that $u=0$ in $D_j$. By Green formula, we immediately get $$0=\int_{\Omega}-(\Delta_{\mathbb{B}^n}-\frac{(n-1)^2}{4}-\mu^2)u\cdot vdV_{\mathbb{B}^n}-\int_{\Omega}-(\Delta_{\mathbb{B}^n}-\frac{(n-1)^2}{4}-\mu^2)v\cdot udV_{\mathbb{B}^n}$$
since $v$ and its derivative are equal to zero on the boundary $\partial\Omega$. This together with the $u$ satisfying equation $$\big(-\Delta_{\mathbb{B}^n}-\frac{(n-1)^2}{4}\big)u=\mu^2 u,\ \ x\in \mathbb{B}^n$$
 yields that $$\int_{\Omega}u'udV_{\mathbb{B}^n}=0.$$ Then we accomplish the proof of Theorem \ref{density1}.
\end{proof}

\subsection{Density theorem for time-harmonic acoustic waves in hyperbolic space}
We consider the inverse scattering problem for time-harmonic acoustic waves in an
inhomogeneous medium. The direct scattering problem we are
now concerned with is, given the refractive index
$$q(x)=q_1(x)+i q_2(x),$$
where $q(x)$ is piecewise continuous in $\mathbb{B}^n$ such that
$$V(x)=q(x)-1$$ has the compact support to determine $u$ such that
\begin{equation}\label{acous1}\begin{cases}
\left(-\Delta_{\mathbb{B}^n}-\frac{(n-1)^2}{4}\right)u-\mu^2qu=0,\ \ x\in \mathbb{B}^n\\
u(x)=e_{-2\mu, \xi}(x)+u^s(x)\\
\frac{\partial u^s}{\partial \rho}-\left(\mu \mathrm{i}-\frac{n-1}{2}\right)\tanh(\frac{\rho}{2})u^s=o\left(\frac{1}{\sinh^{\frac{n-1}{2}}(\rho)}\right).
\end{cases}\end{equation}
 Obviously $u^s$ satisfies equation
\begin{equation}\begin{cases}
\left(-\Delta_{\mathbb{B}^n}-\frac{(n-1)^2}{4}\right)u^s-k^2 u^s=k^2V(x)u(x),\ \ x\in \mathbb{B}^n\\
\frac{\partial u^s}{\partial \rho}-\left(k \mathrm{i}-\frac{n-1}{2}\right)\tanh(\frac{\rho}{2})u^s=o\left(\frac{1}{\sinh^{\frac{n-1}{2}}(\rho)}\right).
\end{cases}\end{equation}
Applying the Green-representative formula in Remark \ref{acouse}, we immediately derive that
$$u^s(x)=\mu^2\int_{\mathbb{B}^n}G_{-\mu \mathrm{i}}(\rho(x,y))V(y)u(y)dV_\mathbb{B}(y)$$
and $$u(x)=e_{-2\mu, \xi}(x)+\mu^2\int_{\mathbb{B}^n}G_{-\mu \mathrm{i}}(\rho(x,y))V(y)u(y)dV_\mathbb{B}(y).$$
Now, we state a density theorem.

\begin{thm}\label{density2}
Let $B$ and $D$ be two open hyperbolic balls centered at the origin and containing the support of $V=q-1$ such that $\bar{B}\subset D$. Then the set of total fields $\{u(.,\xi),\ \xi\in \mathbb{S}^{n-1}\}$ satisfying equation \eqref{acous1} is complete in the closure of
$$H:=\{v\in H^2(D),\ \ \left(-\Delta_{\mathbb{B}^n}-\frac{(n-1)^2}{4}\right)v-\mu^2q(x)v=0\ {\rm in}\ D\}$$
with respect to $L^2(B)$ norm.
\end{thm}
\begin{proof}
Consider the mapping $A:H\rightarrow H^2_{loc}(\mathbb{R}^n)$ defined by the volume potential
$$(Av)(x):=\int_{B}G_{-\mu \mathrm{i}}(\rho(x,y))\overline{\left(I+\mu^2T_m^{*}\right)^{-1}v(y)}dV_{\mathbb{B}^n}(y),\ \ x\in \mathbb{B}^n,$$
where $T_V^{*}: L^2(B)\rightarrow L^2(B)$ denotes the adjoint of the Lippman-Schwinger operator  $T_V: L^2(B)\rightarrow L^2(B)$ given by
$$T_V(u)=-\int_{B}G_{-\mu \mathrm{i}}(\rho(x,y))V(y)u(y)dV_{\mathbb{B}^n}(y),\ \ x\in \mathbb{B}^n.$$
Note that for $v\in H$, we have that $Av\in H^2(B)$. For the density
$$W:=\left(I+\mu^2T_{V}^{*}\right)^{-1}v$$ of this potential, we have
\begin{equation}\label{density}
v(x)=W(x)-\overline{\mu^2V(x)}\int_{B}\overline{G_{-\mu \mathrm{i}}(\rho(x,y))}V(y)dV_{\mathbb{B}^n}(y),\ \ x\in \mathbb{B}^n.
\end{equation}
It is not difficult to check that $Av$ satisfies equation
$$\left(-\Delta_{\mathbb{B}^n}-\frac{(n-1)^2}{4}\right)Av-\mu^2Av=\overline{V},\ \ x\in B.$$ This together with \eqref{density} gives that
$$\bar{v}=\left(-\Delta_{\mathbb{B}^n}-\frac{(n-1)^2}{4}\right)Av-\mu^2Av-\mu^2V(x)Av=\left(-\Delta_{\mathbb{B}^n}-\frac{(n-1)^2}{4}\right)Av-k^2qAv.$$
Multiplying this equation by $w\in H$ and then integrating, appying the Green's theorem and using the fact that both $Av$ and $w$ solves the Helmholtz equation of hyperbolic space in $D\setminus \overline{B}$, we deduce that
\begin{equation}\label{integral}\begin{split}
\int_{B}w\bar{v}dV_g(x)&=-\int_{B}w\left(\left(\Delta_{\mathbb{B}^n}+\frac{(n-1)^2}{4}\right)Av+\mu^2qAv\right)\\
&=\int_{\partial D}\left(Av\frac{\partial w}{\partial \nu}-w\frac{\partial Av}{\partial \nu}\right)d\sigma_{\mathbb{B}^n}.
\end{split}\end{equation}
Now, let $v\in \hat{H}$ i.e. $v\in L^2(B)$ is the $L^2$ limit of a sequence $v_j$ from $H$. Assume that
$$(u(\cdot, \xi),v)=\int_{B}\overline{v(x)}u(x)dV_{\mathbb{B}^n}(x)=0,\ \ \xi\in \mathbb{S}^{n-1}.$$
Then from the hyperbolic Lippmann-Schwinger equation, we obtain that
$$0=\left((I+\mu^2T_V)^{-1}e_{-2\mu, \xi},v\right)=\left(e_{-2\mu, \xi},(I+\mu^2T_V^{*})^{-1}v\right).$$
Then it follows that
\begin{equation}\begin{split}
(Av)_{\infty}^{\mu}(\hat{x})&=\frac{A_{n,-\mu \mathrm{i}}}{2^{\frac{n}{2}-\frac{1}{2}}}\int_{B}\frac{\left(\cosh(\frac{\rho(y)}{2})\right)^{2\mu \mathrm{i}}}{\cosh^{n-1}(\frac{\rho(y)}{2})}\left(1-2\hat{x}\cdot y+|y|^2\right)^{\mu i-\frac{n-1}{2}}\overline{V(y)}dV_{\mathbb{B}^n}(y)\\
&=C_{n,\mathrm{i}}\int_{B}e_{-2\mu, \hat{x}}\overline{V(y)}dV_{\mathbb{B}^n}(y)\\
&=0.
\end{split}\end{equation}
This together with hyperbolic Rellich Lemma implies that $Av=0$ in $\mathbb{B}^n\setminus B$. Since $I+\mu^2T_V^{*}$ has a bounded inverse, by the Cauchy-Schwartz inequality, the mapping $A$ is bounded from $L^2(B)$ into $H^2(K)$ for each compact set $K\subset \mathbb{B}^n\setminus B$. Hence, inserting $v_j$ into \eqref{integral} and passing to the limit $j\rightarrow +\infty$ yields
$$\int_{B}w\overline{v}dV_{\mathbb{B}^n}=0$$
for all $w\in H$. Now inserting $w=v_j$ in this equation and again passing to the limit $j\rightarrow \infty$, we obtain
$$\int_{B}|v|^2dV_{\mathbb{B}^n}=0.$$
Then we accomplish the proof of Theorem \ref{density2}.
\end{proof}

\section{Hyperbolic inverse obstacle scattering}\label{sec:inverse-obstacle}
Assume that $\Omega$ is a bounded domain in the hyperbolic space $\mathbb{B}^n$. Assume that $u$ satisfies equation
\begin{equation}\begin{cases}
&-\Delta_{\mathbb{B}^n}u-\frac{(n-1)^2}{4}u=\mu^2u,\ \ x\in \mathbb{B}^n\setminus \Omega\\
&\mathcal{B}u=0,\ \ x\in \partial \Omega\\
&\displaystyle{ \frac{\partial (u-u^i)}{\partial \rho}-(\mu \mathrm{i}-\frac{n-1}{2})\tanh\left(\frac{\rho}{2}\right)(u-u^i)=o\left(\frac{1}{\sinh^{\frac{n-1}{2}}(\rho)}\right)},
\end{cases}\end{equation}
where $\mathcal{B}$ denotes the boundary operator (Soft, hard or impedance obstacle). The existence and uniqueness of the scattering from obstacle has been proved in Theorem \ref{obstacle1}. Recall that we have shown that the total wave field $u$ can be written as
$$u(x)=e_{-2\mu, \xi}(x)+\frac{\left(\cosh(\frac{\rho(x)}{2})\right)^{2\mu \mathrm{i}}}{\cosh^{n-1}(\frac{\rho(x)}{2})}u_{\infty}^{\mu}(\hat{x}, \xi)+O\left(\frac{1}{\sinh^{\frac{n+1}{2}}(\rho(x))}\right).$$
in Remark \eqref{rem4.5}. The inverse obstacle problem in hyperbolic space addressed in this work concerns the simultaneous recovery of the obstacle $\Omega$ from the measurement of far-field pattern $u^{\mu}_{\infty}(\hat{x}, \xi):$
$$u^{\mu}_{\infty}(\hat{x}, \xi),\ \ (\hat{x},\xi)\in \mathbb{S}^{n-1}\times \mathbb{S}^{n-1}\rightarrow \Omega \subseteq \mathbb{B}^n.$$
We proved
\begin{thm}\label{inverse-obstacle}
If we know the far-field pattern $u^{\mu}_{\infty}(\hat{x}, \xi)$ for all incident direction $\xi$ and observable direction $\hat{x}$, then we can uniquely determine the obstacle $\Omega$.
\end{thm}
\begin{proof}
Without loss of generality, we only consider the soft case.
Assume that $u_1$ and $u_2$ satisfies equation
\begin{equation}\begin{cases}
&-\Delta_{\mathbb{B}^n}u_1-\frac{(n-1)^2}{4}u_1=\mu^2u_1,\ \ x\in \mathbb{B}^n\setminus \Omega_1\\
&u_1=0,\ \ x\in \partial \Omega_1\\
&\displaystyle{ \frac{\partial (u_1-u^i)}{\partial \rho}-(\mu \mathrm{i}-\frac{n-1}{2})\tanh\left(\frac{\rho}{2}\right)(u_1-u^i)=o\left(\frac{1}{\sinh^{\frac{n-1}{2}}(\rho)}\right)},
\end{cases}\end{equation}
and
\begin{equation}\begin{cases}
&-\Delta_{\mathbb{B}^n}u_2-\frac{(n-1)^2}{4}u_2=\mu^2u_2,\ \ x\in \mathbb{B}^n\setminus \Omega_2\\
&u_2=0,\ \ x\in \partial \Omega_2\\
&\displaystyle{ \frac{\partial (u_2-u^i)}{\partial \rho}-(\mu \mathrm{i}-\frac{n-1}{2})\tanh\left(\frac{\rho}{2}\right)(u_2-u^i)=o\left(\frac{1}{\sinh^{\frac{n-1}{2}}(\rho)}\right)},
\end{cases}\end{equation}
respectively. Let $u_{\infty,1}^{\mu}(\hat x,\xi)$ and $u_{\infty,2}^{\mu}(\hat x,\xi)$ denote the far field pattern of $u_1$ and $u_2$ respectively. Our aim is to show that if $u_{\infty,1}^{\mu}(\hat x,\xi)=u_{\infty,2}^{\mu}(\hat x,\xi)$ for any incident direction $\xi$ and observable direction $\hat x$, then $\Omega_1=\Omega_2$.
\medskip

For any $z\in \mathbb{B}^n\setminus (\Omega_1 \cup \Omega_2)$, let
 $\Omega_z$ denote the bounded domain satisfying $\overline{\Omega_1}\cup \overline{\Omega_2}\subseteq \Omega_z$ and $z\in \mathbb{B}^n\setminus \Omega_z$. Let $G_{-\mu \mathrm{i}}(\rho(x,z))$ denote the Green function of the hyperbolic Helmholtz operator $-\Delta_{\mathbb{B}^n}-\frac{(n-1)^2}{4}-\mu^2$ with singularity at $z$. Then $G_{-\mu \mathrm{i}}(\rho(x,z))$ is a smooth function in
 $\Omega_z$ and satisfying equation $$-\Delta_{\mathbb{B}^n}G_{-\mu \mathrm{i}}(\rho(x,z))-\frac{(n-1)^2}{4}G_{-\mu \mathrm{i}}(\rho(x,z))=\mu^2G_{-\mu \mathrm{i}}(\rho(x,z)),\ \ x\in \Omega_z.$$
  According to the Theorem \ref{density1}, we may assume that there exists $g_z^k(\xi) \in L^2(\mathbb{S}^{n-1})$ such that $$\int_{\mathbb{S}^{n-1}}g_z^k(\xi)e_{-2\mu, \xi}(x)d\xi$$ strongly converges to $G_{-\mu \mathrm{i}}(\rho(x,z))$ in $C^{\infty}(\Omega_z)$. Denote by $u^{i,k}_z(x)=\int_{\mathbb{S}^{n-1}}g_z^k(\xi)e_{-2\mu, \xi}(x)d\xi$. Assume that $u_{z,1}^k$ and $u_{z,2}^k$ satisfies
equation \begin{equation*}\begin{cases}\label{adeq1}
&-\Delta_{\mathbb{B}^n}u_{z,1}^k-\frac{(n-1)^2}{4}u_{z,1}^k=\mu^2u_{z,1}^k,\ \ x\in \mathbb{B}^n\setminus \Omega_1\\
&u_{z,1}^k=0,\ \ x\in \partial \Omega_1\\
&\displaystyle{ \frac{\partial (u_{z,1}^k-u^{i,k}_z)}{\partial \rho}-(\mu \mathrm{i}-\frac{n-1}{2})\tanh\left(\frac{\rho}{2}\right)(u_{z,1}^k-u^{i,k}_z)=o\left(\frac{1}{\sinh^{\frac{n-1}{2}}(\rho)}\right)},
\end{cases}\end{equation*}
and
\begin{equation*}\begin{cases}\label{adeq2}
&-\Delta_{\mathbb{B}^n}u_{z,2}^k-\frac{(n-1)^2}{4}u_{z,2}^k=\mu^2u_{z,2}^k,\ \ x\in \mathbb{B}^n\setminus \Omega_2\\
&u_{z,2}^k=0,\ \ x\in \partial \Omega_2\\
&\displaystyle{ \frac{\partial (u_{z,2}^k-u^{i,k}_z)}{\partial \rho}-(\mu \mathrm{i}-\frac{n-1}{2})\tanh\left(\frac{\rho}{2}\right)(u_{z,2}-u^{i,k}_z)=o\left(\frac{1}{\sinh^{\frac{n-1}{2}}(\rho)}\right)}.
\end{cases}\end{equation*}
Denote by $$u_{z,1}^{s,k}(x)=u_{z,1}^k(x)-u^{i,k}_z(x),\ \ u_{z,2}^{s,k}(x)=u_{z,2}^k(x)-u^{i,k}_z(x).$$
Recall Remark \ref{rem4.5}, we can similarly prove that
$$u_z^1(x)=u^{i,k}_z(x)+\frac{\left(\cosh(\frac{\rho(x)}{2})\right)^{2\mu \mathrm{i}}}{\cosh^{n-1}(\frac{\rho(x)}{2})}u_{\infty,z,1}^{k,\mu}(\hat{x})+O\left(\frac{1}{\sinh^{\frac{n+1}{2}}(\rho(x))}\right)$$
and
$$u_z^2(x)=u^{i,k}_z(x)+\frac{\left(\cosh(\frac{\rho(x)}{2})\right)^{2\mu \mathrm{i}}}{\cosh^{n-1}(\frac{\rho(x)}{2})}u_{\infty,z,2}^{k,\mu}(\hat{x} )+O\left(\frac{1}{\sinh^{\frac{n+1}{2}}(\rho(x))}\right),$$

where \begin{equation*}\begin{split}
&u_{\infty,z,1}^{k,\mu}(\hat{x})\\
&\ \ =\left(\int_{0}^{\pi}(\sin t)^{-2\mu \mathrm{i}}dt\times \frac{A_{n,-\mu \mathrm{i}}}{2^{\frac{n-1}{2}}}\right)\left(\int_{\partial B_{\mathbb{B}^n}(0,R)}\frac{\left(\cosh(\frac{\rho(y)}{2})\right)^{2\mu \mathrm{i}}}{\cosh^{n-1}(\frac{\rho(y)}{2})}\left(1-2\hat{x}\cdot y+|y|^2\right)^{\mu \mathrm{i}-\frac{n-1}{2}}\frac{\partial u_{z,1}^{s,k}}{\partial \rho(y)}d\sigma_{\mathbb{B}^n}(y)\right.\\
&\ \ \ \ + \left.\int_{\partial B_{\mathbb{B}^n}(0,R)}\frac{\partial}{\partial\rho(y)}
\left(\frac{\left(\cosh(\frac{\rho(y)}{2})\right)^{2\mu \mathrm{i}}}{\cosh^{n-1}(\frac{\rho(y)}{2})}
\left(1-2\hat{x}\cdot \tanh(\frac{\rho(y)}{2})+|\tanh(\frac{\rho(y)}{2})|^2\right)^{\mu \mathrm{i}-\frac{n-1}{2}}\right)d\sigma_{\mathbb{B}^n}(y)\right)
\end{split}\end{equation*}
and
\begin{equation*}\begin{split}
&u_{\infty,z,2}^{k,\mu}(\hat{x})\\
&\ \ =\left(\int_{0}^{\pi}(\sin t)^{-2\mu \mathrm{i}}dt\times \frac{A_{n,-\mu \mathrm{i}}}{2^{\frac{n-1}{2}}}\right)\left(\int_{\partial B_{\mathbb{B}^n}(0,R)}\frac{\left(\cosh(\frac{\rho(y)}{2})\right)^{2\mu \mathrm{i}}}{\cosh^{n-1}(\frac{\rho(y)}{2})}\left(1-2\hat{x}\cdot y+|y|^2\right)^{\mu \mathrm{i}-\frac{n-1}{2}}\frac{\partial u_{z,2}^{s,k}}{\partial \rho(y)}d\sigma_{\mathbb{B}^n}(y)\right.\\
&\ \ \ \ + \left.\int_{\partial B_{\mathbb{B}^n}(0,R)}\frac{\partial}{\partial\rho(y)}
\left(\frac{\left(\cosh(\frac{\rho(y)}{2})\right)^{2\mu \mathrm{i}}}{\cosh^{n-1}(\frac{\rho(y)}{2})}
\left(1-2\hat{x}\cdot \tanh(\frac{\rho(y)}{2})+|\tanh(\frac{\rho(y)}{2})|^2\right)^{\mu \mathrm{i}-\frac{n-1}{2}}\right)d\sigma_{\mathbb{B}^n}(y)\right)
\end{split}\end{equation*}

Since $u_{\infty,1}^{\mu}(\hat x, \xi)=u_{\infty,2}^{\mu}(\hat x, \xi)$ for any $\xi \in \mathbb{S}^{n-1}$, then it is not difficult to check that
$u_{\infty,z,1}^{k,\mu}(\hat{x})=u_{\infty,z,2}^{k,\mu}(\hat{x})$. By hyperbolic Rellich theorem, we deduce that
$u_{z,1}^{s,k}(x)=u_{z,2}^{s,k}(x)$ in $\mathbb{B}^n\setminus (\overline{\Omega_1}\cup \overline{\Omega_2})$. Then $$u_{z,1}^k(x)-u^{i,k}_z(x)|_{\partial \Omega_2}=u_{z,2}^k(x)-u^{i,k}_z(x)|_{\partial\Omega_2}.$$ In view of $G_{-\mu \mathrm{i}}(\rho(x,z))=u^{i}_z$ in $\Omega_z$ and $u_{z}^2|_{\partial \Omega_2}=0$, we derive that \begin{equation}\label{infinty}
u_{z,1}^k(x)-u^{i,k}_z|_{\partial \Omega_2}=u_{z,2}^k(x)-u^{i,k}_z|_{\partial \Omega_2}=-u^{i,k}_z|_{\partial \Omega_2}=-G_{-\mu \mathrm{i}}(\rho(x,z))|_{\partial \Omega_2}.
\end{equation}
Let $z\rightarrow z^{*}$ for some $z^{*}\in \partial \Omega_2 $, we derive that $$(u_{z,1}^k-u^{i,k}_z)(z^*)=(u_{z,2}^k-u^{i,k}_z)(z^*)=-G_{-\mu \mathrm{i}}(\rho(z^{*},z))\rightarrow -\infty$$ as $z$ approaches to $z^{*}$. On the other hand, if $\Omega_2 \neq \Omega_1$, we may choose $z^{*}\in \partial\Omega_2$ such that $d(z^{*}, \partial\Omega_1)=d_0>0$. Since $u_{z,1}^k(x)$ satisfies equation \eqref{adeq1}, we can derive that in $B_{\mathbb{B}^n}(0, R)\setminus\Omega_1$, $u_z^{1}-u_z^{i}$ satisfies equation
\begin{equation*}\begin{cases}\label{adeq3}
&-\Delta_{\mathbb{B}^n}(u_{z,1}^k-u^{i,k}_z)-\frac{(n-1)^2}{4}(u_{z,1}^k-u^{i,k}_z)=\mu^2(u_{z,1}^k-u^{i,k}_z),\ \ x\in B_\mathbb{B}^n(0, R)\setminus \Omega_1\\
&(u_{z,1}^k-u^{i,k}_z)(x)=-u^{i,k}_z(x)=G_{-\mu i}(\rho(x,z)),\ \ x\in \partial \Omega_1.\\
\end{cases}\end{equation*}
Observing that $G_{-\mu i}(\rho(x,z))|_{\partial \Omega_1}$ is bounded as $z\rightarrow z^{*}$ since $\lim\limits_{z\rightarrow z^{*}} d(z, \partial\Omega_1)=d(z^{*}, \partial\Omega_1)=d_0>0$ and $G_{-\mu i}(\rho(x,z))$ being bounded for $x$ away from $z$, applying standard elliptic estimate, we derive that $u_{z,1}^k-u^{i,k}_z$ is bounded in $L^{\infty}(B_{\mathbb{B}^n}(0, R))$ when $z$ approaches to $z^{*}$. Especially, this deduces that $(u_{z,1}^k-u^{i,k}_z)(z^*)$ is bounded when $z$ approaches to $z^{*}$, which is a contradiction with \eqref{infinty}. Then we accomplish the proof of Theorem \ref{inverse-obstacle}.
\end{proof}
\section{Hyperbolic inverse potential scattering}\label{sec:inverse-medium}
We consider the inverse scattering problem for time-harmonic acoustic waves in an
inhomogeneous medium in hyperbolic space. The direct scattering problem we are now
concerned with is, given the refractive index $q(x)$ with support of $1-q(x)$ contained in $\Omega \in \mathbb{B}^n$ and imaginary part being non-negative to determine $u$ such that
\begin{equation}\begin{cases}
&\Delta_{\mathbb{B}^n} u+\mu^2q(x)u=0,\ \ x\in \mathbb{B}^n,\\
&u^s=u-u^{i}\ \  {\rm satisfies\ SRC}\  \\
&u^{i}=e_{-2\mu, \xi}(x),\ \hat{d}\in \mathbb{S}^{n-1}.\\
\end{cases}
\end{equation}
The existence and uniqueness of the solution has been established in Theorem \ref{thm4} of Chapter 4:
$$u(x)=e_{-2\mu, \xi}(x)-\int_{\mathbb{B}^n}G_{-\mu \mathrm{i}}(\rho(x,y))\mu^2(1-q(y))u(y)dV_{\mathbb{B}^n}(y),\ \ x\in \mathbb{B}^n.$$
Furthermore, it was also shown that $u^s$ has the expression of far-field pattern:
$$u^s(x)=\frac{\left(\cosh(\frac{\rho(x)}{2})\right)^{2\mu \mathrm{i}}}{\cosh^{n-1}(\frac{\rho(x)}{2})} u^{\mu}_{\infty}(\hat{x}, \xi)+O\left(\frac{1}{\sinh^{\frac{n+1}{2}}(\rho(x))}\right)$$
with the far field pattern $u^s_{\infty}(\hat{x}, \mu, \xi)$ given by
$$u^{\mu}_{\infty}(\hat{x},\xi)=-\mu^2\frac{A_{n,-\mu \mathrm{i}}}{2^{\frac{n-1}{2}}}\int_{0}^{\pi}(\sin t)^{-2\mu \mathrm{i}}dt\int_{\mathbb{B}^n}\frac{\left(\cosh(\frac{\rho(y)}{2})\right)^{2\mu \mathrm{i}}}{\cosh^{n-1}(\frac{\rho(y)}{2})}\left(1-2\hat{x}\cdot y+|y|^2\right)^{\mu \mathrm{i}-\frac{n-1}{2}}(1-q(y))u(y)dV_{\mathbb{B}^n}(y).$$ The inverse medium problem in hyperbolic space addressed in this work concerns the simultaneous recovery of the medium $q(x)$ from the measurement of far-field pattern $u^s_{\infty}(\hat{x}, \mu, \xi):$
$$u^{\mu}_{\infty}(\hat{x},\xi),\ \ (\hat{x},\xi)\in \mathbb{S}^{n-1}\times \mathbb{S}^{n-1}\rightarrow q(x), \ x\in \Omega\subseteq \mathbb{B}^n.$$
We proved
\begin{thm}\label{inverse-medium}
If we know the far-field pattern $u^{\mu}_{\infty}(\hat{x},\xi)$ for all incident direction $\xi$ and observable direction $\hat{x}$, then we can uniquely determine $q(x)$.
\end{thm}

\begin{proof}
Assume that $q_1$ and $q_2$ are two refractive indices such that the support of $q_1-1$ and $q_2-1$ is contained in $B_{\mathbb{H}}(0,R)$.
Obviously,
$$u_1(x)=e_{-2\mu, \xi}(x)+\frac{\left(\cosh(\frac{\rho(x)}{2})\right)^{2\mu \mathrm{i}}}{\cosh^{n-1}(\frac{\rho(x)}{2})} u^{\mu}_{\infty,1}(\hat{x}, \xi)+O\left(\frac{1}{\sinh^{\frac{n+1}{2}}(\rho(x))}\right)$$
and
$$u_2(x)=e_{-2\mu, \xi}(x)+\frac{\left(\cosh(\frac{\rho(x)}{2})\right)^{2\mu \mathrm{i}}}{\cosh^{n-1}(\frac{\rho(x)}{2})} u^{\mu}_{\infty,2}(\hat{x}, \xi)+O\left(\frac{1}{\sinh^{\frac{n+1}{2}}(\rho(x))}\right).$$
We also assume that $$u^{\mu}_{\infty,1}(\hat{x},\xi)=u^{\mu}_{\infty,2}(\hat{x},\xi)$$ for all $(\hat{x},\xi)\in \mathbb{S}^{n-1}\times \mathbb{S}^{n-1}$. Then
$u_1(x)-u_2(x)$ satisfies equation
$$\big(\Delta_{\mathbb{B}^n}+\mu^2\big)(u_1-u_2)=0,\ \ x\in \mathbb{B}^n\setminus B_{\mathbb{H}}(0,R)$$ with
$u_1(x)-u_2(x)=O\left(\frac{1}{\sinh^{\frac{n+1}{2}}(\rho(x))}\right)$ when $\rho(x)\rightarrow +\infty$.
Then by Rellich-theorem in hyperbolic space, we have
$$u_1(x,\mu, \xi)=u_2(x,\mu, \xi)\ ,\ \  x \in \mathbb{B}^n\setminus B_{\mathbb{H}}(0,R)$$
for any $\xi\in \mathbb{S}^{n-1}$. Hence $u=u_1-u_2$ satisfies the boundary conditions
$$u=\frac{\partial u}{\partial \nu_{\mathbb{B}^n}}=0,\ \ x\in \partial B_{\mathbb{H}}(0,R)$$
and the differential equation
$$\Delta_{\mathbb{B}^n} u+\frac{(n-1)^2}{4}u+\mu^2qu=\mu^2(q_2-q_1)u_2(x).$$
From this and differential equation for $u_1(x)$, we have
$$\mu^2u_1u_2(q_2-q_1)=u_1\Delta_{\mathbb{B}^n} u-u\Delta_{\mathbb{B}^n} u_1.$$
Now, Green's theorem and the boundary values imply that
$$\int_{B_{\mathbb{H}}(0,R)}u_1u_2(q_1-q_2)dV_{\mathbb{B}^n}(y)=0.$$
\end{proof}
By Theorem \ref{density2}, we derive that $$\int_{B_{\mathbb{H}}(0,R)}v_1(y)v_2(y)(q_1-q_2)dV_{\mathbb{B}^n}(y)=0$$ for
any all solutions $v_1$ and $v_2$ of $$\Delta_{\mathbb{B}^n} v_1+\frac{(n-1)^2}{4}v_1+\mu^2qv_1=0,\ \ \Delta_{\mathbb{B}^n} v_2+\frac{(n-1)^2}{4}v_2+\mu^2qv_2=0.$$

We now prove that $q_1$ must be equal to $q_2$. Let $w_1(x)=\big(\frac{2}{1-|x|^2}\big)^{\frac{n}{2}-1}v_1(x)$ and $w_2(x)=\big(\frac{2}{1-|x|^2}\big)^{\frac{n}{2}-1}v_2(x)$, through
conformal laws:$$\Delta_{\mathbb{B}^n}u+\frac{n(n-2)}{4}u=\big(\frac{2}{1-|x|^2}\big)^{-\frac{n}{2}-1}\Delta_{\mathbb{R}^n}(\big(\frac{2}{1-|x|^2}\big)^{\frac{n}{2}-1}u),$$
we derive that $w_1(x)$ and $w_2(x)$ satisfies equation
\begin{equation}\label{ca1}
\Delta_{\mathbb{R}^n}w_1+\big(\frac{1}{4}+\mu^2q\big)w_1\big(\frac{2}{1-|x|^2}\big)^2=0,\ \ x\in B^n(0,\tanh(\frac{R}{2}))
\end{equation}
and
\begin{equation}\label{ca2}
\Delta_{\mathbb{R}^n}w_2+\big(\frac{1}{4}+\mu^2q\big)w_2\big(\frac{2}{1-|x|^2}\big)^2=0,\ \ x\in B^n(0,\tanh(\frac{R}{2})).
\end{equation}
Combining the above computations, we conclude that $$\int_{B^n(0,\tanh(\frac{R}{2}))}w_1(x)w_2(x)(q_1-q_2)\big(\frac{2}{1-|x|^2}\big)^{2n-2}dx=0$$
for all solutions $w_1$ and $w_2$ of \eqref{ca1} and \eqref{ca2}. In order to prove that $q_1-q_2=0$, we need the following lemma.
\begin{lem}

Let $\Omega\subseteq \mathbb{B}^n$ be a bounded open set, and let $q\in L^{\infty}(\Omega)$. There is a constant $C_0$ depending only on $\Omega$ and $n$ such that for any
$z\in \mathbb{C}^n$ satisfying $z\cdot z=0$ and $|z|\geq \max\{C_0\|q\|_{L^{\infty}},\ 1\}$, and for any function $c\in H^2(\Omega)$ satisfying
$$z\cdot \nabla_{\mathbb{R}^n}c=0,\ \ x\in \Omega,$$
then the equation $\Delta_{\mathbb{R}^n}u+qu=0$ in $\Omega$ has a solution
\begin{equation}\label{CGO}
u(x)=e^{iz\cdot x}(c+r),
\end{equation}
where $r\in H^2(\Omega)$ satisfies $$\|r\|_{L^2(\Omega)}\leq \frac{C}{|Rez|}\|\Delta_{\mathbb{R}^n}c+qc\|_{L^2(\Omega)}.$$
\end{lem}
Equation \eqref{ca1} and \eqref{ca2} can be written as \begin{equation}\label{ca3}
\Delta_{\mathbb{R}^n}w_1+(\frac{1}{4}+\mu^2)qw_1\big(\frac{2}{1-|x|^2}\big)^2=0,\ \ x\in B^n(0,\tanh(\frac{R}{2})).
\end{equation}
For any $y\in \mathbb{R}^n$ and $\rho>0$, we can choose vectors $a$, $b\in \mathbb{R}^n$ such that $\{y,a,b\}$ is an orthogonal basis in $\mathbb{R}^n$ with the properties $|a|=1$ and $|b|^2=|y|^2+\rho^2$. Then for $$z_1=y+\rho a+ib,\ \ z_2=y-\rho a-ib,$$
we have that $$z_j\cdot z_j=|Rez(z_j)|^2-|Imz_j|^2+2i Rez_j\cdot Im z_j=|y|^2+\rho^2-b^2=0, \ j=1,2.$$
From the above lemma, we can choose $w_1=e^{iz_1x}(1+r_{z_1}(x))$ and $w_2=e^{iz_2x}(1+r_{z_2}(x))$, then $$\int_{B^n(0,\tanh(\frac{R}{2}))}w_1(x)w_2(x)(q_1-q_2)\big(\frac{2}{1-|x|^2}\big)^{2n-2}dx=0$$ becomes
$$\int_{B^n(0,\tanh(\frac{R}{2}))}(q_1-q_2)\big(\frac{2}{1-|x|^2}\big)^{2n-2}e^{2iy\cdot x}(1+r_{z_1}(x))(1+r_{z_2}(x))dx=0.$$
Let $\rho\rightarrow +\infty$, we derive that
$$\int_{B^n(0,\tanh(\frac{R}{2}))}(q_1-q_2)\big(\frac{2}{1-|x|^2}\big)^{2n-2}e^{2iy\cdot x}dx=0,$$
which implies that $(q_1-q_2)\big(\frac{2}{1-|x|^2}\big)^{2n-2}=0$. Then we accomplish the proof of Theorem \ref{inverse-medium}.

\vskip0.5cm

\noindent\textbf{Acknowledgement:} The authors would like to thank Prof. Qiaohua Yang for helpful discussion on the Green function of Helmholtz operator in hyperbolic space, which inspired us a lot.

\end{document}